\documentclass{IEEEtran}
\usepackage{cite}
\usepackage{amsmath,amssymb,amsfonts}
\usepackage{algorithmic}
\usepackage{graphicx}
\usepackage{algorithm,algorithmic}
\usepackage{hyperref}
\hypersetup{hidelinks=true}
\usepackage{textcomp}
\usepackage{booktabs}
\def\BibTeX{{\rm B\kern-.05em{\sc i\kern-.025em b}\kern-.08em
    T\kern-.1667em\lower.7ex\hbox{E}\kern-.125emX}}
\usepackage{adjustbox}
\usepackage{subfigure}
\usepackage{makecell}
\usepackage{multirow}
\usepackage{mathrsfs}
\usepackage{dsfont}
\newtheorem{thm}{\bf Theorem}[section]
\newtheorem{cor}[thm]{\bf Corollary}
\newtheorem{lem}[thm]{\bf Lemma}
\newtheorem{rem}[thm]{\bf Remark}

\newtheorem{prop}[thm]{\bf Proposition}

\newtheorem{deff}[thm]{\bf Definition}
\newtheorem{eg}[thm]{\bf Example}
\newcommand{\A}{\mathcal{A}}
\renewcommand{\ll}{\mathcal{L}}
\renewcommand{\lll}{\mathcal{L}_{\lambda}}
\newcommand{\X}{\Omega}

\newcommand{\D}{\mathcal{D}}
\renewcommand{\S}{\mathcal{S}}
\newcommand{\I}{\mathcal{I}}
\newcommand{\J}{\mathcal{J}}

\newcommand{\B}{\mathcal{B}}
\newcommand{\G}{\mathcal{G}}

\newcommand{\rr}{\mathcal{R}}
\renewcommand{\S}{\mathcal{S}}

\newcommand{\K}{\mathcal{K}}
\newcommand{\T}{\mathcal{T}}
\renewcommand{\O}{\mathcal{O}}
\renewcommand{\Sigma}{C}
\newcommand{\F}{\mathcal{F}}
\newcommand{\R}{\mathbb{R}}
\newcommand{\N}{\mathbb{N}}
\newcommand{\Z}{\mathbb{Z}}
\newcommand{\C}{\mathbb{C}}

\newcommand{\Kb}{K}

\newcommand{\Lb}{L}
\newcommand{\Xb}{X}

\newcommand{\tomg}{\tilde{\omega}}
\newcommand{\xb}{x}
\newcommand{\vb}{\mathbf{\xi}}
\newcommand{\xeq}{x_{\operatorname{eq}}}
\newcommand{\roa}{\mathcal{D}}

\newcommand{\Zk}{\mathcal{Z}}
\newcommand{\zk}{\mathfrak{z}}

\newcommand{\hrr}{\widehat{\mathcal{R}}}
\newcommand{\eps}{\varepsilon}

\newcommand{\slim}{\mathop{\mathrm{s\text{-}\lim}}}
\newcommand{\ra}{\rightarrow}
\newcommand{\sra}{\rightharpoonup}
\newcommand{\zero}{\mathbf{0}}
\newcommand{\id}{\operatorname{I}}
\newcommand{\dom}{\operatorname{dom}}
\newcommand{\inte}{\operatorname{int}}

\renewcommand{\Re}{\operatorname{Re}}
\newcommand{\set}[1]{\left\{#1\right\}}

\newcommand{\norm}[1]{\left\Vert #1 \right\Vert}
\newcommand{\abs}[1]{\left\vert #1 \right\vert}
\newcommand{\taus}{T}
\newcommand{\trans}{{\mathsf{T}}}
\newcommand{\dd}{\mathrm{d}}
\newcommand{\rhop}{\rho_{\operatorname{pse}}}
\newcommand{\sss}{\scriptscriptstyle}

\renewcommand{\dagger}{{\scriptscriptstyle +}}

\newcommand{\BlackBox}{\rule{1.5ex}{1.5ex}}    
\newcommand{\pfbox}{\hfill\BlackBox}    
\newcommand{\Qed}{\hfill$\diamond$}

\newcommand{\ymr}[1]{\textcolor[rgb]{1,0.5,0}{#1}}
\newcommand{\ym}[1]{{\color{black} #1}}
\hypersetup{breaklinks=true}

\newenvironment{pf}{\par\noindent{\bf Proof:\ }}{\hfill\BlackBox\\[2mm]}

\usepackage{hyperref}
 \hypersetup{
    colorlinks=true,
    linkcolor=blue,
    filecolor=magenta,      
    urlcolor=black,
    pdftitle={Overleaf Example},
    pdfpagemode=FullScreen,
    }
    \usepackage[misc]{ifsym}

\begin{document}
\title{Resolvent-Type Data-Driven Learning of %Infinitesimal 
Generators for Unknown Continuous-Time Dynamical Systems}
\author{Yiming Meng$^\star$,  Ruikun Zhou$^\star$, Melkior Ornik, \IEEEmembership{Senior Member, IEEE}, and Jun Liu, \IEEEmembership{Senior Member, IEEE}
\thanks{This research was supported in part by NASA under grant numbers
80NSSC21K1030 and 80NSSC22M0070, by the Air Force Office of Scientific
Research under grant number FA9550-23-1-0131, and in part by
the Natural Sciences and Engineering Research Council of Canada and the
Canada Research Chairs program.}
\thanks{$^\star$Equal contribution.}
\thanks{Yiming Meng is with  the
Coordinated Science Laboratory, University of Illinois Urbana-Champaign,
Urbana, IL 61801, USA.
        {\tt\small ymmeng@illinois.edu}.
}
\thanks{Ruikun Zhou is with the Department of Applied Mathematics,
University of Waterloo, Waterloo ON N2L 3G1, Canada {\tt\small ruikun.zhou@uwaterloo.ca}.}
\thanks{Melkior Ornik is with the Department of Aerospace Engineering and the
Coordinated Science Laboratory, University of Illinois Urbana-Champaign,
Urbana, IL 61801, USA.
        {\tt\small mornik@illinois.edu}.}
\thanks{
Jun Liu is with the Department of Applied Mathematics,
University of Waterloo, Waterloo ON N2L 3G1, Canada {\tt\small j.liu@uwaterloo.ca}.}
}

\maketitle

\begin{abstract}
 A semigroup characterization, or equivalently, a characterization by the %infinitesimal 
  generator, is a classical technique used to describe continuous-time nonlinear dynamical systems. In the realm of data-driven learning for an unknown nonlinear system, one must estimate the generator of the semigroup of the system's transfer operators (also known as  the semigroup of Koopman operators) based on discrete-time observations and verify convergence to the true generator in an appropriate sense. As the generator encodes essential  \ym{instantaneous}   transitional information of the system, challenges arise for some existing methods that rely on accurately estimating the time derivatives of the state with constraints on the observation rate. Recent literature develops a technique that avoids the use of time derivatives by employing the logarithm of a Koopman operator. However, the validity of this method has been demonstrated only within a restrictive function space and requires knowledge of the operator's spectral properties. 
In this paper, we propose a resolvent-type method for learning the system generator to relax the requirements on the observation frequency and overcome the constraints of taking operator logarithms. %We formally verify the convergence rate and propose data-driven algorithms embedded within the Koopman operator learning framework. This approach is compatible with other system verification tools using the same set of training data. 
We also provide numerical examples to demonstrate its effectiveness in applications of %such as 
system identification and constructing Lyapunov functions. %, and finding ergodicity measures.

%Notice of Previous Publication. This manuscript substantially improves the work of [1] in the following aspects: 1) All key convergence proofs are completed, providing step-by-step reasoning on the development of the proposed method; 2) The key step of the learning algorithm is provided, explicitly demonstrating the dependence on tuning parameters; 3) A novel modification is derived to relax the constraint on the observation rate; 4) Numerical simulations on representative systems are thoroughly analyzed.
\end{abstract}

\begin{IEEEkeywords}
Unknown nonlinear systems,   infinitesimal generator,  operator-logarithm-free,   observation sampling rate, convergence analysis.
\end{IEEEkeywords}

\section{Introduction}
\label{sec:introduction}
\IEEEPARstart{O}{perator} 
 semigroups are essential for characterizing Markov processes, whether random or deterministic. In the scope of %discussing 
 \ym{the approximation of a Markov process}, the convergence of an approximating sequence of processes should be studied. The technique involves studying the convergence of the processes' (infinitesimal) generators in an appropriate sense. Such convergence implies the convergence of the processes' transition semigroups, which, in turn, implies the convergence of the processes to the true Markov process \cite{ethier2009markov}.
 
 %of Markov processes, the convergence of their (infinitesimal) generators in an appropriate sense should be studied. Such convergence implies the convergence of the corresponding transition semigroups, which  in turn  implies the convergence of the processes \cite{ethier2009markov}.

For a continuous-time nonlinear dynamical system, the operator semigroup, denoted by $\set{\K_t}_{ t\geq 0}$, is known as the semigroup of transfer operators, which are also called Koopman operators \cite{koopman1931hamiltonian, haseli2025recursive}. The corresponding infinitesimal generator, $\ll$, is a differential operator (most likely unbounded \cite[Chapter VIII]{reed1980methods}) that provides essential \ym{instantaneous} transitional information for the dynamical system \cite{pavliotis2008multiscale}. 
Particularly, many mass transport-related PDEs, such as transport equations \cite{evans2022partial}, Lyapunov equations \cite{vannelli1985maximal}, and Hamilton-Jacobi equations \cite{mitchell2005time, bardi1997optimal, crandall1983viscosity}, are defined by the generator. Their solutions serve as valuable tools for discovering physical laws \cite{sjoberg1995nonlinear, haber1990structure} and verifying dynamical system properties related to stability \cite{lin1996smooth, teel2000smooth, meng2022smooth}, safety \cite{ames2016control, ames2019control, meng2023lyapunov}, and kinetic and potential energy \cite{pavliotis2008multiscale, evans2022partial}, with applications across various fields including mechanical, space, and power systems, as well as other physical sciences.

To learn unknown continuous-time dynamical systems, one must approximate the entire semigroup, or equivalently its generator, from discrete trajectory data and demonstrate convergence to the true system in an appropriate sense. Considering possible nonlinear effects, challenges arise in retrieving transient system transitions from unavoidable discrete-time observations. Particularly in real-world applications, due to the limitations of sensing apparatuses in achieving high-frequency observations, it is difficult to accurately evaluate quantities related to transient transitions, such as those involving time derivatives. 

%To enhance Koopman-based learning methods while aligning with generator learning procedures, our goal  is to propose a generator learning scheme  that relaxes the observation rate and demonstrate its convergence to the true generator. 
To enhance learning accuracy, we propose a Koopman-based generator learning scheme that relaxes the observation rate compared to existing methods and demonstrates its convergence to the true generator. Below, we review some crucial results from the literature that are pertinent to the work presented in this paper.

%%In this paper, we propose a generator learning scheme, named the \textit{resolvent-type method}, that is robust to the observation rate and demonstrate its convergence to the true generator.

\subsection{Related Work}

The equivalence between the vector field and the infinitesimal generator of the system semigroup has been shown in \cite{koopman1931hamiltonian, mezic2005spectral}, indicating an equivalence in terms of conventional system (model) identification and generator (or semigroup) characterization.

For autonomous dynamical systems, direct methods such as  Bayesian approaches \cite{pan2015sparse} and the sparse identification of non-linear dynamics (SINDy) algorithm \cite{brunton2016discovering} have been developed to identify state dynamics with a known structure   \cite{sjoberg1995nonlinear, haber1990structure}, relying on nonlinear parameter estimation \cite{nash1987nonlinear, varah1982spline} and static linear regression techniques. However, direct methods   require that time derivatives of the state can be accurately estimated, an assumption that may not be robustly satisfied due to potential challenges such as low sampling rates, noisy measurements, and short observation periods.

%%However, direct methods have a potential shortcoming in their inability to perform effectively when the prior knowledge of the model structure is incorrect. They tend to adjust the parameters of an incorrect nominal model to fit the data and reduce the mean square error, which can result in an identified model that does not reflect the true physical meaning. For instance, the intelligent transportation community is motivated to identify effective models for non-open-source commercial autonomous vehicle models to achieve more accurate traffic predictions. The parameterized models used for this purpose often rely on simplified physical laws, whereas the true physical laws may be far more complex. Evidence from \cite[Table III (IDM)]{gunter2019modeling}   shows that the parameters fitting the model do not make logical sense—indicating that when the speed is slightly different from the desired speed, the acceleration/deceleration tends to be unnaturally abrupt—despite the small error in parameter fitting.  Direct methods also require that time derivatives of the state can be accurately estimated, an assumption that may not be robustly satisfied due to potential challenges such as low sampling rates, noisy measurements, and short observation periods. 

Indirect characterization methods can leverage the learning structure of Koopman operators \cite{otto2021koopman}. By applying certain transformations to the learned Koopman operator, the image of the generator related to the instantaneous state dynamics can be expressed as a linear combination of dictionary functions. Such a framework is based solely on observed state data snapshots, which enhances the robustness of characterization through linear least squares optimization over an expressive dictionary that includes nonlinear functions.  An additional benefit of using indirect methods is that they can streamline the procedures for generator characterization and solving transport-related PDEs using the same set of dictionary functions, including neural network bases \cite{mauroy2013spectral, mauroy2016global, deka2022koopman, meng2023learning}. In this introduction, we summarize two widely used Koopman-based characterization methods. 

Finite-difference method (FDM) is a numerical technique used to solve differential equations by approximating derivatives with finite differences. In the context of Koopman-based learning, the generator on any observable   function is approximated by the rate of its evolution along the trajectory, using a forward difference method with a small discretized time scale. Similar to direct methods, this approximation scheme also relies on the notion of evolution speed or time derivatives. As anticipated, its precision    heavily depends on the size of the temporal discretization \cite{bramburger2024auxiliary,klus2020data,nejati2021data,wang2022data, otto2024learning}.

Recent studies \cite{mauroy2019koopman, klus2020data,drmavc2021identification, black2023safe} have facilitated another indirect learning of the generator using the logarithm of the learned Koopman operator. This approach can potentially circumvent the need for high observation rates and longer observation periods, as it does not require the estimation of time derivatives. %, where $\ll$ is the generator and $\K_t$ is the Koopman operator at $t$.  
Heuristically, researchers tend to represent the Koopman operator $\K_t$ by an exponential form of its  generator $\ll$ as $\K_t=e^{t\ll}$, leading to the converse representation $\ll = \frac{1}{t}\log(\K_t)$ for any $t>0$. However, problems arise given the following concerns: 1) representing Koopman operators in exponential form requires the boundedness of the generator $\ll$; 2) the operator logarithm is a single-valued mapping only within a specific sector of the spectrum; 3) for general systems that fall short of the aforementioned restrictions, it is unclear how the data-driven approximation of the logarithm of Koopman operators converges to the true generator.

As a complement to the work in \cite{mauroy2019koopman, drmavc2021identification}, recent studies \cite{zeng2024sampling, zeng2023generalized} have rigorously investigated the sufficient and necessary conditions under which the Koopman-logarithm based generator learning method can guarantee learning accuracy. To provide the sampling rate, the theorem %heavily 
relies on the concept of a ``generator-bounded space'', which remains invariant under the Koopman operator, and where the generator is bounded when restricted to it. However, the mentioned concept   is  less likely to be verifiable for unknown systems. %Though the results appear to be sound, the conditions are less likely to be verifiable for unknown systems.

%Besides addressing the aforementioned issues, our goal in this paper is to propose a more robust generator learning scheme. 

%Besides addressing the aforementioned issues, our goal in this paper is to propose a more robust generator learning scheme. 

\subsection{Contributions}

To address the aforementioned issues, we aim to propose an operator logarithm-free generator learning scheme, named the \textit{resolvent-type method} (RTM), which is robust to the choice of the dictionary of observable functions and does not require knowledge of spectral properties of the Koopman operators. A brief discussion on the resolvent of the Koopman generator can be found in~\cite{susuki2021koopman}. However, to the best of the authors' knowledge, the current work is the first to utilize a resolvent-based method to identify the generator and, consequently, the vector fields. In summary, the main contributions are:
\begin{enumerate}
    \item We demonstrate the converse relationship between the Koopman operators and the generator in more general cases where the generator may be unbounded. Specifically, we draw upon the rich literature to propose a finite-dimensional approximation based on Yosida’s approximation for these cases.
    \item We provide the analytical convergence   of this approximation and demonstrate the theoretical feasibility of a data-driven adaptation.
    \item We adapt the Koopman operator learning structure for generator learning based on 1), and demonstrate that a modification is needed to accommodate a low observation rate constraint.
    \item We provide numerical experiments and demonstrate the effectiveness of the proposed approach by comparing it with the FDM and logarithm-based methods.
\end{enumerate}
%This scheme will be compatible with the current advances in \cite{meng2023learning} for Koopman-based construction of maximal Lyapunov functions. It is important to note that the method in \cite{meng2023learning} assumes full knowledge of the equilibrium point and acknowledges that verification of the constructed Lyapunov function depends on information about the actual system transitions, which may diminish its predictive value in stability analysis. 

\textbf{Notice of Previous Publication}. This manuscript substantially improves the work of \cite{meng2024koopman} in the following aspects: 1) All key convergence proofs are completed, providing step-by-step reasoning on the development of the proposed method; 2) The key step of the learning algorithm is provided, explicitly demonstrating the dependence on tuning parameters; 3) A novel modification is derived to relax the constraint on the observation rate; 4) Numerical simulations on representative systems are thoroughly studied.

The rest of the paper is organized as follows. In Section~\ref{sec: pre}, we present some preliminaries on dynamical systems and the corresponding semigroup. % and revisit the basic concept of operator convergence. 
In Section \ref{sec: characterize}, we provide a Yosida-type approximation for the infinitesimal generator in the strong operator topology and \ym{justify its robustness, particularly in terms of the class of dictionary functions, as the foundation for subsequent data-driven approximations.} %explicitly show the analytical convergence rate, which will foster a better understanding for the development of data-driven approximations. 
In Section \ref{sec: finite-dim}, we verify the feasibility of using a data-driven approximation by providing a finite-horizon, finite-dimensional reduction. % based on the analysis in Section \ref{sec: characterize}. 
\ym{A modification of this approximation is discussed in Section \ref{sec: quad} to improve precision under potential sampling rate and horizon constraints.  The data-driven algorithm is then detailed in Section \ref{sec: alg}}, followed by numerical simulations in Section \ref{sec: num} that demonstrate the performance of the proposed method.   % and conclude the paper in Section \ref{sec: conclusion}.

\subsection{Notation}
We denote by $\R^d$ the Euclidean space of dimension $d>1$, and by $\R$ the set of real numbers. We use $\abs{\;\cdot\;}$ to denote the Euclidean norm. %For any $x\in\R^n$, denote $x_i$ by its $i$-th element.  
%For $x\in\R^n$ and $r\ge 0$, we denote the ball of radius $r$ centered at $x$ by $\B(x, r)=\set{y\in\R^n:\,\abs{y-x}\le r}$, where $\abs{\;\cdot\;}$ is the Euclidean norm. %For a closed set $A\subset\R^n$ and $x\in\R^n$, we denote the distance from $x$ to $A$ by $\abs{x}_{A}=\inf_{y\in A}\abs{x-y}$ and $r$-neighborhood of $A$ by $\B(A, r)=\cup_{x\in A}\B(x, r)=\set{x\in\R^n:\,\abs{x}_A\le r}$. 
For a set $A\subseteq\R^d$, $\overline{A}$ denotes its closure, 
$\inte(A)$ denotes its interior  and $\partial A$ denotes its boundary.  %%For two sets $A,B\subseteq\R^n$, the set difference is defined by $A\setminus B=\set{x:\,x\in A,\,x\not\in B}$. 
%We use the Frobenius norm  $\|\cdot\|_F$ as the metric for finite-dimensional matrices. 
We also use $\cdot^\trans$ and $\cdot^+$ to denote the matrix/vector transpose and pseudoinverse, respectively. %%Given $a,b\in\R$, we define $a\wedge b:=\min(a,b)$. %denote by $M_i$ its $i^{\text{th}}$ row and by $M_{ij}$ its entry at $i^{\text{th}}$ row and $j^{\text{th}}$ column. Given a matrix $M$, 

For any complete normed function spaces \ym{$(\mathcal{V}, \|\cdot\|_{\scriptscriptstyle \mathcal{V}})$ and $(\mathcal{W}, \|\cdot\|_{\scriptscriptstyle \mathcal{W}})$} of the real-valued observable functions, and for any bounded linear operator \ym{$\B: \mathcal{V}\ra \mathcal{W}$}, we define the operator (uniform) norm as \ym{$\|\B\|:=\sup_{\scriptscriptstyle\|h\|_{\scriptscriptstyle\mathcal{V}}=1}\|\B h\|_{\scriptscriptstyle \mathcal{W}}$}. 
We denote by $\id$ the identity operator on any properly defined spaces. Let $\mathcal{C}(\Omega)$ be the set of continuous functions with domain $\Omega$, \ym{endowed with the uniform   norm $\|\cdot\|:=\sup_{x\in\X}|\cdot(x)|$}. We denote the set of $i$-times continuously
differentiable functions by $\mathcal{C}^i(\Omega)$. \ym{For convergence analysis, we also write $a\lesssim b$ if there exists a $C>0$ (independent of $a$ and $b$) such that $a\leq Cb$, particularly when 
$C$ is lengthy or its exact value is not necessarily required.}

For the reader's convenience, we also provide a partial list of notations related to the technical details, which will be explained later in the article.

\begin{itemize}
    \item[-] $M$: the number of sampled initial conditions;
    \item[-] $N$: the number of dictionary observable test functions;
    %\item[-] $\Zk_N(x)$: the dictionary  of  $N$ observable test functions;
    \item[-] $\K_t, \ll$: the Koopman operator and infinitesimal generator;
    \item[-] $K, L$: the matrix version   of a $\K_t$ (with a fixed $t$) and the $\ll$ after data fitting, with dimension $N\times N$;
        \item[-] $\taus$: terminal observation instance for data collection;
        \item[-] $T_s$:  running time for performance verification;
    \item[-] $\gamma$: observation/sampling frequency;
    \item[-] $\Gamma$: number of snapshots  for each trajectory during data collection, which is equal to $\gamma\taus$;
 
     \item[-] $\tau$: the sampling period, %initial snapshot moment,
     which is equal to $1/\gamma$. 

\end{itemize}
From Section \ref{sec: finite-dim} to \ref{sec: num}, we will frequently use notations that represent approximations of certain quantities. %Heuristically,
We use   $Q_{\text{sub}}^{\text{sup}}$ to denote an approximation of quantity  $Q$ that depends specifically on parameters appearing in the subscripts and superscripts, and use    $\widetilde{Q}$ to denote 
  an approximation of quantity $Q$ without emphasizing its dependent parameters. 
%%We denote the set of locally and uniformly Lipschitz continuous functions by $\loclip(\Omega)$ and $\lip(\Omega)$. 

%When making general statements for $v\in C^1(\R^n)$ with $n\geq 1$, we denote $D$ as its gradient operator (or as the derivative when $n=1$). 

\section{Preliminaries}\label{sec: pre}
\subsection{Dynamical Systems}
Given a pre-compact state space $\X\subseteq\R^d$, 
we consider a continuous-time nonlinear dynamical system of the form 

\begin{equation}\label{E: sys}
   \dot{x} = f(x), %\dot{\xb}(t) = f(\xb(t)),\;\;\xb(0)=x\in\X,\;t\in[0,\infty), 
\end{equation}
where %$x$ denotes the initial condition, and 
the vector field $f:\X\ra\X$ is assumed to be locally Lipschitz continuous.

Given an initial condition $x_0$, on the maximal  interval of existence $I\subseteq [0,\infty)$, the unique forward flow map (solution map) $\phi: I\times \X\rightarrow \X$ should satisfy 1) $\partial_t(\phi(t,x_0)) = f(\phi(t,x_0))$, 2) $\phi(0, x_0)=x_0$, and 3) $\phi(s, \phi(t,x_0))=\phi(t+s, x_0)$ for all $t,s\in I$ such that $t+s\in I$.
\iffalse
\begin{equation}
    \begin{cases}
        &\partial_t(\phi(t,x_0)) = f(\phi(t,x_0)), \\
        & \phi(0, x_0)=x_0, \\
        & \phi(s, \phi(t,x_0))=\phi(t+s, x_0), \;\forall t,s\in\mathcal{I}
    \end{cases}
\end{equation}\fi

%Without loss of generality, 
Throughout the paper, we will assume that the 
maximal interval of existence of the   flow map to the initial value problem   \eqref{E: sys} is $I=[0,\infty)$. %For simplicity of notation, we may also write the solution as $\phi(t)$ if the initial condition is not emphasized. 

\begin{rem}
  The above assumption is equivalent to assuming that the system exhibits forward invariance w.r.t. the set $\X$. However, this is usually not the case for general nonlinear systems. In this paper, if the system dynamics violate the above assumption, we can adopt the approach outlined in \cite[Section III.B]{meng2023learning} to recast  the dynamics within the set $\overline{\X}$. In other words, we constrain the vector field $f$ such that $f(x)=0$ for any $x\in\partial\X$, while $f$ remains unchanged within the open domain $\X$. 
  This modification ensures that the system data is always collectible within $\overline{\X}$. \Qed %, which will be beneficial for learning the system transitions in this region of interest. \Qed
\end{rem}

\subsection{Koopman Operators and the  Infinitesimal Generator}

        Let us now consider a complete normed  function space $(\F, \|\cdot\|_{\scriptscriptstyle\F})$ of the  real-valued observable functions\footnote{We clarify that in the context of operator learning, the term ``observable functions,'' or simply ``observables,'' commonly refers to ``test functions'' for operators, rather than to the concept of ``observability''  in control systems.} $h: \X\rightarrow\mathbb{R}$. %%For any bounded linear operator $B:\F\ra\F$, we define the operator (uniform) norm as $\|B\|:=\sup_{\|h\|_{\F}=1}\|Bh\|_\F$. 

\begin{deff}\label{def: semigroup}
A one-parameter family $\{\S_t\}_{t\geq 0}$ 
of bounded linear operators from $\F$ into $\F$ is a semigroup % ofbounded linear operators on $\F$ 
if
\begin{enumerate}
    \item $\S_0=\id$; 
    \item $\S_t\circ \S_s=\S_{t+s}$ for every $t,s\geq 0$. 
\end{enumerate}
In addition, a semigroup $\{\S_t\}_{t\geq 0}$ is a  $\mathcal{C}_0$-semigroup if 
$\lim_{t\downarrow 0} \S_t h = h$ for all $h \in \F$. %% and moreover  a  $\mathcal{C}_0$-semigroup of contractions if $\norm{\S_t}\leq 1$ for $t\geq 0$. 

The (infinitesimal) generator $\A$ of $\{\S_t\}_{t\geq 0}$ is defined by
\begin{equation}\label{E: generator}
    \A h(x):= \lim_{t\downarrow 0}\frac{\S_t h(x)-h(x)}{t},
\end{equation}
where the observables are within the natural domain of $\A$, defined as $\dom(\A)=\set{h\in\F:\lim_{t\downarrow 0}\frac{\S_t h-h}{t}\;\text{exists}}. $
\Qed
\end{deff}

It is a well-known result that \ym{$\dom(\A)$ is dense in $\F$, i.e.}, $\overline{\dom(\A)} = \F$.

The evolution of observables of system \eqref{E: sys} restricted to $\mathcal{F}$ is governed by the family of Koopman operators, as defined below. Koopman operators also form a linear $\mathcal{C}_0$-semigroup, allowing us to study nonlinear dynamics through the infinite-dimensional lifted space of observable functions, which exhibit linear dynamics.

\begin{deff}
    \label{def: Koopman} % (semi)-group of
The Koopman operator family $\{\K_t\}_{t\geq 0}$ of system \eqref{E: sys} is a collection of maps $\mathcal{K}_t: \mathcal{F} \rightarrow \mathcal{F}$  defined by
\begin{align}
\mathcal{K}_t h = h\circ \phi(t, \cdot), \quad h \in \mathcal{F}
\end{align}
for each $t\geq 0$, where $\circ$ is the composition operator. The  generator  $\ll$ of $\{\K_t\}_{t\geq 0}$ is defined accordingly as in \eqref{E: generator}. \Qed
\end{deff}

Due to the (local) Lipschitz continuity   of $f$ in \eqref{E: sys}, and considering that  observable   functions are usually continuous,  we will focus on $\K_t: \mathcal{C}(\X)\ra \mathcal{C}(\X)$ for the rest of the paper. \ym{For  \eqref{E: sys}, there exist  constants $\omega\geq 0$ and $\Sigma\geq 1$ such that $\norm{\K_t}\leq C
   e^{\omega t}$ for all $t\geq 0$ \cite[Theorem 1.2.2]{pazy2012semigroups}. }

%%See also the fundamental properties of $\{\K_t\}$ in Proposition \ref{prop: fact1} and Remark \ref{rem: smooth_version} in Appendix \ref{sec: facts}. 

\ym{It can also be verified that $\mathcal{C}^1(\Omega)\subseteq\dom(\ll)\subseteq \mathcal{C}(\Omega)$. In general, $\dom(\ll)$ depends on the regularity and degeneracy   of 
$f$, as well as the geometry of the flow, which determine whether 
$f\cdot\nabla h$ ($h\in\dom(\ll)$) exists in a classical, weak, or viscosity sense. For instance, 
$h\in\dom(\ll)$ need not be differentiable on $E:=\{x\in\R^d: f(x)=0\}$
 but it may be $\mathcal{C}^1$ on 
$\Omega\setminus E$ with 
$f\cdot\nabla h$ extending continuously to $A$. In the extreme case 
$f(x)\equiv0$, we have  
$\dom(\ll)=\mathcal{C}(\Omega)$.}
 
 \ym{For an unknown 
$f$, it is impossible to recover 
$\ll$ on the entire domain $\dom(\ll)$. We therefore restrict our attention to a \textit{core} of $\dom(\ll)$, based on the following concept.}
\begin{deff}\label{def: core}
\ym{Consider $\A:\dom(\A)\subseteq\F\to\F$ and $\B:\mathcal{D}\subseteq\F\to\F$. 
We say that $\B$   admits a closure  w.r.t. $\A$, denoted as $\overline{\B} =\A$,  if  $\A|_{\mathcal{D}} = \B$  and $\overline{\mathcal{D}}=\dom(\A)$ w.r.t. the graph norm of $\A$.  
In this case, we also say that $\mathcal{D}$ is a core of $\A$. \Qed}
\end{deff}

\ym{For $\F = \mathcal{C}(\Omega)$  endowed with the uniform norm $\|\cdot\|:= \sup_{x\in\X}|\cdot(x)|$, the graph norm of $\ll$ is naturally defined as  
\begin{equation}\label{E: norm_derivative}
    \|\cdot\|_{\sss\ll} := \|\cdot\| + \|\ll(\cdot)\|. 
\end{equation}
%Let $\ll:=\ll|_{\scriptscriptstyle\mathcal{C}^1(\Omega)}$.  
Then,  by Proposition \ref{prop: core},  we have $\overline{\ll|_{\scriptscriptstyle\mathcal{C}^1(\Omega)}} = \ll$, and $\mathcal{C}^1(\Omega)$ is the {core} of $\ll$. For the remainder of the paper, we work with $\mathcal{C}^1(\Omega)$-class dictionary of observables, where approximation on this core is effectively equivalent to approximation on the full domain for the purpose of learning  $\ll$ with both approximation guarantees and numerical tractability.} Moreover,  for all $h\in \mathcal{C}^1(\X)$ \cite{pavliotis2008multiscale}, the generator of the Koopman operators is  $\ll h(x) = \nabla h(x) \cdot f(x)$.

 \ym{\begin{deff}\label{def: dictionary}
 For operator learning, we define
     \begin{equation}\label{E: dict}
    Z_N(x):=[\zk_{\scriptscriptstyle  0}(x), \zk_{\scriptscriptstyle  1}(x),  \cdots, \zk_{\scriptscriptstyle  N-1}(x)]^\trans, \;N\in\N, 
\end{equation}
 as the the vector of dictionary functions,   where $\{\zk_i\}\subseteq\mathcal{C}^1(\Omega)$.    We denote by $\Zk_N:=\operatorname{span}\{\zk_i: i = 1, 2, \cdots, N\}$ the span of dictionary functions. 

 For any   linear operator $\B:\mathcal{C}(\X)\ra\mathcal{C}(\X)$, we adopt the notational conventions $\B(\D) :=\{\B h: h\in\D\}$ for any $\D\subseteq \mathcal{C}(\X)$ and $\B Z_N(x) = [\B\zk_0(x), \B\zk_1(x), \cdots, \B\zk_{\scriptscriptstyle N-1}(x)]^\trans$. \Qed
 \end{deff}
 
 We make the   standing assumption on the dictionary that  $\Zk_1\subseteq\Zk_2\subseteq\cdots$ and 
 $\overline{\cup_{n\geq 1}\Zk_N} = \mathcal{C}(\Omega)$ to ensure density. %for the density requirement of Galerkin-type approximations.  
 Denoting $E_N(h):=\inf_{v\in\Zk_N}\|h-v\|$ as the best-approximation error, then 
 $E_N(h)\xrightarrow{N\ra\infty} 0$ for any $h\in \mathcal{C}(\X)$.}

   %However, we will pursue this direction in future work. 

%%Given that the observable functionsare continuously differentiable, the generator of Koopman operators is such that $\ll h(x) = \nabla h(x) \cdot f(x)$ for all $h\in C^1(\X)$. 

\begin{rem}\label{rem: FDM}
    In view of Koopman-based indirect approximation of the generator, for each $h\in \mathcal{C}^1(\X)$, we have $\K_\tau h(x) -h(x)    =  \int_0^\tau \ll h(\phi(s, x)) ds   \approx  \int_0^\tau \ll h(\phi(0, x)) ds = \ll h(x) \tau$, 
\iffalse
\begin{equation*}
    \begin{split}
        \K_\tau h(x) -h(x)  & =  \int_0^\tau \ll h(\phi(s, x)) ds   \approx  \int_0^\tau \ll h(\phi(0, x)) ds = \ll h(x) \tau,\\
    \end{split}
\end{equation*}\fi
where the approximation is achieved through a small terminal time $\tau$. Then, the derivative (of any $h$) along the trajectory is approximated by $\ll h(x) \approx \frac{\K_\tau h(x) - h(x)}{\tau}.$ Note that the finite-difference method (FDM) $\ll   \approx  \frac{\K_\tau - \id}{\tau}$ %for some $\tau>0$
\iffalse
\begin{equation}\label{E: fdm}
    \ll   \approx  \frac{\K_\tau - \id}{\tau}, \;\;\tau>0
\end{equation}\fi
simply follows \eqref{E: generator} for sufficiently small $\tau>0$. % without taking the limit $\tau\downarrow 0$. 
Through this approximation scheme of the time derivative, it can be anticipated that the precision heavily depends on the size of $\tau$ \cite{bramburger2024auxiliary,nejati2021data,wang2022data}. We will revisit the FDM in the numerical examples to compare it with the method proposed in this paper. \Qed
\end{rem}

%%There is literature that, although focused on stochastic differential equations \cite{nejati2021data, wang2022data}, relies on the definition \eqref{E: generator} to learn the generator. Despite the theoretical soundness, sampling trajectory information within an arbitrarily small time horizon is not practical. %We will elaborate on how our proposed method could possibly avoid a high sampling rate later in this paper.

\subsection{Representation of Semigroups} % and the Generator}
In this subsection, we introduce basic operator topologies and explore how a semigroup $\set{\S_t}_{t\geq 0}$ can be represented through its generator $\A$.

\begin{deff}[Operator Topologies]
Consider Banach spaces \ym{$(\mathcal{V}, \|\cdot\|_{\scriptscriptstyle \mathcal{V}})$ and $(\mathcal{W}, \|\cdot\|_{\scriptscriptstyle \mathcal{W}})$}  %$(\F, \|\cdot\|_\F)$. 
    Let $\B: \mathcal{V}\ra\mathcal{W}$ and $\B_n:\mathcal{V}\ra\mathcal{W}$, for each $n\in\N$, 
be  linear operators. 
\begin{enumerate}
    \item The sequence $\{\B_n\}_{n\in\N}$ is said to converge to $\B$ \emph{uniformly}, denoted by $\B_n\ra \B$, if $\lim_{n\ra\infty}\norm{\B_n-\B}=0$. We also write $\B=\lim_{n\ra\infty}\B_n$. 
    \item The sequence $\{\B_n\}_{n\in\N}$ is said to converge to $\B$ \emph{strongly}, denoted by $\B_n\sra \B$, if \ym{$\lim_{n\ra\infty}\|\B_nh-\B h\|_{\scriptscriptstyle \mathcal{W}}=0$ for each $h\in\mathcal{V}$}. We also write $\B=\slim_{n\ra\infty}\B_n$. \Qed
\end{enumerate} 
\end{deff}

\begin{rem}
   In analogy to the pointwise convergence of functions,   the strong topology is the coarsest topology such that $\B\mapsto \B h$ is continuous in $\B$ for each fixed $h\in\mathcal{V}$. \Qed
\end{rem}

%Suppose that 
If $\A$ is a
bounded linear operator that generates $\{\S_t\}$, then $\S_t = e^{t\A}$ for each $t$ in the uniform operator topology.  %In the case where $A$ is unbounded, 
Otherwise, \cite[Chap. I, Theorem 5.5]{pazy2012semigroups} (also rephrased as Theorem \ref{thm: H-Y} in this paper) still provides an interpretation for the  sense in which  $\S_t$ ``equals" $e^{t\A}$.

We revisit some facts to show the above concepts of equivalence, particularly in the context where $\A$ is unbounded.

\begin{deff}[Resolvents]
Let $\A:\dom(\A)\subseteq\F\ra\F$ be a linear, not necessarily bounded, operator. Then the resolvent set is defined as 
\begin{small}
    \begin{equation*}
\rho(\A):=\set{\lambda\in\C:  \ym{ \lambda\id-\A  \;\text{is bijective and\;}(\lambda\id-\A)^{-1} \;\text{is bounded}}} . 
    %\rho(\A):=\set{\lambda\in\C: \lambda\id-\A\;\text{is bijective and} \;(\lambda\id-\A)^{-1} \;\text{is continuous} }. 
\end{equation*}
\end{small}

\noindent The resolvent  operator is defined as
\begin{small}
\begin{equation*}
    \rr(\lambda; \A):=(\lambda \id-\A)^{-1}, \;\lambda\in\rho(\A). \tag*{\(\diamond\)}
\end{equation*}
\end{small}
\end{deff}
%\ymr{given that $\rho(\A)\neq\emptyset$}. 
\noindent   We further define the \textit{Yosida approximation} of  $\A$  as %The family ${\rr(\lambda; \A)}_{\lambda \in \rho(\A)}$ consists of bounded linear operators \cite[Chap. I,  Theorem 4.3]{pazy2012semigroups}.  
\begin{small}
    \begin{equation}
   \A_\lambda  :=\lambda \A \rr(\lambda;\A) = \lambda^2\rr(\lambda;\A)-\lambda \id. %,\;\;\lambda>0.
\end{equation}
\end{small}

\noindent Note that ${\rr(\lambda; \A)}_{\lambda \in \rho(\A)}$ and $\set{ \A_\lambda}_{\lambda\in\rho(\A)}$ are families of bounded linear operators \cite[Chap. I,  Theorem 4.3]{pazy2012semigroups}, and $e^{t \A_\lambda}$ is well-defined for each $\lambda\in\rho(\A)$. 

\ym{The characterization of $\A$ as the infinitesimal generator of a $\mathcal{C}_0$-semigroup is typically formulated in terms of conditions on the resolvent of $\A$, where 
$\A$ must be closed and $\rho(\A)\neq\emptyset$. We also have the following theorem for semigroup approximation.}

\begin{thm}\label{thm: H-Y} \cite[Chap. I, Theorem 5.5]{pazy2012semigroups}
    Suppose that $\set{\S_t}_{t\geq 0}$ is a $\mathcal{C}_0$-semigroup on $\F$ and $\A$ is the generator. Then, $\S_t  = \slim_{\lambda\ra\infty} e^{t \A_\lambda}$ for all $t\geq 0$. 
    %$$\S_t  = \slim_{\lambda\ra\infty} e^{tA_\lambda},\;t\geq 0.$$
    \iffalse
    If $\set{\S_t}_{t\geq 0}$ is moreover a $C_0$ semigroup  of contractions, then on $\dom(A)$, 
    \begin{equation}
        A = \slim_{\lambda\ra\infty}A_\lambda.
    \end{equation}\fi
\end{thm}

%In other words, it is not possible to represent the semigroup directly using the exponential form of its unbounded generator. Instead, we resort to an asymptotic approximation, as described in Theorem\ref{thm: H-Y}. 

\section{The Characterization of the Infinitesimal
Generators}\label{sec: characterize}
For  system \eqref{E: sys}, we are able to represent $\K_t$ as $e^{t\ll}$ for bounded $\ll$, or, by Theorem \ref{thm: H-Y},  as the strong limit of $e^{t\ll_\lambda}$, where $\ll_\lambda$ is the Yosida approximation of   $\ll$.  

For the converse representation of $\ll$ based on $\{\K_t\}_{t\geq 0}$,  it is intuitive to take the operator logarithm such that $\ll = (\log \K_t)/t$. When $\ll$ is bounded, its spectrum's sector should be confined to make the logarithm a single-valued mapping \cite{zeng2024sampling, zeng2023generalized}. However, for an unbounded $\ll$, there is no direct connection $\ll$ and $\{\K_t\}_{t\geq 0}$.  In this subsection, we review how $\ll$ can be properly approximated based on $\{\K_t\}_{t\geq 0}$. 

\ym{Recall that $\norm{\K_t}\leq C
   e^{\omega t}$ for some    $\omega\geq 0$ and $\Sigma\geq 1$   for all $t\geq 0$. We first examine the classical Hille-Yosida theorem \cite[Theorem 3.1, Chapter I]{pazy2012semigroups}, where $\ll$ is the strong limit of the Yosida approximation:}
   \iffalse
First, we examine the base case where $\{\K_t\}_{t\geq 0}$ is a contraction semigroup. In this situation, by the famous Hille-Yosida theorem \cite[Theorem 3.1, Chapter I]{pazy2012semigroups}, the resolvent set is such that $\rho(\ll)\supseteq \R^+$, and for each $\lambda\in\rho(\ll)$, we have $\norm{R(\lambda; \ll)}\leq 1/\Re(\lambda)$. Consequently, 
by \cite[Lemma 1.3.3]{pazy2012semigroups}, $\set{\lll}_{\lambda>0}$ converges in a strong sense to $\ll$ (w.r.t. $|\cdot|_\ll$), i.e.,\fi
\begin{small}
    \begin{equation}
    \slim\limits_{\lambda\ra\infty}\lll = \slim\limits_{\lambda\ra\infty}\lambda^2R(\lambda;\ll)-\lambda \id = \ll.
\end{equation}
\end{small}

Although it has been demonstrated extensively in   literature, to motivate further developments and provide insight  into how data-driven approaches can be integrated into the approximation scheme, we present the following theorem. %, which establishes an explicit convergence rate. 
The proof is provided in Appendix \ref{sec: proof_fundamental}.

%%To accommodate more general cases, we propose the following extension of the Yosida approximation on $\R^+$. Although it has been demonstrated in extensive literature that $\slim_{\lambda\ra\infty}\lll = \ll$ for $\set{\lll}_{\lambda>\omega}$, to better understand how data-driven approaches can be integrated into the approximation scheme,  we prove the following theorem and demonstrate an explicit convergence rate.

\begin{thm}\label{thm: conv}
    For system \eqref{E: sys}, consider  $\set{\lll}_{\lambda>\tomg}$, where $\tomg\geq \omega\geq0$. Then, $\slim_{\lambda\ra\infty}\ll_{\lambda}=\ll$. Moreover, for any $h\in \mathcal{C}^2(\X)$, %there exists an $\widetilde{\Sigma}>0$ such that 
    \begin{small}
        \begin{equation*}
            \|\lll h-\ll h\|\lesssim %\widetilde{\Sigma}
            (\lambda-\tomg)^{-1} (\|h\|_{\sss\ll}+\|\ll h\|_{\sss\ll}).
        \end{equation*}
    \end{small}
\end{thm}

%\subsection{Representation of  Resolvent Operators}

Motivated by representing $\ll$ by $\set{\K_t}_{t\geq 0}$ and  the Yosida approximation for $\ll$ on $\{\lambda>\omega\}$, we  establish a connection between $\rr(\lambda;\ll)$ and $\set{\K_t}_{t\geq 0}$. 

\begin{prop}\label{prop: R_form}
    Let $\rr(\lambda)$ on $\mathcal{C}(\X)$ be defined by
    \begin{small}
            \begin{equation}
        \rr(\lambda)h:=\int_0^\infty e^{-\lambda t}(\K_t h) dt. 
    \end{equation}
    \end{small}

\noindent     Then, for all $\lambda>\omega$,  
    \begin{enumerate}
        \item $\rr(\lambda)(\lambda \id-\ll)h = h$ for all $h\in \dom(\ll)$;
        \item $(\lambda \id-\ll)\rr(\lambda)h = h$ for all $h\in \mathcal{C}(\X)$.
    \end{enumerate}
\end{prop}

The proof follows the standard procedures of calculus and dynamic programming, and is completed in  Appendix D-A of the supplementary material. %\ref{sec: proof_fundamental}. 

\ym{The family of $\{\rr(\lambda)\}$ serves as  pseudo-resolvent operators. % and satisfies the resolvent identity shown in Proposition \ref{prop: pseudo}. 
It becomes a true resolvent family, i.e., satisfies the commutative property between  $\rr(\lambda)$ and $\lambda\id-\ll$, only when the inverse mapping of $\rr(\lambda)$ is defined on its range $\dom(\ll)$ rather than on the entire $\mathcal{C}(\X)$.} 

\ym{Note that the pseudo-resolvent need not be a resolvent for the restriction of $\ll$\ to a subdomain, not even for the \emph{pre-generator}  $\ll|_{\scriptscriptstyle\mathcal{C}^1(\Omega)}$. In particular, if $\ll|_{\scriptscriptstyle\mathcal{C}^1(\Omega)}$ is not closed, %on $\mathcal{C}^1(\X)$. 
Proposition \ref{prop: R_form}-2) (asserting surjectivity of 
$\lambda\id-\ll|_{\scriptscriptstyle\mathcal{C}^1(\Omega)}$) may fail; in that case  one may only have $\overline{(\lambda\id-\ll)(\mathcal{C}^1(\X))} = \mathcal{C}(\X)$. We illustrate this fact in Example \ref{eg: empty_resolvent}.} 
    %Even though $\rr(\lambda)$ is well defined on $\mathcal{C}(\X)$, the commutative property between $\rr(\lambda)$ and $\lambda\id-\ll$ only holds on $\mathcal{C}^1(\X)$. Consequently, $R(\lambda;\ll)$ and $R(\lambda)$ are equivalent only on $\mathcal{C}^1(\X)$. This domain aligns with the valid domain where $\ll=\slim_{\lambda\ra\infty}\lll$.

\ym{
\begin{eg}\label{eg: empty_resolvent}
Consider   system $\dot{x} = -x^3$ on $\X = (-a, a)$ for any $a > 0$, whose solution is given by $\phi(t, x) = \frac{x}{\sqrt{1+2tx^2}}$. It can be verified that $\mathcal{C}^1(\X) \subsetneq \dom(\ll)$, with a counterexample given by the continuous function $h(x) = |x|^{1/2}$, which lies in $\dom(\ll)$ but not in $\mathcal{C}^1(\X)$.
   % Consider system $\dot{x} = -x^3$ on  $\X = [-a, a]$ for any $a>0$, whose solution is given by $\phi(t, x) = \frac{x}{\sqrt{1+2tx^2}}$. It can be verified that $\mathcal{C}^1(\X)\subsetneq\dom(\ll^\dagger)$, as one can verify that the continuous function $h(x)=|x|^{1/2}$ is within $\dom(\ll^\dagger)$ but $h\notin\mathcal{C}^1(\X)$. 
   Indeed, $\ll h(x) %= \lim_{t\ra 0}\frac{h(\phi(t, x))-h(x)}{t} 
   = \lim_{t\downarrow 0}\sqrt{|x|}\cdot\frac{(1+2tx^2)^{-1/4}-1}{t} = -\frac{|x|^{5/2}}{2}$ for all $x\in\X$.  
   
   Now, let $g: = (\lambda\id-\ll)h \in\mathcal{C}(\X)$. Then, by Proposition \ref{prop: R_form}-1), we have $\rr(\lambda)g = h\notin\mathcal{C}^1(\X)$. Hence $(\lambda\id-\ll|_{\scriptscriptstyle\mathcal{C}^1(\Omega)}) \rr(\lambda)g$ is not well defined, and Proposition \ref{prop: R_form}-2) does not hold for $\ll|_{\scriptscriptstyle\mathcal{C}^1(\Omega)}$. This fact implies that $\ll|_{\scriptscriptstyle\mathcal{C}^1(\Omega)}$ is not closed on $\mathcal{C}^1(\X)$, and one cannot characterize $(\lambda I - \ll|_{\scriptscriptstyle\mathcal{C}^1(\Omega)})^{-1}$ via $\mathcal{R}(\lambda)$.
 Accordingly, $\rho(\ll|_{\scriptscriptstyle\mathcal{C}^1(\Omega)})=\emptyset$; since if it is not, $(\lambda\id-\ll|_{\scriptscriptstyle\mathcal{C}^1(\Omega)})^{-1}$ is continuous for some $\lambda$, which implies that $\ll|_{\scriptscriptstyle\mathcal{C}^1(\Omega)}$ is closed and contradicts the discussion above.  \Qed %
\end{eg}}

\begin{rem}\label{rem: inversibility}
   \ym{Although  
$\rr(\lambda) = (\lambda\id-\ll)^{-1}$ from Proposition \ref{prop: R_form}, we do not attempt to learn $\ll$ directly from this inversion relationship  with $\rr(\lambda)$ in this paper.  This contrasts with \cite{kostic2024learning} in the following ways. The approach in \cite{kostic2024learning} considers stochastic systems and assumes the resolvent  is compact, implying a discrete spectrum. It also employs a Dirichlet-form energy-based metric, which implicitly requires $\ll$ to be self-adjoint. Under these assumptions, 
$\ll$ admits an exact spectral decomposition sharing eigenfunctions $\{\varphi_i\}$ with $\rr(\lambda)$. Based on this fact,  if $\{\beta_i\}$ are the learned eigenvalues of $\rr(\lambda)$,   the corresponding eigenvalues $\{\alpha_i\}$ of $\ll$ satisfy $(\lambda-\alpha_i)^{-1}\varphi_i\approx\beta_i\varphi_i$. 

While such assumptions are reasonable for stochastic systems (owing to the diffusion term), they are generally not met for deterministic ODE systems (see also Section~\ref{sec: sub_compact}). At best, one can approximate $\ll$ \emph{strongly} by a spectral decomposition restricted to a subdomain $\mathcal{D}\subsetneq \dom(\ll)$. However, without extra information about $f$, this restriction $\ll|_{\scriptscriptstyle\mathcal{D}}$ (and hence its  approximation $\widetilde{\ll|_{\scriptscriptstyle\mathcal{D}}}$) may fail to be closed \cite{zeng2024sampling, liu2025properties}. In the same spirit as Example~\ref{eg: empty_resolvent}, it is dangerous to assume that the resolvent of $\ll|_{\scriptscriptstyle\mathcal{D}}$ exists and to represent $(\lambda\id-\widetilde{\ll|_{\scriptscriptstyle\mathcal{D}}})^{-1}$ by a spurious eigen-decomposition, not to mention forcing the data to match a supposed relation between the pseudo-resolvent $\rr(\lambda)$ and that spurious decomposition of $(\lambda\id-\widetilde{\ll|_{\scriptscriptstyle\mathcal{D}}})^{-1}$ for an eigenvalue approximation as \cite{kostic2024learning}. \Qed
}
\end{rem}

\ym{In light of Remark~\ref{rem: inversibility}, where a strong-sense approximation of $\ll$ may not be feasible from the information in $\rr(\lambda)$ via  direct inversion, we do not pursue this inverse strategy and instead adopt the Yosida-type approximation, whose existence and boundedness are always guaranteed.} To use the approximation, we replace $R(\lambda;\ll)$ with $R(\lambda)$. We can then immediately conclude the following representation.

\begin{cor}
    For each 
    $\lambda>\omega$, 
    \begin{small}
            \begin{equation}\label{E: app_L}
        \lll =\lambda^2\int_0^\infty  e^{-\lambda t}\K_tdt -\lambda\id
    \end{equation}
    \end{small}

\noindent and $\lll\sra\ll$ on  \ym{$\dom(\ll)$ as $\lambda\ra\infty$}. 
\end{cor}

\ym{The remainder of the paper builds on this Yosida-type approximation \eqref{E: app_L} and  the evaluation of 
$\rr(\lambda)$ using Koopman-related information in the integral to learn 
$\ll_\lambda$ from data.}

\section{Koopman-Based Finite-Dimensional Approximation of Resolvent Operators}\label{sec: finite-dim}
%To employ \eqref{E: app_L}, the current form of $\int_0^\infty  e^{-\lambda t}(\K_t\cdot)dt$ is not advantageous. 

\ym{The rest of the paper focuses on learning a bounded operator $\rr(\lambda)$ — and hence $\ll_\lambda$— using the dictionary of observable functions $Z_N$ introduced in Definition \ref{def: dictionary}. 
 This process is conceptually similar to the conventional Koopman learning framework, and also resembles many recent works on learning Koopman-related bounded operators, only if input data from the domain and output data from the range are accessible. 
 
 To make the learning of $\ll_\lambda$ effective, it is necessary to preprocess the integral representation of 
$\rr(\lambda)$ in \eqref{E: app_L}. In this section, we show that $\{\rr(\lambda)\}$ can, in principle, be approximated by  finite-rank operators, and the resulting form serves as an ``interface” for generating output data. }

\subsection{Finite Time-Horizon Approximation of the Resolvent}
Observing the form of \eqref{E: app_L}, we first define the following truncation integral operator for $\rr(\lambda)$. 
\begin{deff}\label{def: T_t}
    For any $h\in \mathcal{C}(\X)$ and $\taus\geq0$, we define $\rr_{\lambda,\taus} : \mathcal{C}(\X)\ra \mathcal{C}(\X)$ as
\begin{equation}
     \rr_{\lambda,\taus}  h(x):=\int_0^{T}e^{-\lambda s}\K_{s} h(x) ds.  
\end{equation} 
\end{deff}

We aim to demonstrate that for any sufficiently large $\lambda\in\R$, the aforementioned truncation of the integral will not significantly ``hurt'' the accuracy of the approximation \eqref{E: app_L}. 

\begin{thm}\label{thm: conv_t}
    Let $\taus\geq 0$ and $\lambda>\omega$ be fixed. Then, $\ym{E_{\operatorname{tr}}(\lambda,T)}:=\|\lambda^2 \rr_{\lambda,\taus}  - \lambda\id  -\lll \|\leq  \frac{\Sigma\lambda^2}{\lambda-\omega}e^{-\lambda \taus}$ on \ym{$\dom(\ll)$}. 
\end{thm}
\begin{pf}
\iffalse
    It is clear that $\T_0h(x)=h(x)$. We also have the following identities. 
            \begin{equation*}
        \begin{split}
            &\T_s\circ\T_th(x) \\
            =  & \int_0^s e^{-\lambda r}\K_r\int_0^t e^{-\lambda \sigma}h(\phi(\sigma, x))d\sigma dr\\
            =  & \int_0^s e^{-\lambda r}\int_0^t e^{-\lambda \sigma}h(\phi(\sigma, \phi(r,x)))d\sigma dr\\
            =  & \int_0^s e^{-\lambda r}\int_0^t e^{-\lambda \sigma}h(\phi(\sigma+r, x))d\sigma dr\\
            =  & \int_0^s \int_0^t e^{-\lambda (\sigma+r)}h(\phi(\sigma+r, x))d\sigma dr\\
        \end{split}
    \end{equation*}\fi
Note that, for any $\lambda>\omega$, 
\begin{small}
    \begin{equation*}
\begin{split}
        \|\rr(\lambda)\|\   \leq \int_0^\infty e^{-\lambda t} \|\K_t\| dt  \leq \int_0^\infty \Sigma e^{-\lambda t}e^{\omega t}  dt = \frac{\Sigma}{\lambda-\omega}.
\end{split}
\end{equation*}
\end{small}

\noindent Therefore, for any $h\in C(\X)$, 
\begin{small}
        \begin{equation}
        \begin{split}
           & \|\rr_{\lambda,\taus}  h-\rr(\lambda;\ll)h\|   = \sup_{x\in\X}|e^{-\lambda \taus} \rr(\lambda)h(\phi(\taus,x))|\\
            \leq & e^{-\lambda \taus}\|\rr(\lambda)\|\|h\|
             \leq   e^{-\lambda \taus}\frac{\Sigma}{\lambda -\omega}\|h\|
        \end{split}
    \end{equation}
\end{small}

\noindent and $\sup_{\|h\| =1}\|\lambda^2 \rr_{\lambda,\taus} h - \lambda h -\lll h\| \leq  \frac{\Sigma\lambda^2}{\lambda-\omega}e^{-\lambda \taus}$, 
which completes the proof. 
\end{pf}

% We notice that $\lim_{\lambda\ra\infty}\frac{\lambda^2}{\lambda-\omega} e^{-\lambda \taus} \leq \lim_{\lambda\ra\infty}\frac{2}{\taus(\lambda-\omega) e^{\lambda \taus}}$
\iffalse
\begin{equation}
    \begin{split}
        \lim_{\lambda\ra\infty}\frac{\lambda^2}{\lambda-\omega} e^{-\lambda \tau} = \lim_{\lambda\ra\infty}\frac{\lambda^2}{(\lambda-\omega)e^{\lambda \tau}}\leq \lim_{\lambda\ra\infty}\frac{2}{\tau(\lambda-\omega) e^{\lambda \tau}}
    \end{split}
\end{equation}\fi
The truncation error \ym{$E_{\operatorname{tr}}$} %$\frac{\Sigma\lambda^2}{\lambda-\omega}e^{-\lambda \taus}$ 
in Theorem \ref{thm: conv_t} demonstrates an exponential decaying rate %in the uniform sense 
as $\lambda\rightarrow\infty$ for any fixed $\taus>0$. \ym{We also attempt to use a relatively small 
$T$ to reduce the data size for the integral evaluation, as seen in   Section \ref{sec: quad}.}

\subsection{Finite-Rank Approximation of the Resolvent}\label{sec: sub_compact}
%%Based on \eqref{E: approx_t}, it suffices to learn the operator $\rr_{\lambda,\taus}$ for any fixed $\taus>0$ and to  predict the image of the operator when acting on some $\mathcal{C}(\X)$ function.  
In favor of a learning-based approach based on a dictionary of a finite number of observable test functions serving as basis functions, we verify if
$\rr_{\lambda,\taus}$ is representable as a finite-rank operator.  The proofs are completed in Appendix \ref{sec: proof_characterize}. 

We first look at the following property of $\rr_{\lambda,\taus}$. 

\begin{prop}\label{prop: compact_resolvent}
  For any $\lambda>\omega$, the operator $\rr_{\lambda,\taus}$ is compact if and only if $\K_s$ is compact for any $s\in(0, \taus]$. 
\end{prop}

It is worth noting that $\K_t$ of \eqref{E: sys} is not necessarily compact for each $t>0$. To show that $\K_t(B_r) \subseteq \mathcal{C}(\X)$ is relatively compact, where $B_r = \{h \in \mathcal{C}(\X): \|h\| \leq r\}$ for some $r>0$, one needs to verify the equicontinuity within $\K_t(B_r)$. However, this is not guaranteed. As a counterexample, we set $\X:=(-1, 1)$, $h_n(x) = \sin(nx) \in B_1$ (or similarly, for the Fourier basis), and define $\phi(t, x)= e^{-t}x$ for all $x\in\X$. Then, the sequence ${h_n\circ \phi(t,\cdot)}$ for each $t$ does not exhibit equicontinuity due to the rapid oscillation as $n$ increases.

For the purpose of approximating $\ll_\lambda$ strongly, we aim to find a compact approximation of $\set{\K_t}_{t>0}$ that enables a finite-rank representation of $\rr_{\lambda,\taus}$ in the same sense.

\begin{prop}\label{prop: smoothen}
   \ym{Consider a smooth mollifier $\eta\in C_0^\infty(\X)$ with compact support $\mathbb{B}(0;1)$. For each $\eps>0$, let 
$\eta_\eps(x):= \frac{1}{\eps^n}\eta\left(\frac{x}{\eps}\right)$ with $\int_\X \eta(y)dy = 1$. For all $h\in\mathcal{C}(\X)$ and $\eps>0$, set $\J_\eps h(x) = \int_\X \eta_\eps(x-y)h(y)dy$. Define $\K_t^\eps : = \J_\eps\K_t$.} Then,   for each $t>0$, $\{\K_t^\eps\}_{\eps>0}$ is a family of compact linear operator and satisfies $\K_t^\eps \sra \K_t$ on $\mathcal{C}(\X)$ as $\eps\ra 0$.  %there exists a family of compact linear operator $\{\K_t^\eps\}_{\eps>0}$, such that $\K_t^\eps \sra \K_t$ on $\mathcal{C}(\X)$ as $\eps\ra 0$. 
In addition, for each $t, s>0$, there exists a family of compact linear operator $\{\K_t^\eps\}_{\eps>0}$, such that $\K_t^\eps \circ\K_s^\eps \sra \K_t \circ \K_s$ on  $h\in \mathcal{C}(\X)$ as $\eps\ra 0$. 
\end{prop}

    We omit the proof  as it  follows similar arguments presented in \cite[Section V.A]{meng2023learning}. Heuristically, \ym{the bump functions $\eta_\eps$ converge weakly to Dirac measures centered at their respective flow locations,   distributing point masses that $\{\K_t\}$
 transports along the trajectories.}
%Heuristically, the semigroup $\{\K_t\}_{t>0}$ characterizes the flow of point masses of the system  \eqref{E: sys}, where each point mass is considered to be distributed by a Dirac measure centered at its respective point in the flow. The compact approximation is achieved through the introduction of smooth mollifiers. 
Taking advantage of the compactness approximation, the following statement demonstrates the feasibility of approximating $\rr_{\lambda,\taus}$ by a finite-rank operator. 

\begin{cor}\label{cor: finite_1}
For any fixed $T>0$, for any arbitrarily small $\delta>0$, there exists a sufficiently large $N$ and a finite-dimensional approximation $\rr_{\lambda,\taus}^N$ such that 
$\|\rr_{\lambda,\taus}^N  h-\rr_{\lambda,\taus} h\|<\delta, \;\;h\in \mathcal{C}(\X).$
\end{cor}

\ym{The following exemplifies the construction of a finite-rank approximation using the dictionary functions (Definition \ref{def: dictionary}). 
\begin{deff}\label{def: projection}
    Consider $g\in \mathcal{C}(\X)\subset  L^2$. We define the Gram matrix $\bar{X}$ of the dictionary functions as $\bar{G}_{i,j}= \langle\zk_i, \zk_j\rangle_{\sss L^2}$,   $\bar{H}_g := \langle Z_N, g\rangle_{\sss L^2}$, and    the projection of $g$ on to $\Zk_N$ as $\Pi^Ng : = Z_N(x)^\trans \bar{G}^{\dagger}\bar{H}_g$. 
%We also define $Z_N(x) = Z_N(x)$.  
For a bounded linear operator $\B:\mathcal{C}(\X)\ra\mathcal{C}(\X)$, %define $\hat{Y} _\B= \int_\X Z_N(y)\B Z_N(y)dy$, $\hat{B} =\hat{X}^\dagger \hat{Y}_\B$,  and 
     we define %the projection of  $\J_\eps\B$ onto $\Zk_N$ as $\widehat{\Pi}^N \B h = Z_N^\eps(x)^\trans \hat{X}^{\dagger}\hat{Y}_{\B h}$. 
       $\Pi^N \B Z_N(x)^\trans = Z_N(x)^\trans \bar{G}^{\dagger}\bar{H}_{\B}$ with   $\bar{H} _\B= \int_\X Z_N(y)(\B Z_N(y))^\trans dy$. \Qed
\end{deff}

}

\begin{rem}\label{rem: construction}
%\begin{rem}\label{rem: finite-rank}
\ym{The fundamental properties of $\Pi^N$ are provided in the supplementary material. Note that it is common in the literature to construct  $K_t^N:=\Pi^N\K_t$, which is   proved to satisfy $\|\K_t^N(\cdot) - \K_t(\cdot)\|_{L^2}\xrightarrow{N\ra\infty}0$.  %%It then follows that the approximation error $|\K_t^N h-\K_t^\eps h|\lesssim \inf_{v\in\Zk_N}|h-v|_\infty$, which is small given $\Zk_N$ is sufficiently dense in $\mathcal{C}(\X)$. 
%As in \cite{williams2015data}, one can  consider a Hilbert space $\mathcal{H}:=L^2(\X)\supseteq \mathcal{C}(\X)$ \ym{and express $\K_t^N h(x)=Z_N(x)\hat{K}_t\theta$ as a Galerkin approximation for any $h(x) = Z_N(x)\theta$, where $\hat{K}_t = \hat{X}^{-1}\hat{Y}$ with $\hat{X}_{ij} = \langle \zk_i, \zk_j\rangle$ and $\hat{Y}_{ij} = \langle \zk_i,  \K_t\zk_j\rangle$. 
However, we cannot derive convergence w.r.t. the uniform norm $\|\cdot\|$ from $L^2$-convergence due to the possible lack of compactness of $\K_t(\Zk_N)$. Indeed, for any $g\in \mathcal{C}(\X)$, we have $\|\Pi^Ng-g\| \leq (1+\|\Pi^N\|)|E_N(g)$ by Lemma \ref{lem: error_sup} of the supplementary material. Although $E_N(g)$ converges as assumed, obtaining convergence of the r.h.s. is not straightforward if $\|\Pi^{N}\|$ is not uniformly bounded.

Nonetheless,  the analysis in Lemma \ref{lem: multi_converge} offers an alternative perspective. %By introducing an arbitrarily small smoothing parameter $\eps$, 
One can work with the compact operator family $\{\J_\eps\K_t\}$ and their finite-rank approximations $\Pi^{\eps, N}\K_t$,  which ensures convergence w.r.t. the uniform norm %$\|\cdot\|$ 
as $N\uparrow\infty$ for any fixed $\eps$. Furthermore, $\Pi^{\eps, N}\K_t$ converges   to $\Pi^N\K_t$ as $\eps\downarrow 0$ for any fixed $N$. Through this two-layer convergence, where $\eps$ can be scheduled to go to $0$ faster than the growth rate of $\|\Pi^N\|$, we obtain   the desired  convergence   $\K_t^N \sra\K_t$ w.r.t. $\|\cdot\|$.
}   \Qed
\end{rem} 

In view of the proof of Corollary \ref{cor: finite_1} and Remark \ref{rem: construction}, the finite-horizon finite-rank approximation 
$\rr_{\lambda,\taus}^N$ of $\rr(\lambda)$  exists and can be constructed by  the same procedure as $\K_t^N$.   \ym{Note that $\K_t^N$ is also the limit of the EDMD data-driven version as the number of samples tends to infinity.  In this view, the theoretical existence of the finite-rank form $\K_t^N$ (and hence $\rr_{\lambda, T}^N$) also plays the role of a hypothetical interface   for data-driven computation.  }

\section{Quadrature Approximation of the Resolvent and Reformulation of $\ll_\lambda$}\label{sec: quad}

\ym{Building on the feasibility of the finite-horizon finite-rank approximation $\rr_{\lambda, T}^N$ of $\rr(\lambda)$ in Section \ref{sec: finite-dim}, % from the preceding section, 
we investigate how to reformulate $\ll_{\lambda}$ when acting on the dictionary functions $Z_N(x)$, so that a directly accessible form of the output can be expressed in terms of the Koopman-related flow information. 

\subsection{Initial Attempt and   Conflict with   Practical Constraints}

We attempt to directly work on the form $\lambda^2 \rr_{\lambda,\taus}^N  -\lambda\id$ for sufficiently large $\lambda$ to generate output evaluation for $\ll_\lambda$ within a small time horizon corresponding to all observables \cite{meng2024koopman}.  The overall error consists of the analytical error of the Yosida approximation (Theorem \ref{thm: conv}) and the truncation error $E_{\text{tr}}$ (Theorem \ref{thm: conv_t}); the latter is exponentially faster compared to the former, so its contribution is comparatively negligible for sufficiently large $\lambda$. However, the output information requires numerical quadrature, using discrete-time observations of Koopman flow to approximate the integral of $\rr_{\lambda, T}^N$.}

%%Recalling the convergence rate of $\lll\sra\ll$  in Theorem \ref{thm: conv}, \ym{we see that the convergence in Theorem \ref{thm: conv_t} is exponentially faster, so its contribution to the error is negligible compared to that in Theorem \ref{thm: conv}.} This allows us to use $\ll_{\lambda,\taus}:=  \lambda^2 \rr_{\lambda,\taus}  -\lambda\id$, and hence   its finite-rank form (with a sufficiently large $N$)\begin{equation}\label{E: approx_t}\ll_{\lambda,\taus}^N:=  \lambda^2 \rr_{\lambda,\taus}^N  -\lambda\id, \;\;\lambda\ra\infty,  \end{equation}to approximate $\ll$ within a small time horizon.  

 %%In view of Remark \ref{rem: finite-rank}, once we collect the outputs corresponding to all observable functions at sampled initial conditions under the operation of $\ll_{\lambda, T}$, we can follow the exact same procedure to process the data and obtain convergence to $\ll_{\lambda,\taus}^N$. 

 Given the Laplace transform type integral, we employ the Gauss–Legendre quadrature rule to secure accuracy, using $\Gamma = \gamma T$ snapshots for each trajectory during data collection under observation frequency $\gamma$. \ym{Although it is beyond the scope of this work to evaluate the error bound of the quadrature rule for different choices of dictionary functions, we exemplify in Appendix D  that, for the monomial basis with maximal degree $\mathsf{N}$, an estimate of the error bound is proved to be $E_{\scriptscriptstyle\sum} \lesssim  T^{2\Gamma + 1} \left(\frac{\mu+\mathsf{N}L_f}{8\Gamma^2}\right)^{2\Gamma}$.  Other choices of dictionary functions follow a similar procedure and exhibit a similar pattern of error bounds, but involve nontrivial analysis; a detailed analysis is left for future work.
 
 %Now let $\rr_{\lambda,T}^{N,\text{quad}}$ denote the quadrature approximation of $\rr_{\lambda,T}^N$, which is  a finite weighted sum of Koopman flow evaluations at discrete observation times, whose existence and convergence is guaranteed by virtue of   Section~\ref{sec: sub_compact}. 
 By the spirit of the above analysis}, for larger $\lambda$ under a fixed small $\taus>0$, the numerical integration requires a larger $\Gamma$, which corresponds to a higher sampling frequency.  Conversely, if $\Gamma$ is restrictive in practice (for example, in automatic vehicles or intelligent transportation systems, where sensors can collect state data at no more than $\gamma = 100$ \cite{meng2024koopmanitsc}), $\lambda$ cannot be large to ensure accuracy when using numerical integration. %Making $\lambda$ and $\Gamma$ both sufficiently large cannot be achieved without modifications. 

 \subsection{Modified Evaluation Structure}
 To resolve this conflict, we employ the first resolvent identity $[(\lambda-\mu)\rr(\mu)+\id]\rr(\lambda) = \rr(\mu)$, which connects the two resolvent operators corresponding to  $\lambda$ and $\mu$ in the resolvent set.
This identity can be reformulated as \begin{small}
    \begin{equation}\label{E: formula_2}
    \begin{split}
    & [(\lambda-\mu)\rr(\mu)+\id] \ll_\lambda h(x) \\
  = &   \lambda ^2 \rr(\mu) h(x) - \lambda(\lambda-\mu)\rr(\mu)h(x)-\lambda h(x)\\
  = & \lambda\mu\rr(\mu)h(x) -\lambda h(x), 
    \end{split}
    \end{equation}
\end{small}

\noindent %thereby connecting the two resolvent operators corresponding to different values $\lambda$ and $\mu$ in the resolvent set. 
 allowing us to learn a resolvent operator with a small $\mu$ (thus enabling the use of a small $\Gamma$ for numerical integration) and then infer the $\ll_\lambda$ with a large $\lambda$.
\ym{The form of \eqref{E: formula_2} requires transferring the output of $\ll_\lambda$ through  
$(\lambda-\mu)\rr(\mu)+\id$\ and matching it with the  r.h.s.,  which again necessitates a prior finite-horizon quadrature evaluation $\rr_{\mu,T}^{\text{quad}}$ of $\rr(\mu)$. 

\begin{rem}\label{rem: fail}
    %To provide an interface through 
    To enable learning via an analytical finite-horizon finite-rank approximation $\ll_{\lambda}^N$ of  $\ll_\lambda$, we first attempt to follow the same spirit of Section~\ref{sec: sub_compact}  %define $\ll_{\lambda}^N := \Pi^N\ll_\lambda$, a
    and evaluate  $\ll_\lambda$ %\eqref{E: formula_2}  as 
via $ [(\lambda-\mu)\rr_{\mu,T}^{\text{quad}}+\id] \ll_{\lambda}^N h(x)  
  =  \lambda\mu\rr_{\mu,T}^{\text{quad}}h(x) -\lambda h(x)$.
In this formulation,  both sides incur truncation error $E_{\text{tr}}$ and quadrature  error $E_{\sss \sum}$, with an additional projection error    on l.h.s.,  %with each absolutely continuous to $\frac{\mu^2}{\mu-\omega}e^{-\mu T}$,  $\max_{h\in\K_t Z_N}\mathcal{O}_N(h)$, and $\mathcal{O}_{\sss \sum}$,
compared to  \eqref{E: formula_2}. Note that the latter two errors can be made sufficiently small with a suitable dictionary and small 
$\mu$, whereas the $E_{\text{tr}}$ may remain large due to the conflict between choosing small $T$ and small $\mu$.  \Qed
\end{rem}

In view of Remark \ref{rem: fail}, directly replacing 
$\rr(\mu)$ with $\rr_{\mu,T}^{\text{quad}}$ in \eqref{E: formula_2} may cause significant truncation error.
To keep the quadrature error small by using appropriately small $T$, $\Gamma$, and $\mu$  while simultaneously reducing   truncation error, we proceed with the following fact  for a further modification.

\begin{prop}\label{prop: contraction}
     For any $h\in \mathcal{C}(\Omega)$ and a fixed $T>0$, we have $\rr(\mu)h(x)=\rr_{\mu, T}h(x)+e^{-\mu T}\K_T\rr(\mu) h(x)$. %let $I(x)=\int_0^T e^{-\mu t} \K_th(x)dt$. Then $$\rr_{\mu}h(x)=I(x)+e^{-\mu T}\K_T\rr_\mu h(x).$$
     Moreover, by defining operator $\mathcal{T}_h$ as  $\mathcal{T}_h (\cdot) =\rr_{\mu,T}h +e^{-\mu T}\K_T (\cdot) $, it 
      possesses the contraction mapping property.
 \end{prop}

We then  evaluate
$\rr(\mu)\zk_i$    by using an ansatz $Z_N^\trans\zeta_i\in\Zk_N$, where 
$\zeta_i$ is determined by plugging the ansatz into the above equation. %with the (fixed) finite-horizon evaluation 
%$\rr_{\mu, T}^N$. 
%Solving for $\zeta_i$ then yields an error estimated as follows.
The solvability is established as follows.

\begin{lem}\label{lem: unique_proj}
   Recall $\Pi^N$ from Definition \ref{def: projection}. For a fixed $T>0$, let $\K_T^N=\Pi^N\K_T$. For each $\zk_i\in Z_N$,  let $V_i(x) = \rr(\mu)\zk_i(x)$, let   $E_{N}^i(x): = \K_T^NV_i(x)-\K_TV_i(x)$, and denote $E_N^\text{max}:=\max_{i}\|E_{N}^i\|$.  Then, there is $\tilde{\omega}$ such that $\tilde{\omega}\lesssim \omega + E_N^\text{max}$; and for all $\mu\in(\tilde{\omega},\infty)$, there is a unique solution $V_i^N=Z_N^\trans\zeta_i\in\Zk_N$ to $V_i^N(x) = \rr_{\mu,T}^{\text{quad}}\zk_i(x) + e^{-\mu T}\K_T^N V_i^N(x)$. Furthermore, $\|V_i^N-V_i\|\lesssim E_{\sss \sum}+ E_N^\text{max}$.
\end{lem}

%This property implies that one can evaluate 
The above statement provides a way to evaluate $\rr(\mu)Z_N(x)$ via $V^N(x) := [V_i^N(x), i\in\{0, 1, \cdots, N-1\}]^\trans$, and, equivalently, yields a finite-rank representation of $\rr(\mu)$ given by %$\rr_{\mu}^NZ_N(x)$ Now, let $\bar{H}_{\sss\rr_\mu} = \int_\X Z_N(x)V(x)^\trans dx$.  Then we can construct a finite-rank representation of $\rr(\mu)$ by defining 
\begin{small}
    \begin{equation}
        \rr_{\mu}^NZ_N(x) = V^N(x) = Z_N(x)^\trans[\zeta_0, \zeta_1, \cdots, \zeta_{N-1}]^\trans.
    \end{equation}
\end{small}

\noindent %as in Definition \ref{def: projection}, with 
It is clear that the approximation error depends continuously on $\max_i\|V_i^N - V_i\|$. Unlike $\rr_{\mu, T}^N$,  the construction of $\rr_{\mu}^N$ eliminates  $E_{\text{tr}}$ and does not produce large errors for well-chosen small values of $\mu$, $\Gamma$, and $T$. We next build on $\rr_{\mu}^N$ to obtain a finite-rank approximation of $\ll_\lambda$ based on \eqref{E: formula_2}.

\begin{thm}\label{thm: approx_modify}
    Let $\widehat{\A}:=(\lambda-\mu)\rr_\mu^N+\id$ and $\widehat{\B}:=\lambda\mu\rr_\mu^N-\lambda\id$. Let $\bar{G}$, $\bar{H}_{\widehat{\A}}$, and $\bar{H}_{\widehat{\B}}$ be defined as in Definition \ref{def: projection}. Let $\hat{A} = \bar{G}^{\dagger}\bar{H}_{\widehat{\A}}$, $\hat{B} = \bar{G}^{\dagger}\bar{H}_{\widehat{\B}}$, and   $\hat{A}_\delta^+ = (\hat{A}^\trans \hat{A}+\delta \id)^{-1}\hat{A}^\trans$ for any $\delta>0$. Define $\hat{L}_\lambda:=\hat{A}^{\dagger}_\delta\hat{B}$ and   $\ll_\lambda^N h_\theta = Z_N(x)^\trans \hat{L}_\lambda \theta$ for any $h_\theta=Z_N(x)^\trans\theta\in\Zk_N$. %  with $\theta\in\operatorname{ran}(\hat{A})$ and $\hat{B}\hat{A}_\delta^+\theta\in\operatorname{ran}(\hat{A}^\trans)$. 
    For  $E_A:=E_{\sss \sum}+ E_N^\text{max}<1$, there exists a $\delta = \mathcal{O} (E_A^{1/3})$ and a projection residual $E_N^\Pi\xrightarrow{N\ra\infty} 0$ such that $\|\ll_\lambda^Nh_\theta-\ll_\lambda h_\theta\| 
      \lesssim  E_N^\Pi+  E_A^{1/3}$. 
\end{thm}

\begin{rem}
   In fact, the proof also implies that for any fixed $\delta>0$, the $1/3$-Hölder bound improves to a linear  bound. The above analysis provides only a worst-case heuristic for tuning $\delta$, since the rank of $\hat{A}$ may fail to be preserved as $N\ra\infty$ and  errors vanish. In the ideal case where this issue does not arise, one can use 
    \begin{small}
            \begin{equation}\label{E: formula_L}
        \ll_\lambda^Nh_\theta   = Z_N(x)^\trans \hat{L}_\lambda\theta, \;\;\hat{L}_\lambda:=\hat{A}^{\dagger}\hat{B}, \;h_\theta = Z_N(x)^\trans\theta
    \end{equation}
    \end{small}

    \noindent and perform a perturbation analysis  \cite{wedin1973perturbation} directly on it, in which case the convergence rate becomes linear as well. We omit the proof because of repetition. To better convey the key idea of the learning procedure,  the rest of this paper uses the formula \eqref{E: formula_L} for simplicity. Readers should be aware that a $\delta$-regularization is required to guarantee worst-case theoretical convergence when $\hat{A}$ is rank-deficient. However, it is not always necessary in practice, as demonstrated by data-driven numerical experiments that exhibit ideal convergence behavior even without the regularization. \Qed
\end{rem}

We have   reprocessed   $\rr(\mu)$ with improved evaluation precision  %based on its finite-horizon truncation and quadrature, 
in order to adapt to the reformulation    \eqref{E: formula_2}. The resulting formula~\eqref{E: formula_L} serves as a direct interface for generating training data, ensuring that the learning of $\mathcal{L}_\lambda$ remains effective even under sampling-frequency constraints.  

%We have preprocessed the integral representation of $\rr(\mu)$ into a form that acts as a direct interface for generating training data, ensuring that the learning of $\ll_\lambda$ via \eqref{E: formula_L} remains effective even under sampling-frequency constraints.
}

 %With this form in place, the procedure reduces to one similar to the approaches in \cite{mauroy2019koopman,zeng2024sampling,meng2024koopman}.

\iffalse
 This property implies that one can evaluate the output of $\rr(\mu)$ using the (fixed) finite-horizon evaluation $\rr_{\mu, T}$, together with a  dynamic programming   to update the output from any initial guess, thereby achieving an acceptable error within a finite iteration. We state the following classical result and omit the proof.

 \begin{cor}
    For any $h\in \mathcal{C}(\Omega)$ and a fixed $T>0$, let $V^{(n)}(x)=\rr^{(n)}(\mu)h(x)$ be the output of $\rr(\mu)$ of the $n$-th iteration, where $V^{(n)}(x)$ follows the iteration
    $$V^{(n+1)}(x) = \rr_{\mu, T}h(x)+ e^{-\mu T}  \K_TV^{(n)}(x).$$
    Then, as $n\ra\infty$, $V^{(n)}(x)\ra\rr(\mu) h(x)$. The error is such that $|V^{(n)}-\rr(\mu)h(x)|_\infty\leq \frac{e^{-n(\omega-\mu)T}}{1-e^{-(\omega-\mu)T}}|V^{(1)}-V^{(0)}|_\infty$. 
\end{cor}\fi

\section{Data-Driven Algorithm}\label{sec: alg}

We continue to discuss the data-driven learning based on the approximations discussed in Section \ref{sec: finite-dim} and \ref{sec: quad}. %Recall Definition \ref{def: T_t}, Eq. \eqref{E: approx_t}, and Corollary \ref{cor: finite_1}.  
\ym{Similar to the learning of Koopman-related operators \cite{williams2015data, mauroy2019koopman, meng2023learning, zeng2024sampling, meng2024koopman}, it suffices to estimate $\hat{L}_\lambda$ (defined with respect to the matrix-valued $L^2$-inner product) from \eqref{E: formula_L} using a Monte Carlo–style data-driven integration $L_\lambda$ based on sampled initial conditions $\{x^{(m)}\}_{m=0}^{M-1}\subseteq\X$ of \eqref{E: sys}. The convergence of $L_\lambda$ to $\hat{L}_\lambda$ is identical to \cite{williams2015data} as $M \to \infty$, and we   do not repeat the analysis.}
\iffalse
Similar to the learning of Koopman-related operators \cite{williams2015data, mauroy2019koopman, meng2023learning, zeng2024sampling,meng2024koopman}, obtaining a fully discretized version  $\Lb$  of the bounded linear operator $\lambda^2\rr_{\lambda,\taus}^N-\lambda \id$ based on the training data typically relies on the selection of a discrete dictionary of continuously differentiable observable test functions, denoted by
\begin{equation}\label{E: dict}
    \Zk_N(x):=[\zk_{\scriptscriptstyle  0}(x), \zk_{\scriptscriptstyle  1}(x),  \cdots, \zk_{\scriptscriptstyle  N-1}(x)], \;N\in\N. 
\end{equation}\fi
%$$ for $N\in\N$. 
%Then, the following should hold: 
\iffalse
the learned eigenvalues $\mu_i$ can be used for approximate the true eigenvalues of $\lambda^2\T_t^N-\lambda \id$ (and hence those of $\ll$), and the  eigenfunction $\zeta_i$ of $\ll$ can be approximated by   $\Zk_N(x)\vb_i$.\fi 
\ym{Consequently,   
for sufficiently large $\lambda$ and any $h_\theta\in\Zk_N$ (i.e. $h_\theta(x)=Z_N(x)^\trans\operatorname{\theta}$ with some   vector $\theta$),  
\begin{small}
    \begin{equation}\label{E: approx_Lh}
   \ll h_\theta(\cdot)\approx \ll_\lambda h_\theta\approx \ll_{\lambda}^N h_\theta(\cdot) \approx Z_N(\cdot)^\trans(\Lb_\lambda\theta).
\end{equation}
\end{small} 

\noindent  In this section, we present the procedure  for obtaining $L_\lambda$. } %In this section, we modify the existing Koopman learning technique to obtain $\Lb$ such that   \eqref{E: approx_eigen} and \eqref{E: approx_Lh} hold. 

\begin{rem}
Denoting $\widetilde{\ll}h_\theta(x):=Z_N(x)^\trans(\Lb_\lambda\theta)$,   the   approximation in \eqref{E: approx_Lh}    also guarantee the   convergence of $\{e^{t \widetilde{\ll}}\}$ to the original semigroup $\{\K_t\}_{t\geq 0}$ for any $h_\theta(x)=Z_N(x)^\trans \theta$. Indeed, %one can show that, for any $\tomg\geq \omega$, we have $(\ll-\tomg\id)_{\lambda-\tomg} + \tomg\id = \ll_\lambda + O_2$, where $O_2$ comes from the proof of Theorem \ref{thm: conv}. For simplicity, 
    letting $\mathcal{O} = \widetilde{\ll} - \ll_\lambda$, we have  
    \begin{small}
            \begin{equation}
        \begin{split}
            &\|e^{t\widetilde{\ll}}h_\theta - \K_th_\theta\|\\  \leq & \|e^{t\mathcal{O}}\| \|e^{t\ll_\lambda}h_\theta - \K_th_\theta\| +  \|\K_t\|\|e^{t\mathcal{O}}h_\theta-h_\theta\|\\
             \leq &\|e^{t\mathcal{O}}\| \|e^{t\ll_\lambda}h_\theta - \K_th_\theta\| +  t\|\K_t\|e^{t\|\mathcal{O}\|}\|\mathcal{O}h_\theta\|. 
        \end{split}
    \end{equation}
    \end{small}

\noindent %Note that $\mathcal{O}$ and $\K_t$ (for any $t$) are  bounded operators. 
As $\lambda \ra\infty$, the first part goes to $0$ by Theorem \ref{thm: H-Y}. The second part goes to $0$ by  the sequence of approximations (in a strong sense) in \eqref{E: approx_Lh}. \Qed
\end{rem}

\subsection{Data Collection, Pre-processing, and the Algorithm}
The training data is obtained in the following way. 
\subsubsection{Features data}
By randomly sampling $M$ initial conditions $\{x^{(m)}\}_{m=0}^{M-1}\subseteq\X$, we stack the \textit{features} into $\Xb$:
\begin{equation}\label{E: stack_X}
    \Xb=[Z_N(x^{(0)}), Z_N(x^{(1)}),\cdots, Z_N(x^{(M-1)})]^\trans.
\end{equation}

\subsubsection{Numerical integration}
\ym{Now, fix a $T$. For any fixed $\mu > 0$ and each $x^{(m)}$, we use discrete-time observations (snapshots) of   Koopman flow $\zk_i(\phi(t, x^{(m)}))$ to approximate $\rr_{\sss{\mu,T}}\zk_i$ using numerical quadrature techniques.}
%For any fixed $\mu>0$ and each $x^{(m)}$ we use the discrete-time observations (snapshots) of the Koopman flow $\phi(\taus,x^{(m)})$ to evaluate the $\rr_{\mu, T}$ using numerical quadrature techniques for approximation. given a dictionary $\Zk_N$ of the form \eqref{E: dict}, for each $\zk_i\in\Zk_N$ and each $x^{(m)}$, we consider %$\zk_i(x)$ as the features and  $\lambda^2\rr_{\lambda,\taus}\zk_i(x^{(m)})-\lambda\zk_i(x^{(m)})=\lambda^2\int_0^\taus e^{-\lambda s} \zk_i(\phi(s,x^{(m)})) ds-\lambda\zk_i(x^{(m)})$ as the labels.
%To implement this numerical quadrature, we 
Recall the observation  rate of $\gamma$ over the interval $[0, T]$ for each flow map $\phi(\cdot, x^{(m)})$, 
%To implement the %compute the integral, we employ  numerical quadrature,  techniques for approximation. This approach inevitably requires discrete-time observations (snapshots) with an observation/sampling rate of $\gamma$ Hz within the interval $[0,\taus]$ of each flow map $\phi(\taus,x^{(m)})$. 
and the number of snapshots as $\Gamma:=\gamma\taus$.  We stack the snapshot data of the integrand in the following intermediate matrix:
\begin{small}
    \begin{equation}\label{E: stack_U}
\begin{split}
      U^{(m)}
   & =  [Z_N(\phi(0,x^{(m)})),  \cdots, e^{-\frac{\mu k\taus}{\Gamma}} Z_N(\phi(\frac{k\taus}{\Gamma},x^{(m)})),\\
    & \cdots e^{- \mu\taus} Z_N(\phi(T,x^{(m)}))]^\trans, \; m\in\{0, 1, \cdots, M-1\}.
\end{split}
\end{equation}
\end{small}

\noindent Denote by $\mathcal{G}^\lambda_{[0, \taus]}(v)$, or   $\mathcal{G}(v)$ for brevity, the Gauss–Legendre quadrature\footnote{We omit the details of implementing the Gauss–Legendre quadrature numerically, as the algorithms are well-established in this field.} based on the vector of snapshot points $v$ within $[0, \taus]$, and denote by $U^{(m)}[:, j]$ the $j^{\text{th}}$ column of $U^{(m)}$. The stack of $\rr_{\mu, T}^{\text{quad}}Z_N(x^{(m)})$ is given by % \footnote{The superscripts and subscripts in the notation indicate the dependence of the parameters, reflecting that the image function values of the generator are only achieved approximately. The $\lambda$ and $\taus$ determine the analytical error, while the $\Gamma$ determines the numerical error. When the dependence on any parameter is not emphasized, we can use shorthand notation without the corresponding superscript and subscript.} 
%\begin{equation}\label{E: stack_Y}
   %\Yb \;(=\Yb_{\lambda,\taus}^\Gamma):=\I_{\lambda,\taus}^\Gamma-\lambda\Xb,
%\end{equation}
% where
 \begin{small}
     \begin{equation}\label{E: num_inte}
    % \begin{split}
       \I_{\mu, T}^{\text{quad}}   =
        \begin{bmatrix}
            \mathcal{G}(U^{(0)}[:, 0])&   \cdots & \mathcal{G}(U^{(0)}[:, N-1])\\
            \vdots &   \ddots & \vdots\\
            \mathcal{G}(U^{(m)}[:, 0])&    \cdots & \mathcal{G}(U^{(m)}[:, N-1])\\
             \vdots &   \ddots & \vdots\\
             \mathcal{G}(U^{(M-1)}[:, 0])&   \cdots & \mathcal{G}(U^{(M-1)}[:, N-1])\\
        \end{bmatrix}.
    % \end{split}
\end{equation}
 \end{small}

\ym{\subsubsection{Estimation of $\rr_{\mu}^N$}
We then infer $\rr_{\mu}^NZ_N$ from $\rr_{\mu, T}^{\text{quad}}Z_N$ by Lemma \ref{lem: unique_proj},  which relies on using ansatz $\rr_{\mu}^N\zk_i:=Z_N^\trans \zeta_i$ to solve $\rr_{\mu}^N\zk_i - e^{-\mu T}\K_T\rr_{\mu}^N\zk_i= \rr_{\mu, T}^{\text{quad}}\zk_i$  for all $\zk_i$. Equivalently, the stacked form of these equations evaluated at $\{x^{(m)}\}$ is given by $ X\Xi - e^{-\mu T}\Phi_T\Xi = \I_{\mu, T}^{\text{quad}}$, 
where 
\begin{small}
    \begin{equation} \label{E: Phi_tau}
    \Phi_{\taus} := [Z_N(\phi(\taus,x^{(0)})),  \cdots, Z_N(\phi(\taus,x^{(M-1)}))]^\trans 
\end{equation} 
\end{small}

\noindent and $\Xi\in\R^{N\times N}$ is the stack of weights $\zeta_i$ with the solution given by $\Xi=(X-e^{-\mu T}\Phi_T)^\dagger\I_{\mu, T}^{\text{quad}}$. The stack of $\rr_\mu^NZ_N$ evaluated at $\{x^{(m)}\}$ is then given by $X\Xi$.
 
\subsubsection{Final inference of $L_\lambda$}
We   follow Theorem \ref{thm: approx_modify} and \eqref{E: formula_L} to construct the stacks for $\hat{A}$ and $\hat{B}$.  As in Theorem \ref{thm: approx_modify}, the valuation of $\widehat{\A}$ and $\widehat{\B}$ at $\{x^{(m)}\}$ is such that 
\begin{equation}\label{E: YAYB}
    Y_A = (\lambda-\mu) \Xb \Xi + \Xb \;\text{and}\; Y_B = \lambda \mu \Xb \Xi -\lambda\Xb,
\end{equation}
respectively. Let $A = \left(\Xb^\trans\Xb\right)^\dagger\Xb^\trans Y_A$ and $B = \left(\Xb^\trans\Xb\right)^\dagger\Xb^\trans Y_B$, then $\hat{A}$ and $\hat{B}$ are the limits of $A$ and $B$ as $M\ra\infty$, following the same argument as in \cite{williams2015data}. Then, $L_\lambda = A^\dagger B$ as required. 

%We omit the algorithm for generating the training data $(\Xb, \Yb)$, as it follows straightforwardly from the data stacking shown in   \eqref{E: stack_X} and \eqref{E: stack_Y}.

\subsubsection{Summary of algorithm} A summary of the algorithm for computing $L_\lambda$ is provided in Algorithm \ref{alg: resolvent-learning}. %We appreciate the reader’s patience in following the heavy explanations up to this point. The rest of the paper is lighter and more case-study-oriented, beginning with a summary of the algorithm of computing $L_\lambda$. 

\begin{algorithm}
	\caption{Resolvent-Type Koopman Generator Learning }
    \label{alg: resolvent-learning} 
	\begin{algorithmic}[1]
		\REQUIRE Dictionary $Z_N$, user-defined $\mu$ and $\lambda$, initial conditions $\{x^{(m)}\}_{m=0}^{M-1}\subseteq \X$, $\taus$, $\Gamma$, 
        and snapshots $\phi(\frac{k\taus}{\Gamma},x^{(m)})$ for $k=0,1,\ldots,\Gamma$.
		
		% \FOR{$m$ \textbf{from} $0$ \textbf{to} $M-1$}
		% 	\FOR{$i$ \textbf{from} $0$ \textbf{to} $N-1$}
		% 		\STATE Compute $\zk_i(x^{(m)})$ 
		% 		\STATE Stack $U^{(m)}$ using \eqref{E: stack_U}, and compute $\mathcal{G}(U^{(m)}[:,i])$ for each observable test function $\zk_i$ using \eqref{E: num_inte}.
		% 	\ENDFOR
		% \ENDFOR
		
		\STATE Compute and stack $\Xb$ using \eqref{E: stack_X};
		
		%$$\Xb=[\Zk_N(x^{(0)}), \Zk_N(x^{(1)}),\cdots, \Zk_N(x^{(M-1)})]^\trans$$

        \STATE Compute and stack $U^{(m)}$  using \eqref{E: stack_U};

		\STATE Compute and stack  $\rr_{\mu, T}^{\text{quad}}Z_N(x^{(m)})$  using $\I_{\mu, T}^\text{quad}$ by \eqref{E: num_inte};
		
		\STATE Stack $\Phi_{\taus}$ using \eqref{E: Phi_tau};
		
		\STATE Compute $\Xi=(X-e^{-\mu T}\Phi_T)^\dagger\I_{\mu, T}^{\text{quad}}$
		
		\STATE Compute matrices $Y_A$ and $Y_B$ using \eqref{E: YAYB};
		\STATE Compute $A = \left(\Xb^\trans\Xb\right)^\dagger\Xb^\trans Y_A$, $B = \left(\Xb^\trans\Xb\right)^\dagger\Xb^\trans Y_B$. 
		
		\RETURN $\Lb_{\lambda}= A^\dagger B$ 
	\end{algorithmic}
\end{algorithm}}

\subsection{Discussion of Other Existing Methods}\label{sec: comparison}
Recall Remark \ref{rem: FDM} on the FDM with an expression $\frac{\K_\tau-\id}{\tau} \sra \ll$, where the convergence is well studied in \cite{bramburger2024auxiliary}. We also revisit the benchmark Koopman Logarithm Method (KLM) based on the expression $\ll = \frac{1}{\tau}\log(\K_\tau)$, as described in \cite{mauroy2019koopman}. In this subsection, we discuss the data-driven algorithms for FDM and KLM, in preparation for the comparison in the case studies in Section \ref{sec: num}. 

Observing that the expressions for $\ll$ in FDM and KLM rely on just one moment of the Koopman operator $\K_\tau$, the data-driven versions of these methods are divided into two steps: 1) learning $\K_\tau$; 2) transforming the learned $\K_\tau$ to $\ll$, respectively. To make the  sampling   rate consistent with the RTM, we  set $\tau:= \taus/\Gamma$. The corresponding data-driven approximations  also rely on the selection of a discrete dictionary $Z_N$, and similarly, 
 the approximation is to achieve $ \ll h_\theta(\cdot)\approx Z_N(\cdot)(\Lb_i\theta)$ for $h_\theta(x)=Z_N(x)\theta$ and for any $i\in\set{\operatorname{FDM}, \operatorname{KLM}}$. To do this, we fix a $\tau$,   stack the \textit{features} into $\Xb$ the same way as \eqref{E: stack_X}, and and stack the labels $\Phi_\tau$ same as \eqref{E: Phi_tau}. 
After obtaining the training data $(\Xb, \Phi_\tau)$, we can find $\Kb$ (the fully discretized version of $\K_\tau$) by  $\Kb = \operatorname{argmin}_{A\in\C^{N\times N}}\|\Phi_\tau-\Xb A\|_F = \left(\Xb^\trans\Xb\right)^\dagger\Xb^\trans \Phi_\tau$ \cite{williams2015data}. The data-driven approximation for $\ll$ based on FDM and KLM are immediately obtained using $\Kb$. 
\begin{enumerate}
    \item FDM: $\Lb_{\operatorname{FDM}} = (\Kb-\id)/\tau$, where $\id=\id_{N\times N}$ is the dentity matrix in this expression. 
    \item KLM: $\Lb_{\operatorname{KLM}} = \log(\Kb) /\tau$. 
\end{enumerate}

\iffalse
\begin{rem} \label{rem: log_explain}
    Let $(\mu_i^K, \vb_i^K)_{i=0}^{N-1}$ be the eigenvalues and eigenvectors of $\Kb$. Let $(\rho_i^K, \varphi_i^K)_{i=0}^{N-1}$ be the eigenvalues and eigenfunctions of $\K_\tau$.
    Similar to $\eqref{E: approx_eigen}$ and \eqref{E: approx_Lh}, 
    \begin{equation}\label{E: approx_eigen}
    \mu_i\approx\rho_i,\;\;\varphi_i(x)\approx \Zk_N(x)\vb_i. 
\end{equation}
    for each $i$, we have $\varphi_i^K(x)\approx \Zk_N(x)\xi_i^K$ and $ \K_s h(\cdot)\approx \Zk_N(\cdot)(\Kb\theta)$ for $h(x)=\Zk_N(x)\theta$. \Qed
\end{rem}\fi

   It is worth noting that,  even when $\ll$ can be represented by $(\log{\K_\tau})/\tau$, we cannot guarantee that 
    $\frac{\log{\K_\tau}}{\tau}h(\cdot)\approx \Zk_N(\cdot)(\frac{\log(\Kb)}{\tau}\theta)$ as in KLM, not to mention the case where the above %(l.h.s.) operator 
    logarithm representation does not hold.  As pointed out in \cite[Remark 4.1]{meng2024koopman}, even though the (possibly complex-valued) matrix connecting Koopman eigenfunctions and dictionary functions ensures that any Koopman representation using $Z_N$ can be equivalently expressed in terms of eigenfunctions with cancellation of the imaginary parts, this cancellation effect generally does not hold when applying the matrix logarithm. An exception occurs only when the chosen dictionary is inherently rotation-free w.r.t. the true eigenfunctions \cite{zeng2024sampling}, or when there is direct access to data for $\log(\K_\tau)$ that allows direct training of the matrix. Such conditions, however, contradict our objective of leveraging Koopman data to infer the generator.

    %Denoting $\Phi(\cdot)=[\varphi_0^K(\cdot), \varphi_1^K(\cdot),\cdots,\varphi_{N-1}^K(\cdot)]$ and $\Xi=[\xi_0^K, \xi_1^K, \cdots,\xi_{N-1}^K]$, then it is clear that $\Phi(\cdot) = \Zk_N(\cdot)\Xi$. In view of Remark \ref{rem: finite-rank}, any $\K_\tau h(x)$ for any $h(x)=\Zk_N(x)\theta$ can be approximated using $\Phi$. The (possibly complex-valued) rotation matrix $\Xi$ establishes the connection between finite-dimensional eigenfunctions and dictionary functions through data-fitting, ensuring that any linear combination within $\Zk_N$ can be equivalently represented using $\Phi$ with a cancellation of the imaginary parts. 
    
%This imaginary-part cancellation effect does not generally hold when applying the matrix logarithm. Suppose the imaginary parts account for a significantly large value, the mutual representation of $\Phi$ and $\Zk_N$ does not match in the logarithmic scale. An exception holds unless the chosen dictionary is inherently rotation-free with respect to the true eigenfunctions \cite{zeng2024sampling}, or there is direct access to the data for  $\log(\K_\tau)$ allowing for direct training of the matrix. However, such conditions contravene our objective of leveraging Koopman data to conversely find the  generator. 

In comparison, the FDM and the resolvent-type approach approximate $\ll$ regardless of its boundedness. These two methods enable learning of $\Lb$ without computing the logarithm, thus avoiding the potential occurrence of imaginary parts caused by basis rotation. However, as we will see in Section~\ref{sec: num}, the performance of the FDM  degrades as the sampling frequency decreases.

\section{Case Studies}\label{sec: num}
In this section, we test the effectiveness of the RTM %proposed %resolvent-type 
method and compare its performance to the two benchmark methods, FDM and KLM.  
Particularly, we present numerical examples that apply the learned generator for system identification\footnote{\ym{The primary purpose of the proposed method is Koopman generator learning, and the fairest comparison is against other Koopman-related methods. Regarding system identification using the learned generator, we acknowledge that Koopman-generator-based approaches for this purpose are still in the early stages of development. Nevertheless, we include comparisons with the more mature methods SINDy/WSINDy to demonstrate the potential of our proposed approach in this setting.}} and effectively solve transport-related PDEs. The research code can be found at   \href{https://github.com/RuikunZhou/Resolvent-Type-Operator-Learning}{{https://github.com/RuikunZhou/Resolvent-Type-Operator-Learning}}.

%%Particularly, we perform numerical examples for system identification of the following representative systems \cite{mauroy2019koopman, zeng2023generalized}: the reversed  Van der Pol oscillator (polynomial), two-machine power system (nonpolynomial), a system with rational vector fields (nonpolynomial), Lorenz-63 system (chaotic), and Lorenz-96 system (high-dimensional chaotic). We also show that the RTM is able to solve transport-related PDEs effectively.

%\url{https://github.com/Yiming-Meng/Resolvent-type-Learning-of-Koopman-Generators}. %We would like to show that, with the same set of trajectory data collected within a short time horizon, one can predict long-term trajectory behaviours as well as infinitesimal transitions.

We provide two measurements to demonstrate the error in system identification. 

(1) The root mean squared error (RMSE) of flow 
\begin{small}
    \begin{equation}\label{E: rmse_OBO}
    \mathcal{E}_{\operatorname{RMSE}}^{\text{F}} = \frac{1}{M} \sum_{m=0}^{M-1} \sqrt{\frac{1}{\Gamma_s}\sum_{k=1}^{\Gamma_s}|\phi(t_k, x^{(m)})-\hat{\phi}(t_k, x^{(m)})|^2}
\end{equation}
\end{small}

\noindent measures the root mean square difference between the actual time-series data $\{\phi(t_k, x^{(m)})\}_{k=0}^{\Gamma_s}$ and the estimated data $\{\hat{\phi}(t_k, x^{(m)})\}_{k=0}^{\Gamma_s}$ using the learned generator up to time $T_s$. Here, $\Gamma_s$ represents the number of snapshots used to verify the performance (which is independent of $\Gamma$, used for the data collection procedure), and the time $t_k = kT_s/\Gamma_s$ corresponds to the sampling instances. 

(2) For polynomial models, and if we use monomial dictionary functions, % employing Koopman-based methods, 
we can use the RMSE of the weights assigned to each monomial %basis
\begin{small}
    \begin{equation}
    \mathcal{E}_{\operatorname{RMSE}}^{\text{W}} : = \sqrt{\frac{1}{dN}\sum_{j=1}^{d}\sum_{k=0}^{N-1} |\theta_{k}^j-\hat{\theta}_{k}^j|^2}.
\end{equation}
\end{small}

\noindent Here, $N$ is the size %total number of 
the dictionary; $d$ is the dimension of the system; $\theta_{k}^j$ is the weight for $f_j$ at the $k^{\text{th}}$ dictionary function; $\hat{\theta}_{k}^j$ is similarly defined but for the estimated vector field.

\subsection{Scaled Lorenz-63 System}

%Different from the widely used Lorenz-63 system, 
We first consider the system identification of   Lorenz-63 system~\cite{axenides2015scaling}  by scaling it with a factor of $\epsilon$, which yields
\begin{align}
\nonumber
\dot{\xb}_1 &= \sigma (\xb_2  - \epsilon \xb_1), \\
\dot{\xb}_2 &= \xb_1 (\gamma - \xb_3 ) - \epsilon \xb_2,  \nonumber
\\
\dot{\xb}_3  &= \xb_1  \xb_2  - \epsilon \beta \xb_3 , \nonumber
\end{align}
where $\sigma = 10$, $\gamma = 0.28$, $\beta = \frac{8}{3}$, and $\epsilon = 0.1$. With this so-called $\epsilon-$Lorenz system, the system possesses an attractor on $\X = (-1, 1)^3$. We choose the dictionary of monomials with total number of $N = P\times Q\times J$  and set $\zk_i(x) = x_1^p x_2^q x_3^j$, where $p=i-PQj-Pq$, $q = \lfloor \frac{i-PQj}{P}\rfloor$, and $j = \lfloor \frac{i}{PQ}\rfloor$. %Similar to Section \ref{sec: van der pol}, 

The actual vector field is $f(x) :=[f_1(x), f_2(x), f_3(x)]^\trans= [\sigma (x_2 - \epsilon x_1), x_1(\gamma - x_3) - \epsilon x_2, x_1x_2 -\epsilon\beta x_3]^\trans$, and one can analytically establish that $f_1 = \ll\zk_{\scriptscriptstyle  1}$, $f_2=\ll\zk_{\scriptscriptstyle  P}$, and $f_3 = \ll\zk_{\scriptscriptstyle  P+Q}$. After learning the generator, we use the  approximation $[\widetilde{\ll} \zk_{\scriptscriptstyle  1}, \widetilde{\ll}\zk_{\scriptscriptstyle  P}, \widetilde{\ll}\zk_{\scriptscriptstyle  P+Q}]^\trans$ to conversely obtain $f$. 

\begin{rem}\label{rem: scaling}
  The original system (as studied in \cite{mauroy2019koopman}) has the vector field $\hat{f}(x) = [\sigma (x_2 - x_1), x_1(\gamma - x_3) - x_2, x_1x_2 -\beta x_3]^\trans$.   Denoting the solution to the original Lorenz model as $\hat{\xb}$, it can be verified that 
    $\hat{\xb}_1(t) =  (1/\epsilon) \xb_1(\epsilon t)$, $\hat{\xb}_2(t) = (1/\epsilon^2) \xb_2(\epsilon t)$, and $\hat{\xb}_3(t) = (1/\epsilon^2) \xb_3(\epsilon t)$.  Let $\widehat{\ll}$ be the generator of the original system, one can also verify that $\widehat{\ll} = \epsilon \ll$. It is worth noting that the scaling with $\epsilon>0$ does not affect the topological stability. \Qed
    \iffalse
    \begin{equation*}
        \xb_1(t) = \frac{1}{\epsilon}  \xb_1(\epsilon t), \;\xb_2(t) = \frac{1}{\epsilon^2} \xb_2(\epsilon t), \;\xb_3(t) = \frac{1}{\epsilon^2} \xb_3(\epsilon t). 
    \end{equation*}\fi
\end{rem}

\ym{We choose to use the scaled system to facilitate the prediction of long-term errors in the original unscaled chaotic system within a smaller region of interest and over a short dimensionless observation horizon. Specifically, due to the nature of chaotic systems,} two trajectories starting very close together will rapidly diverge from each other, resulting in completely different long-term behaviors. The practical implication is that long-term prediction becomes impossible in such a system, where small errors are amplified extremely quickly. 
Such a two-point motion divergence phenomenon can be roughly characterized by the maximal Lyapunov exponent $\varrho>0$ \cite{sato1987practical} and quantitatively expressed by $\mathcal{E}(t)\approx \mathcal{E}(0)e^{\varrho t}$, where the process $\set{\mathcal{E}(t)}_{t\geq 0}$ represents the evolution of error w.r.t. any initial condition. Suppose we now consider the original unscaled system (as studied in   \cite{mauroy2019koopman}), and denote $\widehat{\mathcal{E}}(t)$ as the unscaled error. Then, in view of Remark \ref{rem: scaling} on the relation between the solutions of the original and scaled systems, one can obtain $|\widehat{\mathcal{E}}(t/\epsilon)|\geq |\mathcal{E}(t)|/\epsilon$. \ym{%Predicting long-term error is key for characterizing chaotic systems. 
Learning the scaled system with the same sampling frequency $\gamma$ and simulating trajectories up to  time $T_s$ allows us to infer the unscaled system’s error at  $T_s/\epsilon$  (at least $1/\epsilon$ times larger) without actually running the unscaled system for that long. This technique also avoids misleading demonstrations where the original system is simulated only up to a seemingly large time horizon, before significant error growth occurs, giving a falsely optimistic impression of performance.}

%Hence, whatever error is reflected in the scaled system, the error for the unscaled version will be at least $1/\epsilon$ times larger at the extended time scale $t/\epsilon$. 

%%The above scaled system is considered to facilitate easier observation within a smaller region of interest and to maintain a lower growth rate, ensuring that the tuning parameters remain within a reasonable range similar to the previous case study. For the prediction of a chaotic system, controlling the error within an even smaller tolerance level is necessary, and in this regard, RTM significantly outperforms the other two methods.
 %We consider the above scaled system for easier observation within a smaller region of interest and to maintain the growth rate at a smaller value, ensuring that the tuning parameters remain relatively within a reasonable range  as in the previous case study. 
 
%%Fig. \ref{fig: lorenz_RTM} and \ref{fig: lorenz_KLM_FDM}

\begin{figure}
    \centering
    \includegraphics[width=0.85\linewidth, trim={6.3cm 2.3cm 4cm 1cm}, clip]
    {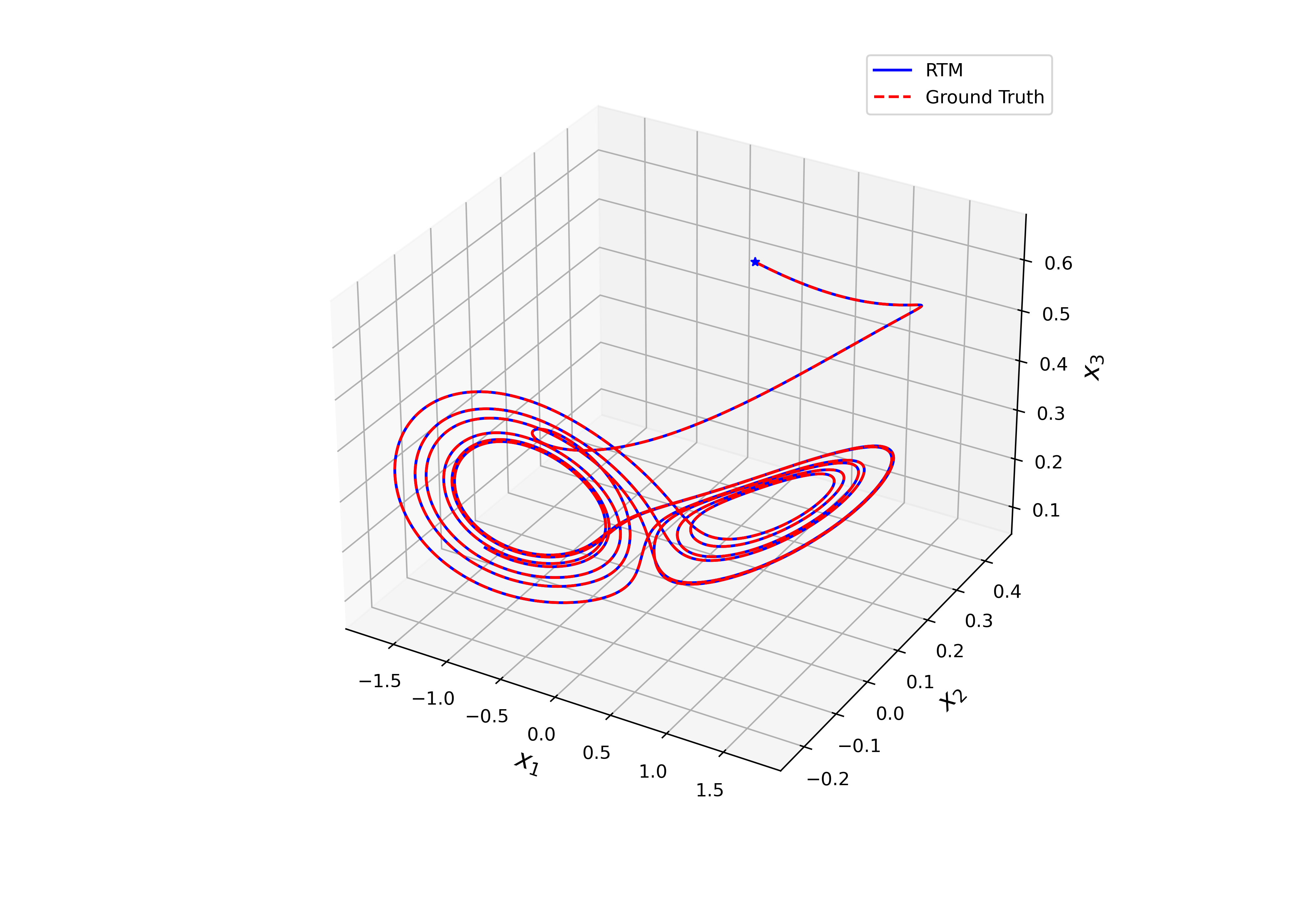}
    \caption{Comparison of the trajectory using RTM with the ground truth for the scaled Lorenz-63 system, where the blue star denotes the initial condition.}
    \label{fig: lorenz_RTM}
\end{figure}

\begin{figure}
    \centering
    \includegraphics[width=0.85\linewidth, trim={6.3cm 2.3cm 4cm 1cm}, clip]{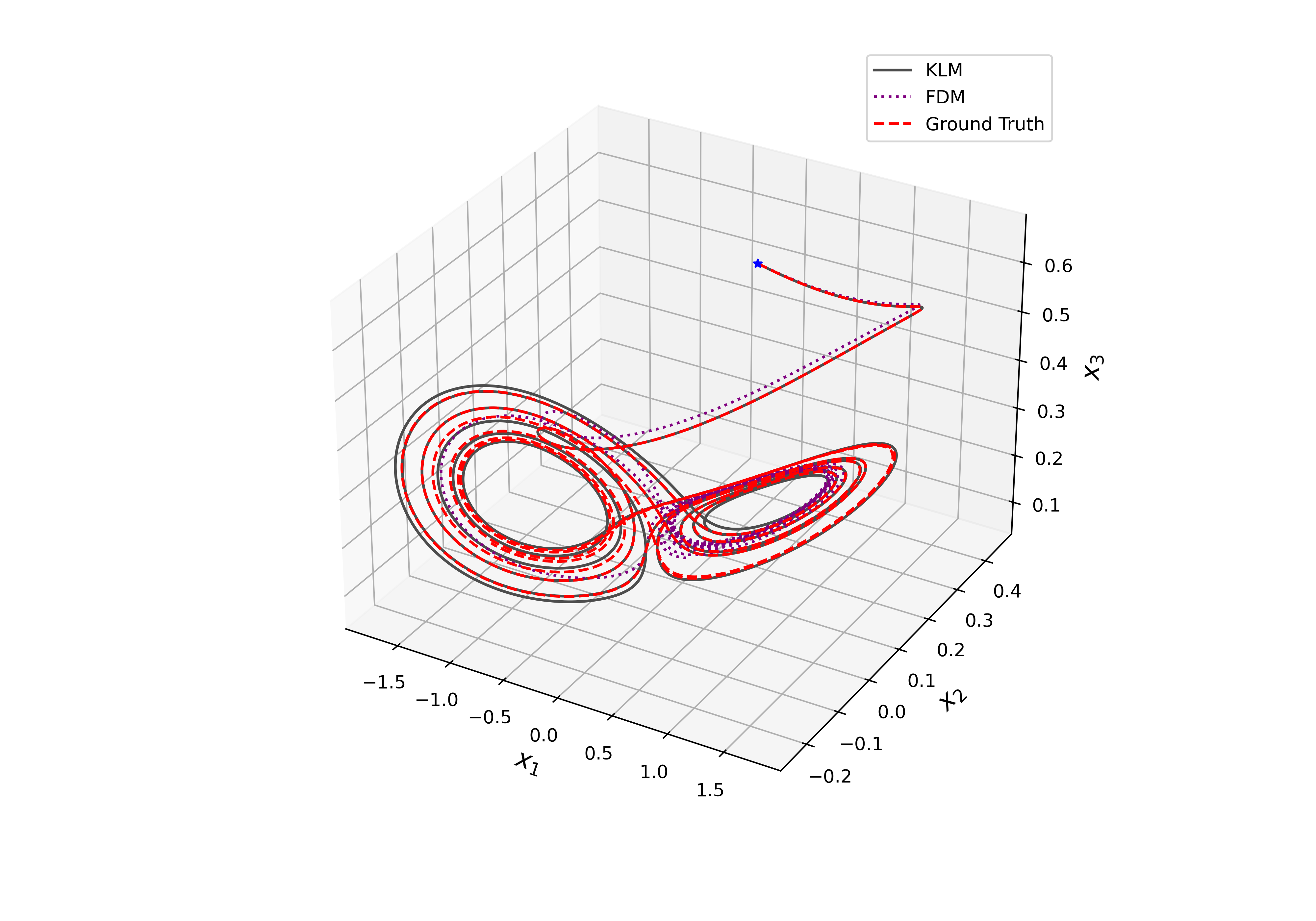}
    \caption{Comparison of the trajectories using KLM and FDM with the ground truth for the scaled Lorenz-63 system, where the blue star denotes the initial condition.}
    \label{fig: lorenz_KLM_FDM}
\end{figure}

In this case study, we set $P = Q=J = 2$ (or $N=2^3)$.  We also set $\mu = 2.5$, $\lambda = 1e8$, and $\taus = 1$ for RTM.  The comparisons for RMSE of weights and flow are summarized in Table~\ref{tab:three poly}. Fig.~\ref{fig: lorenz_RTM} depicts the comparison of the trajectory for $T_s=100$ using the approximated dynamics with RTM ($\gamma = 100$) and ground truth for an initial condition sampled in $\X$, while the comparisons for the ones with KLM and FDM are included in Fig.~\ref{fig: lorenz_KLM_FDM}. Given that the error grows exponentially w.r.t. the top Lyapunov exponent, it is evident that the flow prediction using the dynamics approximated by RTM is highly accurate, successfully exhibiting the attractor (a long-term behavior) for this chaotic system, while the other two methods struggle to predict the trajectory effectively.

\subsection{Polynomial Systems}
\ym{This subsection presents case studies on the system identification of polynomial systems, with a comprehensive comparison among Koopman-based methods and the well-matured SINDy approach. We revisit the scaled Lorenz–63 system and include two additional representative systems}, the reversed Van der Pol oscillator 
\begin{equation*}
\dot \xb_1 =-\xb_2, \quad 
\dot \xb_2 = \xb_1 - (1 - \xb_1^2)\xb_2
\end{equation*} 
and the (unscaled) Lorenz–96 system 
\begin{equation*}
    \dot{\xb}_j = - \xb_{j-2} \xb_{j-1} + \xb_{j-1} \xb_{j+1} - \xb_j+ 0.1, \;j\in\{0, 1, \cdots, 6\}.
\end{equation*}
\ym{
To ensure a fair comparison, we also present a sparse variant of RTM,
referred to as Sparse RTM (SRTM), by enforcing sparsity in the same manner as in SINDy,
namely, using sequential thresholded least squares on a given validation dataset to promote
sparsity in the learned weights. }

We investigate the performance of the methods under $\gamma = 10, 50, 100$, set $\X = (-1, 1)^d$, and use monomial dictionary functions. For RTM and SRTM, we set $\taus = 1$, $\lambda = 1e8$ and set $\mu=2.5$ for all three sampling rates $\gamma$. Detailed comparisons with the other two methods are provided in Table~\ref{tab:three poly}.

For all polynomial systems considered, the results show that RTM and SRTM achieve higher accuracy than KLM and FDM in both weight estimation and flow prediction. Moreover, all three Koopman-based methods improve as the sampling rate increases, with FDM exhibiting the slowest improvement. \ym{Regarding SINDy, the results indicate that although it performs well for low-dimensional polynomial systems, Koopman-based methods, especially RTM/SRTM, consistently outperform it in high-dimensional nonlinear systems such as Lorenz–96, providing superior accuracy in both weight estimation and flow prediction.}

\begin{table*}[htbp]
\centering
\caption{Comparisons of RMSE of weights and flow over 100 trajectories for the three polynomial systems}
\begin{adjustbox}{max width=\textwidth, scale=1}
\begin{tabular}{ccc|c|ccccc|ccccc}
\hline
\hline
\multirow{2}{*}{System} & \multirow{2}{*}{$M$} & \multirow{2}{*}{$N$} & \multirow{2}{*}{\makecell{$\gamma$}} & \multicolumn{5}{c|}{RMSE of weights ($\mathcal{E}_{\operatorname{RMSE}}^{\text{W}}$) }    & \multicolumn{5}{c}{RMSE of flow ($\mathcal{E}_{\operatorname{RMSE}}^{\text{F}}$) }      
\\ \cline{5-14} & & &  & SINDy & FDM & KLM & RTM & SRTM & SINDy & FDM & KLM & RTM & SRTM \\ \hline
\multirow{3}{*}{Reversed Van der Pol}     & \multirow{3}{*}{$10^2$}   & \multirow{3}{*}{9}    &   10 & 1.26e-4 &  3.52e-2  &  5.53e-4  &  \textbf{3.43e-5} & \textbf{3.20e-5} &  2.47e-5 &  1.96e-2  & 9.51e-5  &  \textbf{2.04e-5} & \textbf{1.85e-5} \\ 
\cline{4-14}  &  & & 50 &  1.15e-7 &  7.08e-3  & 2.22e-5   &  \textbf{1.88e-8} & \textbf{1.22e-8} & 2.56e-8 & 4.01e-3    &  3.65e-6 & \textbf{7.27e-9} & \textbf{5.51e-9 } \\ 
\cline{4-14}  &  & & 100 & 6.53e-9 &  3.54e-3  & 6.41e-6   &  \textbf{6.23e-9} & \textbf{5.54e-9} & \textbf{1.54e-9} & 1.83e-3    &  1.07e-6 & 3.37e-9 & 2.34e-9 \\ \hline
\hline
\multirow{3}{*}{\makecell{Scaled Lorenz-63 \\ (3-dimensional)} }   & \multirow{3}{*}{$10^3$}  & \multirow{3}{*}{8}
  & 10  & \textbf{2.45e-4} &   1.82e-1  & 6.99e-3  &  1.74e-3 & 1.43e-3  & \textbf{1.05e-3} & 2.32e-1 &  1.12e-2  &  1.47e-3 &  1.18e-3 \\ 
% \cline{4-10} & &  20   &  2.87  & 2.96   & \textbf{1.52} ($\lambda$ = 30)  &  1.47  & 1.52  &  \textbf{0.78}  \\ 
\cline{4-14} & & &  50  & 3.92e-7 &   3.81e-2  & 2.85e-4   &  \textbf{3.67e-7} & \textbf{1.76e-7}  & 1.43e-6 & 4.32e-2  & 4.78e-4 &  \textbf{3.12e-7} & \textbf{2.32e-7} \\ 
\cline{4-14} & & &  100  & \textbf{2.49e-8} &  1.92e-2  &  7.11e-5 &  \textbf{2.49e-8} & \textbf{1.52e-8} & 8.62e-8 & 2.31e-2  &  1.18e-4   &   \textbf{4.42e-8} & \textbf{2.48e-8} \\ \hline
\hline
\multirow{3}{*}{\makecell{Lorenz-96 \\ (6-dimensional)}}   & \multirow{3}{*}{$5^6$}    & \multirow{3}{*}{64}   
&  10  & 3.03e-1 &  2.18e-2  &  3.03e-4  & \textbf{7.51e-5} &\textbf{4.01e-5} & 3.31e-1 & 2.96e-2  & 3.29e-4 & \textbf{6.40e-5} & \textbf{5.63e-5} \\ 
\cline{4-14} & & & 50  & 2.74e-1 &  4.63e-3 &  1.17e-5   &  \textbf{1.84e-8} & \textbf{6.63e-9}   & 3.60e-1 & 6.08e-3   & 1.42e-5 &  \textbf{1.71e-8} & \textbf{1.58e-8}  \\ 
\cline{4-14} & & &  100  & 2.89e-1 &  2.33e-3  & 2.92e-6  &  \textbf{3.64e-9} & \textbf{3.41e-9}  & 3.74e-1 & 3.05e-3  &  3.60e-6 &   \textbf{4.69e-9} & \textbf{4.68e-9} \\  \hline
\hline
\end{tabular}
\end{adjustbox}
\label{tab:three poly}
\end{table*}
 
\subsection{$7$-D Biochemical System}
\ym{
In this subsection, we investigate a $7$-dimensional biochemical system modeling yeast glycolytic oscillations \cite{ruoff2003temperature}, which has become a standard benchmark for system identification tasks. It was examined in \cite[Appendix B]{brunton2016discovering} as a failure case for SINDy and its variant, Weak SINDy (WSINDy). Here, we slightly modify the 
$\dot{S}_1$ term to ensure that the state space remains invariant and satisfies the assumptions required by the Koopman framework:
}
\begin{small}
    \begin{equation*}
\label{eq:yeast_glycolysis_mod}
\begin{aligned}
\dot S_{1} &= k_{\rm ex}\,(G_{\rm ex} - S_{1}) \;-\; k_{1}\,S_{1}\,S_{6}/(1+(S_{6}/K_{1}\bigr)^{q}),\\
\dot S_{2} &= 2\,v_{1} \;-\; k_{2}S_{2}(N - S_{5}) \;-\; k_{6}S_{2}S_{5},\\
\dot S_{3} &= k_{2}S_{2}(N - S_{5}) \;-\; k_{3}S_{3}(A - S_{6}),\\
\dot S_{4} &= k_{3}S_{3}(A - S_{6}) \;-\; k_{4}S_{4}S_{5} \;-\; \kappa\,(S_{4} - S_{7}),\\
\dot S_{5} &= k_{2}S_{2}(N - S_{5}) \;-\; k_{4}S_{4}S_{5} \;-\; k_{6}S_{2}S_{5},\\
\dot S_{6} &= -2\,v_{1} \;+\; 2\,k_{2}S_{3}(A - S_{6}) \;-\; k_{5}S_{6},\\
\dot S_{7} &= \psi\,\kappa\,(S_{4} - S_{7}) \;-\; k\,S_{7}.
\end{aligned}
\end{equation*}
\end{small}

\ym{The values of the parameters are reported in Table \ref{tab:yeast_glycolysis_params_mod}. In this case, we use monomials up to order $2$, resulting in $N=36$ dictionary functions. The list of the monomials can be found in \cite[Appendix B]{brunton2016discovering}. We salso ample $M=7^7$ initial points randomly across $\X=(0, 0.5)^7$. The parameter $\mu$ is set to 2, and $\lambda = 1e8$ for RTM. The results are summarized in Table~\ref{tab:biochem}. }

\begin{table}[htbp]
  \centering
  \caption{\ym{Yeast glycolysis model parameters for the modified system}}
  \label{tab:yeast_glycolysis_params_mod}
  \begin{tabular}{@{} l c || l c || l c @{}}
    \toprule
    \toprule
    Parameter & Value & Parameter & Value & Parameter & Value \\
     
    $k_{\rm ex}$   & 0.5   & $G_{\rm ex}$ & 0.5   & $k_{1}$     & 100   \\
    $k_{2}$       & 6     & $k_{3}$     & 16    & $k_{4}$     & 100   \\
    $k_{5}$       & 1.28  & $k_{6}$     & 12    & $k$         & 1.8   \\
    $\kappa$     & 13    & $q$         & 4     & $K_{1}$     & 0.52  \\
    $\psi$       & 0.1   & $N$         & 1.0   & $A$         & 4.0   \\
    \bottomrule
    \bottomrule
  \end{tabular}
\end{table}

\begin{table}[h!bt]
\centering
\caption{\ym{Comparisons of RMSE of flow over 100 trajectories for the 7-dimensional biochemical system using monomials}}
\ym{
\begin{tabular}{ccccccc}
\hline
\hline
\multirow{2}{*}{$\gamma$} & \multicolumn{6}{c}{RMSE of flow ($\mathcal{E}_{\operatorname{RMSE}}^{\text{F}}$)} \\ 
\cline{2-7}
& SINDy & WSINDy & FDM & \makecell{KLM \\ (Re)} & \makecell{KLM \\ (Im)} & RTM \\
\hline
10  & - & - & 8.69e-2 & \textbf{2.63e-2} & 9.39e2 & 3.35e-2 \\
\cline{1-7}
50  & - & - & 3.78e-2 & 3.31e-2 & 7.58e-1 & \textbf{2.13e-2} \\
\cline{1-7}
100 & 1.16e4 & 2.02e2 & 3.35e-2 & 4.18e-2 & 7.58e-1 & \textbf{2.33e-2} \\
\hline
\hline
\end{tabular}
}
\label{tab:biochem}
\end{table}

\ym{RTM overall outperforms the other four methods in terms of the RMSE of flow. KLM produces imaginary parts, a phenomenon discussed in Section~\ref{sec: comparison}, even though the error appears small if only the real part is considered.  SINDy/WSINDy   fail to produce meaningful regression results for this high-dimensional system. The weights learned by SINDy lead to numerical instability and computational blow-up, consistent with the limitations reported in the original literature.
}

\subsection{Rational Vector Fields}\label{sec: non-poly}
Consider a system with nonpolynomial vector fields:
$$
\dot{\xb}_1 = -\xb_1 + \frac{4 \xb_2}{1 + \xb_2^2}, \quad
\dot{\xb}_2 = -\xb_2 - \frac{4 \xb_1}{1 + \xb_2^2}.
$$

\subsubsection{Monomial dictionary functions} We first choose monomials $\prod_{i=1}^dx_i^{\alpha_i}$ as the dictionary functions and evaluate their performance, highlighting that even when a sparse representation with monomial dictionary functions fails, the proposed RTM still yields lower errors. The parameter values for RTM are the same as those used in the previous cases. Unlike polynomial systems, where a relatively low order of monomials is sufficient, non-polynomial systems often require a larger $N$ for better approximation accuracy. However, high-order polynomial regression can lead to numerical instability or even computational blow-up due to regression error. To address this dilemma, we employ two monomial sets for the case studies: the first set contains $N=10$ monomials with \begin{small}
    $\sum_{i=1}^d  {\alpha_i} \leq  3$,
\end{small} while the second set includes $N=15$ monomials with \begin{small}
    $\sum_{i=1}^d  {\alpha_i} \leq  4$
\end{small}. We   randomly sample $M=10^2$ initial points  within $\X=(-1, 1)^2$. We also set  $\taus = 1$, $\mu=3.5$, and  $\lambda=1e8$ for RTM.  The results are summarized in Table~\ref{tab:nonpoly_poly}. It is evident that RTM overall outperforms the other four methods in terms of RMSE of flow.  

\begin{table}[h!bt]
\centering
\caption{Comparisons of RMSE of flow over 100 trajectories for the rational dynamical system using monomials}
\ym{
\begingroup
\setlength{\tabcolsep}{4pt} % default is ~6pt; smaller = narrower columns
\renewcommand{\arraystretch}{1.1} % a touch more vertical space for readability
\begin{tabular}{@{}cc*{6}{c}@{}}
\hline
\hline
\multirow{2}{*}{$N$} & \multirow{2}{*}{$\gamma$} & \multicolumn{6}{c}{RMSE of flow ($\mathcal{E}_{\operatorname{RMSE}}^{\text{F}}$)} \\ 
\cline{3-8}
&  & SINDy & WSINDy & FDM & \makecell{KLM \\ (Re)} & \makecell{KLM \\ (Im)} & RTM \\ 
\hline
\multirow{3}{*}{10} & 10  & 1.22e-1 & 3.14e-2 & 1.12e-1 & 1.41e-2 & 7.86e-1 & \textbf{1.39e-2} \\ 
\cline{2-8}
& 50 & 2.76e-1 & 1.85e-2 & 3.00e-2 & 1.41e-2 & 7.77e-1   & \textbf{1.39e-2} \\
\cline{2-8}
& 100 & 3.28e-1 & 2.05e-2 & 2.01e-2 & \textbf{1.43e-2} & 7.77e-1  & 1.57e-2 \\ 
\hline
\multirow{3}{*}{15} & 10  & 1.05e-1 & 4.73e-2 & 1.08e-1 & 1.37e-2 & 7.43e-1 & \textbf{1.36e-2} \\ 
\cline{2-8}
& 50 & 4.45e-1 & 6.78e-2 & 3.10e-2 & 1.49e-2 & 7.93e-1 & \textbf{1.46e-2} \\
\cline{2-8}
& 100 & 4.90e-1 & 1.35e-1 & 2.02e-2 & 1.44e-2 & 7.43e-1 & \textbf{1.36e-2} \\ 
\hline
\hline
\end{tabular}
\endgroup}
\label{tab:nonpoly_poly}
\end{table}

\subsubsection{Random $\operatorname{tanh}$-feature dictionary}

For this nonpolynomial system, as shown in~\cite{zeng2024sampling}, the previously mentioned dilemma of inefficient dictionary functions and potential regression error blow-up caused by high-order monomials can be mitigated by including terms of the form $\{ \frac{x_i}{1+x_j^p} \}$ in the dictionary,   which enables exact model identification. However, this prior knowledge is typically not accessible for most nonlinear systems in practice. To reduce the bias in selecting dictionary functions,  we utilize the
less biased hidden layers of random feature neural networks
as the dictionary functions.

\begin{rem}[Related Work]
In the context of Koopman operator learning, a Koopman autoencoder neural network has been proposed to improve approximation performance, seeking $\Kb = \operatorname{argmin}_{A\in\C^{N\times N}}\|Z-\Xb A\|$. However, for deep neural network dictionaries, a closed-form solution is not available, and training can be computationally expensive. In contrast, it has been shown that shallow $\text{tanh}$-activated neural networks can approximate functions as well as, or better than, deeper ReLU networks, with error decreasing as the number of hidden neurons increases. Their expressive power has also been demonstrated in solving first-order linear PDEs. Motivated by these results, we leverage this expressive feature by using shallow $\tanh$-activated neural networks as the dictionary functions, achieving the linear combination effect of the test and image functions of the learned operator.     \Qed
    %In the context of Koopman operator learning, the authors of ~\cite{azencot2020forecasting} proposed a consistent Koopman autoencoder neural network structure aimed at achieving better approximation performance. Essentially, this autoencoder seeks to find the $\Kb = \operatorname{argmin}_{A\in\C^{N\times N}}\|Z-\Xb A\|_F$ by optimally matching the input and output data. However, for deep neural network-based dictionary functions, a closed-form expression cannot be obtained, as they are not linear combinations. Additionally, the training process tends to be computationally expensive.

    %Comparatively, the work in \cite{de2021approximation} demonstrates that shallow $\tanh$ neural networks suffice to approximate functions at comparable or better rates than much deeper ReLU neural networks. The approximation error decreases when the number of hidden neurons increases~\cite{de2021approximation}. Numerous applications have also been explored in solving first-order linear PDEs, highlighting their expressive power. 
    
   % In the upcoming experiments, we leverage this expressive feature by using shallow $\tanh$-activated neural networks as the dictionary functions, achieving the linear combination effect of the test and image functions of the learned operator.
\end{rem}

By choosing $\tanh$  as the activation function, the dictionary of observable functions mainly consists of $\tanh(xW^\trans+b^\trans)$, where $W\in\R^{\sigma\times d}$ and $b\in\R^{\sigma}$ are the randomly generated weight and bias of the first layer, respectively, and $\tanh$ is applied elementwise.  %The tanh neural networks are widely used to approximate nonlinear functions in various applications, and the approximation error decreases when the number of hidden neurons increases~\cite{de2021approximation}. 
To obtain the expressions of the vector fields, we also need to add $\set{x_j}_{j=1}^d$ to the dictionary (recall that $\ll x_j = f_j$). 
All together, the dictionary $Z_N(x)$ consists of $\sigma$ random feature neural networks appended by $\set{x_j}_{j=1}^d$, with a total of $N=\sigma+d$ elements. 

%To get the expressions of the vector fields, we use $m$ random feature NN and $\xb$ as the dictionary functions, i.e., $N = m+n.$

We use the same parameter settings for RTM as in the cases studied with the monomial dictionary functions. As we cannot compare the RMSE of weights in this case, the comparisons of RMSE of flow with three difference frequencies are summarized in Table~\ref{tab:merged_systems}. It is worth noting that for $N = 52$ and $\gamma = 10$, KLM fails to produce meaningful results, as the trajectory of the approximated dynamics diverges. This is again attributed to the fact that taking the logarithm of the matrix $K$ (as discussed in Section~\ref{sec: comparison}) leads to undesirable results. On the other hand, RTM consistently outperforms the other methods overall. %Like polynomial systems, RTM outperforms the other two methods. 
\begin{table}[ht]
\centering
\caption{CComparison of RMSE of flow over 100 trajectories for two benchmark systems using $\tanh$-activated neural networks.}
\ym{
\begin{tabular}{ccccccc}
\hline
\hline
\multirow{2}{*}{System} & \multirow{2}{*}{$N$} & \multirow{2}{*}{$\gamma$} & \multicolumn{4}{c}{RMSE of flow ($\mathcal{E}_{\operatorname{RMSE}}^{\text{F}}$)} \\ 
\cline{4-7} & & & FDM & \makecell{KLM \\ (Re)} & \makecell{KLM \\ (Im)} & RTM \\
\hline

% ---------- Rational / Nonpoly system ----------
\multirow{9}{*}{Rational}
  & \multirow{3}{*}{22} & 10  & 7.27e-2 & 3.08e-3 & 7.78e-1 & \textbf{2.64e-3}\\
  \cline{3-7} & & 50  & 1.81e-2 & 2.47e-3 & 7.78e-1 & \textbf{2.22e-3} \\
  \cline{3-7} & & 100 & 8.90e-3 & 2.85e-3 & 7.80e-1 & \textbf{2.52e-3} \\
  \cline{2-7}
  & \multirow{3}{*}{52} & 10  & 1.08e-1 & -        & 1.78e36  & \textbf{6.72e-2} \\
  \cline{3-7} & & 50  & 2.42e-2 & 3.96e-2 & 7.90e-1  & \textbf{9.37e-4} \\
  \cline{3-7} & & 100 & 1.74e-2 & 1.10e-1 & 8.32e-1  & \textbf{9.19e-4} \\
  \cline{2-7}
  & \multirow{3}{*}{102} & 10  & 1.01e-1 & 1.43     & 2.60e39 & \textbf{5.48e-3} \\
  \cline{3-7} & & 50  & 2.64e-2 & 9.05e-2 & 7.40e-1 & \textbf{9.33e-5} \\
  \cline{3-7} & & 100 & 1.68e-2 & 2.12e-1 & 7.75e-1 & \textbf{6.18e-4} \\
\hline

% ---------- Two-machine system ----------
\multirow{9}{*}{Two-machine}
  & \multirow{3}{*}{22} & 10  & 1.13e-2 & \textbf{1.00e-3} & 7.80e-1  & 1.06e-3  \\
  \cline{3-7} & & 50  & 2.56e-3 & \textbf{1.01e-3} & 7.84e-1 & 1.09e-3 \\
  \cline{3-7} & & 100 & 1.70e-3 & \textbf{1.03e-3} & 7.82e-1 & 1.05e-3 \\
  \cline{2-7}
  & \multirow{3}{*}{52} & 10  & 1.10e-2 & 3.77e-3 & 7.76e-1 & \textbf{1.47e-4}\\
  \cline{3-7} & & 50  & 4.60e-3 & 2.32e-2 & 7.63e-1 & \textbf{1.07e-5} \\
  \cline{3-7} & & 100 & 7.75e-3 & 1.11e-2 & 7.63e-1  & \textbf{1.07e-5} \\
  \cline{2-7}
  & \multirow{3}{*}{102} & 10  & 1.18e-2 & 1.16e-2 & 7.73e-1 & \textbf{1.56e-4} \\
  \cline{3-7} & & 50  & 5.13e-3 & 8.36e-2 & 8.11e-1  & \textbf{8.08e-6} \\
  \cline{3-7} & & 100 & 2.12e-2 & 1.83e-1 & 8.04e-1 & \textbf{6.64e-6} \\
\hline
\hline
\end{tabular}}
\label{tab:merged_systems}
\end{table}
\subsection{Two-Machine Power System}
\label{sec: two-machine}

%For non-polynomial systems, selecting monomial dictionary functions may be ineffective due to their limited expressiveness. To enhance adaptability in characterizing the generator for more general dynamical systems, we will continue to investigate the impact of using non-polynomial dictionary functions for all three methods in the upcoming experiments.
Consider the two-machine power system \cite{vannelli1985maximal} modeled by $$\dot{\xb}_1 = \xb_2, \;\dot{\xb}_2 =  -0.5\xb_2 - (\sin(\xb_1 +\pi/3)-\sin(\pi/3)). $$
Since the vector field consists of polynomial and periodic terms, using monomial approximations appears to be local and may not be satisfactory. We continue to use random feature NNs as the dictionary functions. %Regrading training data, $M=100^2$ initial points are uniformly sampled within $\X=(-1, 1)^2$. We also set   $\tau_s = 5$ and  $\lambda=1e8$ for RTM.

%Instead, we utilize the less biased hidden layers of random feature neural networks as the dictionary functions. 

%To get the expressions of the vector fields, we use $m$ random feature NN and $\xb$ as the dictionary functions, i.e., $N = m+n.$

\begin{figure}[h!t]
    \centering
    \includegraphics[scale=0.35]{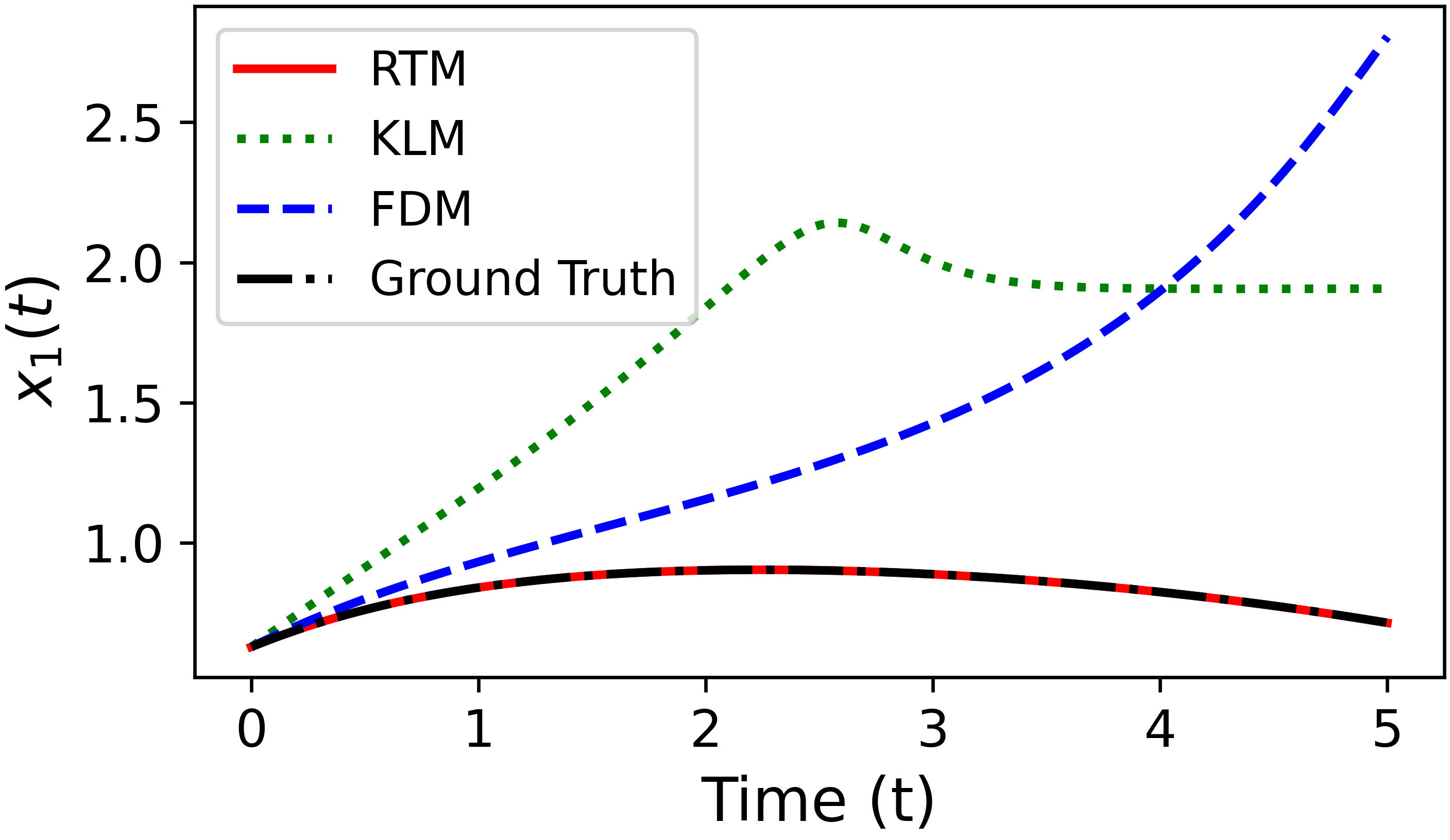}
    \includegraphics[scale=0.35]{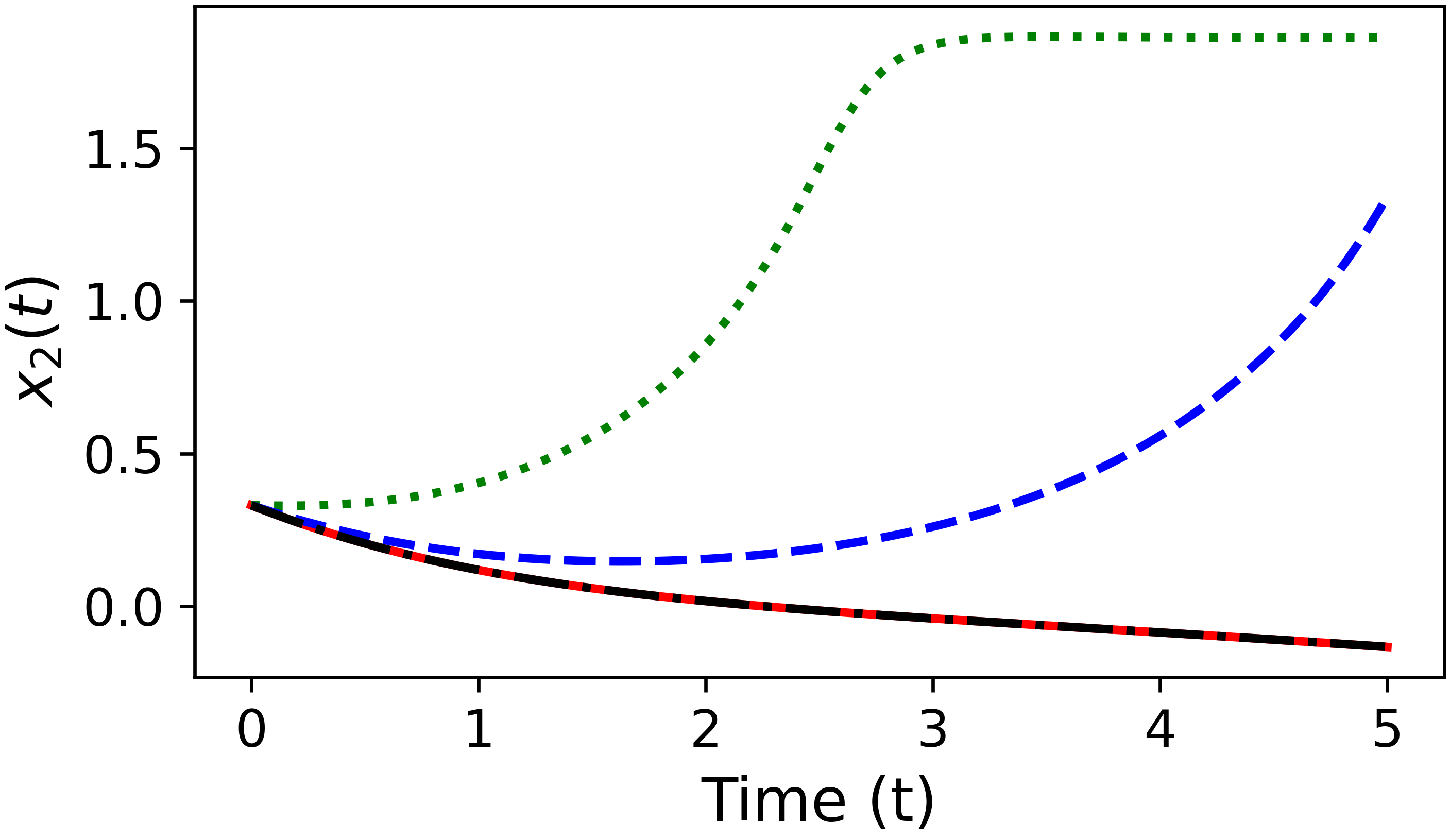}
    %\subfigure[$\sigma = 50, \tau = 0.01s$]{%
        %\includegraphics[scale=0.2]{trajectory_two_machine_comparisons_f_100_m=50.png}
        %\label{fig:tm-a}
    %}
    %\hfill
    % \subfigure[$\sigma = 100, \gamma = 100$ Hz]{%

    % }
    \caption{Comparisons of the trajectories with the approximated dynamics using three different methods and ground truth for the two-machine power system.}
    \label{fig:two-machine}
\end{figure}

We sampled $M  = 10^2$ initial points randomly across $\X=(-1, 1)^2$ with $\taus = 1$. For RTM, we set $\lambda = 1e8$ and $\mu = 2.5$ for all sampling frequencies. As shown in Table~\ref{tab:merged_systems}, RTM consistently outperforms FDM and KLM across nearly all configurations of the two-machine power system. Similar to the previous case study, when using $\tanh$ neural networks, taking the matrix logarithm of the learned Koopman matrix $K$ introduces non-negligible imaginary parts. The trajectories generated by all three methods remain close to the ground truth for $N = 22$ with $\gamma = 100$. However, when $N$ is increased to $102$, the performance of both FDM and KLM deteriorates (particularly KLM) as shown in Fig.~\ref{fig:two-machine}. These results indicate that when prior knowledge of the underlying system is lacking and monomial dictionary functions may not be suitable, which motivates the use of $\tanh$-activated random bases as dictionary functions, both FDM and KLM may struggle to produce accurate results in system identification, whereas RTM is more robust in this regard.

\subsection{Prediction of Region of Attraction}

In this subsection, we demonstrate that both the dataset and dictionary functions %generated during system identification 
can be reused to %solve Zubov's equation and 
predict the region of attraction (RoA) $\roa$ of an equilibrium point $\xeq$ for the system\footnote{Supposing  that $\xeq$ is asymptotically stable, the RoA of $\xeq$ is a set  defined as
$\roa:=\set{x\in\X: \lim_{t\ra\infty}|\phi(t, x)-\xeq|=0}.$}. For simplicity, we can assume that the system has an equilibrium point at the origin. 

 The methodology involves using the learned generator to solve a % construct a transformation of a maximal Lyapunov function based on 
 Zubov's equation of the form 
\begin{equation}\label{E: zubov}
     \ll u(x) = -\alpha|x-\xeq|^2 (1-u(x)),\;\;\alpha>0. 
\end{equation}
The solution should have the following properties: 1) $0<u(x)<1$ for all $x\in \roa\setminus\set{\xeq}$, and 2) $u(x)\ra 1$ as $x\ra y$ for any $y\in \partial \roa$. Therefore, the set $\{x\in\R^n: u(x)< 1\}$  can readily represent the RoA. Moreover, for any approximator $\widetilde{u}$ of $u$ such that $|\widetilde{u}-u|_\infty\leq \eps$, one can show that the set $\{x: \widetilde{u}(x)\leq 1-\eps\}$  is a tight inner approximation of the RoA.

Enlarging the RoA estimate is important for engineering applications as it yields a less conservative safe operational range. Traditional Lyapunov function methods select an ansatz from a finite-dimensional subspace of $\mathcal{C}^1$ satisfying $\mathcal{L}v(x) < 0$ for all $x \in \roa \setminus {\xeq}$, but such an ansatz often gives a conservative result. Zubov-related methods have shown advantages over other inequality-oriented approaches (e.g., \cite{topcu2008local, abate2021fossil, zhou2022neural}), and further details on these PDEs can be found in \cite{camilli2001generalization, liu2023physics}.

\iffalse
\begin{rem}
    Note that the RoA can also be characterized by the function domain $D$ of a maximal Lyapunov function $v$ if the following conditions hold: 1)   $v(\xeq)=0$ and $v(x)>0$ for all $x\in D\setminus \set{\xeq}$; 2) the derivative of $v$ along solutions of (\ref{E: sys}) 
     is well-defined for all $x\in D$ and satisfies the Lyapunov equation
          \begin{equation}\label{eq:DV1}
     \ll v(x)  = - |x-\xeq|^2;    
     \end{equation}
and  3) $v(x)\ra \infty$ as $x\ra \partial D$ or $|x|\ra \infty$. 

It can be verified that $u(x) = 1 - \exp(-\alpha v(x))$ for all $x\in\roa(\A)$, and $u(x)=1$ elsewhere, is the unique bounded viscosity solution that satisfies the required properties \cite{liu2023physics}. We prefer to use $u$ instead of $v$ for characterizing the RoA due to its global boundedness, which allows the extension of the function domain to the entire state space or any desired set where computations take place.  Additional features of the above PDEs can be found in \cite{camilli2001generalization, liu2023physics}. \Qed
\end{rem}\fi

To illustrate the accuracy of the learned generator, as well as the effectiveness of using it to solve Eq. \eqref{E: zubov}, we revisit the reversed Van der Pol oscillator system. The difference is that we choose the dictionary $Z_N(x)$ the same as in Sections \ref{sec: non-poly}, which consists of the %less biased, 
%more expressive 
$\sigma$-dimensional random feature neural network $\tanh(xW^\trans+b^\trans)$ appended with $\set{x_1, x_2}$. In the experiment, we set   $\sigma = 100$, $\mu=10$, and $M=100^2$. 

\begin{rem}
    Due to page limitations, we include only the standard benchmark, the reversed Van der Pol oscillator, which is   used as an example in virtually all Lyapunov function construction tools. Additional system comparisons and detailed discussions are provided in the supplementary material. \Qed %It is worth emphasizing that a method that performs well on another system may not necessarily be suitable for this seemingly “boring” Van der Pol system (see the detailed comparisons in the supplementary material). Moreover, there is no generally accepted standard, not even a partial ordering, for classifying systems as “easy” or “complicated” in terms of ROA estimation. Nevertheless, the proposed learning method demonstrates superior overall performance. 
\end{rem}

With the above settings,  we can achieve a good approximation of the dynamics with the RMSE of flow $\mathcal{E}_{\operatorname{RMSE}}^{\text{F}}$ $= 1.08e$$-5$.  %The approximate phase portrait is plotted in Fig. \ref{fig: phase portrait}. 
We find  an approximate equilibrium point at $[6.81e$$-8$, $-1.25e$$-8]$, which is sufficiently close to the origin with negligible error. In addition, we %\ym{uniformaly sample $100^2$ points on $\X$ and}
observe that the average approximation error of $\ll \tanh(xW^\trans+b^\trans)$  is $1.66e$$-4$, indicating an overall good performance. % of the least-squares regression $\Lb_\lambda$.  

\iffalse
\begin{figure}
    \centering
    \includegraphics[width=0.85\linewidth]{Phase_Portrait_with_m=100.png}
    \caption{Phase portrait of the approximated dynamics, where the black dot denotes the origin which is the equilibrium point of the true dynamics.}
    \label{fig: phase portrait}
\end{figure}
\fi

We aim to find a   single-hidden-layer feedforward neural network of the form $u(x;\theta) = Z_N(x)\theta$
%\begin{equation}\label{eq:elm}
 %V(x;\theta) =  \tanh(xW^\trans+b^\trans)\theta
%\end{equation}
 that approximates the unique bounded viscosity solution of $\ll u(x) = -\alpha|x|^2(1-u(x)$ with $\alpha = 0.1$ and the boundary condition $u(\zero)=0$ using the learned generator.  The problem is then reduced to finding $\theta$ such that $\Zk_N(x)\Lb_\lambda\theta = -\alpha|x|^2(1-    u(x;\theta))$.   %Due to the fact that $\nabla u(x;\theta)\cdot \widetilde{f}(x) \approx \ll u(x;\theta)\approx \Zk_N(x)\Lb_\lambda \theta$  for all $u(x;\theta) = \Zk_N(x)\theta$, where $\widetilde{f}(x)$ is the approximated vector field using the learned generator,  it is equivalent to finding the weights $\theta$ that approximately satisfy 
 %\begin{equation}\label{E: zubov_approx_1}
     %\Zk_N(x)\Lb_\lambda\theta = -\alpha|x|^2(1-    u(x;\theta)). 
 %\end{equation}
 %with the inherited boundary value information.
 \iffalse
 or 
 \begin{equation}\label{E: zubov_approx_2}
     \nabla u(x;\theta)\cdot \widetilde{f}(x) = -\alpha |x|^2(1-  u(x;\theta)),
 \end{equation}
 ensuring that $u(x; \theta)$ is bounded and a viscosity solution at the same time. \fi
 
 %It is worth noting that, without the boundary condition and the constraints on its viscosity property, the problem can be reduced to finding the $\theta$ that minimizes the least-squares error (residual loss)   between the l.h.s. and r.h.s. of Eq. \eqref{E: zubov_approx_1} (or \eqref{E: zubov_approx_2}) at those sampled initial conditions.  However, such a solution may not possess the correct physical meaning, with a simple counterexample being  $u(x;\theta)\approx 1$ for all $x$. 

Inspired by recent research on physics-informed neural networks, it is necessary to introduce and minimize additional loss terms beyond the residual loss that encompass the PDE and supplementary problem information \cite{liu2024tool, liu2023physics, zhou2024physics}. Here, we specifically include loss terms to match $u(\zero; \theta)=0$ and values on $\partial \X$. Clearly, finding $\theta$ that minimizes the weighted sum of these loss terms does not have a closed-form solution. That said, an extreme learning machine (ELM), which seeks solutions of the form $\Zk_N(x)\theta$ for PDEs that are linear in both $u$ and $\nabla u$, transforms the optimization into a linear least-squares problem that can be efficiently solved, exhibiting fast learning speed and strong generalization performance \cite{zhou2024physics}.

\begin{figure}
    \centering
    \includegraphics[width=0.95\linewidth]{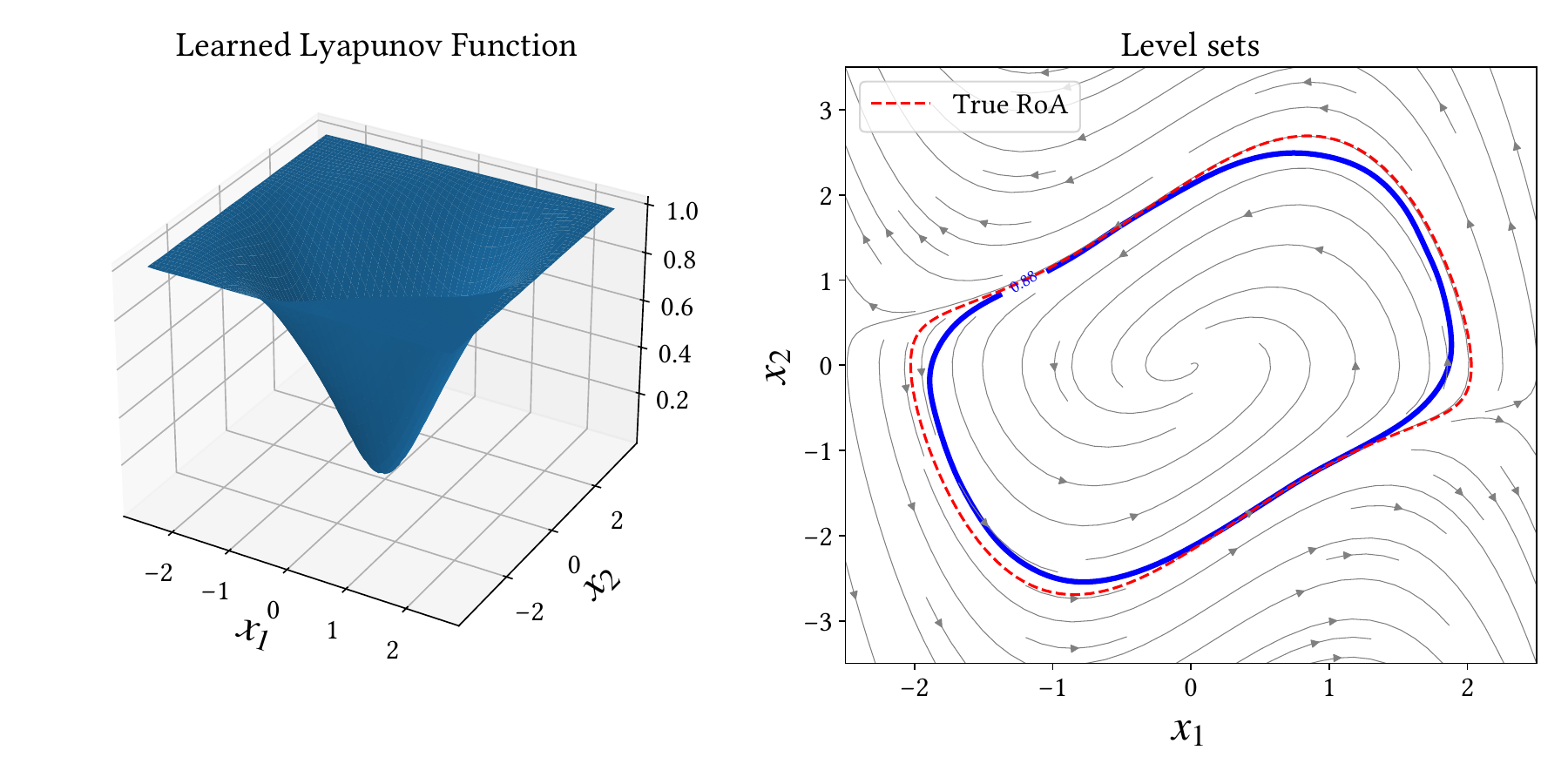}
     \caption{Learned Lyapunov function and the corresponding region of attraction estimate using the random feature neural network dictionary functions for the reversed Van der Pol oscillator.}
    \label{fig:zubov}
\end{figure}

We also encourage readers to refer to \cite{zhou2024learning} for a more systematic study on using learned generators to obtain stability certificates for unknown systems, with additional details on encoding valid physics-informed information for approximating 
$u(x)$ efficiently. The learned Lyapunov function and the corresponding region of attraction estimate can be found in Fig.~\ref{fig:zubov}, where the blue curve represents the region of attraction estimate, which is sufficiently close to the true region of attraction. This result in turn reflects the overall good approximation of the generator.  %Additionally,  the averaged approximation error over the sample points of $\ll u(x)$ is $1.33e$$-4$, which indicates that the proposed method is able to ... %We refer the reader to \cite{zhou2024physics} for more details about solving stability-related PDEs with ELM. 

\section{Conclusion}\label{sec: conclusion}
In this paper, we propose a novel resolvent-type Koopman operator-based learning framework for the generator of unknown systems. The proposed method demonstrates both theoretical and numerical improvements over the Koopman logarithm method \cite{mauroy2019koopman} and the finite difference method \cite{bramburger2024auxiliary}, offering greater adaptability in handling low observation rate scenarios and utilizing less biased random feature neural networks as dictionary functions. Specifically, in     cases where the logarithm of the Koopman operator is invalid for representing the generator, we draw upon the rich literature to propose %a sequence of 
data-driven approximations based on Yosida’s approximation and the properties of the resolvent operator. The analytical convergence   and a modified data-driven learning algorithm are provided to assist users in tuning the parameters. The extended work on control systems \cite{zeng2025data} based on this method also demonstrates better performance than benchmark tools.

The current drawbacks lie in tuning the parameter $\mu$ under the constraints of the observation rate, as well as in the lack of understanding regarding the sampling efficiency of the initial conditions. We will pursue these analyses in future work and provide a more systematic and automated parameter tuning strategy, with the ultimate goal of developing a software toolbox %designing computational tools 
for data-driven system verification of unknown systems. We also expect that the current research and the analysis of data efficiency will lay the foundations for studying control systems, shed light on dimension reduction in generator learning, and be beneficial for developing an online learning adaptation of this method with a provable error bound. Another possible extension direction based on the current research is robustness analysis under noisy and partial observations.

 Additionally, the domain can be extended to $L^\infty(\X)$, which allows us to study the $L^2$-adjoint operator of $\ll$ with domain $L^1(\X)$, leading to more interesting applications such as the Liouville equation and ergodic measure \cite{pavliotis2008multiscale}. We will pursue generator learning on this extended domain in future work.

\appendices

\section{Fundamental Properties of Koopman  Generators}\label{sec: facts}

\renewcommand{\thethm}{\thesection.\arabic{thm}}

\ym{
Note that the flow $\phi(t, x_0)$ is absolutely continuous. For learning the generator, we usually work with a (hypothetical) smoother surrogate of the flow. In particular, there exists a 
$\mathcal{C}^1(\Omega)$ surrogate that approximates the true flow uniformly on any finite time window. For discrete snapshot data, the approximation can be chosen to match the observations to arbitrary accuracy, and, if desired, exactly at finitely many sampled states via an additional small $\mathcal{C}^1$
 perturbation. By abuse of notation, we continue to   write 
$\phi(t,x_0)$ for this $\mathcal{C}^1$
 surrogate\footnote{We do not need to construct this surrogate explicitly. It is used solely to simplify exposition, and since the learned system is approximate in any case, this convention does not affect our conclusions.}.  Additionally, in this setting, for any $h\in \mathcal{C}^1(\X)$, $\nabla(\K_th)(x)$ exists and $\K_t: \mathcal{C}^1(\X)\ra \mathcal{C}^1(\X)$ for all $t>0$. %$\nabla(\K_th)(x) = D\phi(t, x)^\trans\nabla h(\phi(t, x))$

\iffalse
 \begin{rem}\label{rem: smooth_version}
    \Qed
\end{rem}
\fi

\begin{prop}\label{prop: core}
    Working with a sufficiently smooth surrogate $\phi(t, x)$, for any $h\in\dom(\ll)\subseteq\mathcal{C}(\Omega)$, there exists a sequence of $\{h_n\}\subseteq\mathcal{C}^1(\Omega)$ such that $\|h_n-h\|_{\sss\ll}\ra 0$. 
\end{prop}
\begin{pf}
For any $h\in\dom(\ll)$, let $\{g_n\}\subseteq\mathcal{C}^1(\Omega)$ be a sequence such that $\|g_n-h\| \ra 0$. 
Let $h_{n,t}(x) = \frac{1}{t}\int_0^t\K_t h_n(x) ds$ for any $t>0$. Then, $h_{n, t}\in\mathcal{C}^1(\Omega)$ for any $n$ and any fixed $t>0$. Consider $\{t_n\}$ such that $t_n\ra0$ and $\|g_n-h\|/t_n\ra 0$ as $n\ra\infty$, and set $h_n:=h_{n, t_n}$. Then, 
\begin{small}
    \begin{equation}
    \begin{split}
       & \|h_n-h\|\\ \leq &\left\|\frac{1}{t_n}\int_0^{t_n}\K_s(g_n-h)\dd s\right\| + \left\|\frac{1}{t_n}\int_0^{t_n}\K_s h\; \dd s-h\right\|\\
         \leq &\frac{C(e^{\omega t_n}-1)}{\omega t_n}\|g_n-h\| + \left\|\frac{1}{t_n}\int_0^{t_n}\K_s h\; \dd s-h\right\|, 
    \end{split}
\end{equation}
\end{small}

\noindent and $\|h_n-h\| \ra 0$ as $n\ra \infty$. 

Furthermore, by  \cite[Theorem 2.4, Chap 1]{pazy2012semigroups}, we have $\ll\int_0^t\K_sg_nds = \K_sg_n-g_n$. Therefore, 
\begin{small}
    \begin{equation}
    \begin{split}
       & \|\ll h_n-\ll h\| \\ = &\left\|\frac{\K_{t_n}g_n-g_n}{t_n}-\ll h\right\|\\
          \leq & \left\|\frac{(\K_{t_n}-\id)(g_n-h)}{t_n} \right\|+ \left\|\frac{\K_{t_n}h-h}{t_n}-\ll h\right\|\\
          \leq &\frac{Ce^{\omega t_n}+1}{t_n}\|g_n-h\|+ \left\|\frac{\K_{t_n}h-h}{t_n}-\ll h\right\|,
    \end{split}
\end{equation}
\end{small}

\noindent and $ \|\ll h_n-\ll h\|\ra0$ as $n\ra\infty$. The conclusion follows immediately. 
\end{pf}
}

\section{Proofs and Supplementaries for Section \ref{sec: finite-dim}}\label{sec: proof_characterize}

\iffalse
\begin{eg}
    Consider the simple dynamical system $\dot{\xb}= \xb$ and the observable function $v(x)=x^n$ for any $n\geq 1$. Then, analytically, $\phi(s,x)=xe^s$ for all $s\geq 0$ and $\ll v(x) = nx^{n}$. We test the validity of using Eq. \eqref{E: approx_t}. Note that, for sufficiently large $\lambda$, we have  $\lambda^2\int_0^\taus  e^{-\lambda s}(\K_sv(x))ds  = \lambda^2 \int_0^\taus e^{-\lambda s}e^{sn}x^n ds =\lambda^2x^n \int_0^\taus e^{-(\lambda-n) s}ds=  \frac{\lambda^2 x^n}{\lambda -n}(1-e^{-(\lambda -n)\taus})\approx \frac{\lambda^2 x^n}{\lambda-n}$, and $ \lambda^2\rr_{\lambda,\taus} v(x)-\lambda v(x)\approx \frac{\lambda^2 x^n}{\lambda-n}-\lambda x^n =\frac{n\lambda}{\lambda-n}x^n\approx nx^n=\ll v(x).$
    \iffalse
    \begin{equation*}
        \begin{split}
            &\lambda^2\int_0^\taus  e^{-\lambda s}(\K_sV(x))ds  \\= &\lambda^2 \int_0^\taus e^{-\lambda s}e^{sn}x^n ds =\lambda^2x^n \int_0^\taus e^{-(\lambda-n) s}ds\\
            = & \frac{\lambda^2 x^n}{\lambda -n}(1-e^{-(\lambda -n)\taus})\approx \frac{\lambda^2 x^n}{\lambda-n}
        \end{split}
    \end{equation*}
and 
\begin{equation*}
    \begin{split}
       & \lambda^2 \rr_{\lambda,\taus}  V(x)-\lambda V(x)\\
       \approx &\frac{\lambda^2 x^n}{\lambda-n}-\lambda x^n =\frac{n\lambda}{\lambda-n}x^n\approx nx^n=\ll V(x).
    \end{split}
\end{equation*}\fi
With high-accuracy evaluation of the truncated integral, we can achieve a reasonably good approximation. \Qed
\end{eg}\fi

\textbf{Proof of Proposition \ref{prop: compact_resolvent}:}
    %Note that the identity operator $\id$ on $C(\X)$ cannot be compact. Therefore, $s=0$ must be excluded. 
    We assume that $\rr_{\lambda,\taus}$ is compact for $\lambda>\omega$. By \cite[Theorem 3.2]{pazy2012semigroups}, $\K_s$ is continuous w.r.t.  the uniform operator norm $\|\cdot\|$. We can also easily verify the compactness of $\lambda\rr_{\lambda,\taus}\K_s$, for any $s\in(0, \taus]$,  by Definition \ref{def: T_t}. In addition, for any arbitrarily small $t>0$, 
    \begin{small}
            \begin{equation*}
        \begin{split}
            & \|\lambda \rr_{\lambda,\taus}\K_s-\K_s\|\\
            \leq & \int_0^t \lambda e^{-\lambda \sigma}\|\K_{s+\sigma}-\K_s\|d\sigma + \int_t^\taus \lambda e^{-\lambda \sigma}\|\K_{s+\sigma}-\K_s\|d\sigma\\
            \leq & \sup_{\sigma\in(0, t]}\|\K_{s+\sigma}-\K_s\| (1-e^{-\lambda t}) + 2\int_t^\taus \lambda e^{-\lambda\sigma}\|\K_\taus\|d\sigma\\
            \leq & \sup_{\sigma\in (0, t]}\|\K_{s+\sigma}-\K_s\| + \frac{2\Sigma\lambda}{\lambda-\omega}e^{\omega(t+\taus)-\lambda t}. 
        \end{split}
    \end{equation*}
    \end{small}

\noindent By the continuity of $\K_s$, letting $\lambda\ra\infty$ and $t\ra 0$, we have $\K_s\ra \lambda\rr_{\lambda,\taus}\K_s$ in the uniform sense, which shows the compactness. 

To show the converse side, we notice that the operator $\rr_{\lambda,\taus}^t:=\int_t^\tau e^{-\lambda s}\K_s ds$ is always compact for any $t\in(0, \taus]$ given the compactness of $\K_s$ with $s\in(0, \taus]$. However, 
\begin{small}
    \begin{equation*}
    \begin{split}
        \|\rr_{\lambda,\taus}-\rr_{\lambda,\taus}^t\|&\leq \int_0^t e^{-\lambda s}\|\K_s\|ds\leq \int_0^t \|\K_s\|ds\\
        & \leq \int_0^t \Sigma e^{\omega s}ds \leq t \Sigma e^{\omega t}. 
    \end{split}
\end{equation*}
\end{small}

\noindent Letting $t\ra 0$, we see the uniform convergence of $\rr_{\lambda,\taus}^t$ to $\rr_{\lambda,\taus}$, which shows  compactness of $\rr_{\lambda,\taus}$. \pfbox

\textbf{Proof of Corollary \ref{cor: finite_1}}:
    We show the sketch of the proof. Working on the compact family of $\set{\K_s^\eps}_{s\in(0, \taus]}$ for some small $\eps>0$, one can find a sufficiently large $N$ such that 
    \begin{small}
            \begin{equation*}
        \|\K_s^Nh - \K_sh\|<\delta, \;\;h\in C(\X), \;s\in(0, \taus], 
    \end{equation*}
    \end{small}
    
\noindent     where $K_s^N$ is the finite-dimensional representation of $K_s^\eps$. Let $\rr_{\lambda,\taus}^N:=\int_0^\taus e^{-\lambda s} \K_s^Nds$. Then, for any $t\in(0, \taus)$, 
    \begin{small}
            \begin{equation*}
        \begin{split}
           & \|\rr_{\lambda,\taus}^N  h-\rr_{\lambda,\taus}h\|\\
            \leq &\int_0^t e^{-\lambda s} \|\K_s^Nh-\K_s h\|  ds + \int_t^\taus e^{-\lambda s} \|\K_s^Nh-\K_s h\| ds\\
            \leq & C\|h\| t + \sup_{s\in(0, \taus]}\|\K_s^Nh - \K_sh\|\cdot\frac{e^{-\lambda t}-e^{-\lambda \taus}}{\lambda}\\
            < & C\|h\| t + \delta,
        \end{split}
    \end{equation*}
    \end{small}
    
\noindent where $C:=\sup_{s\in (0, t]}\|K_s^N\|+ \Sigma e^{\omega t}$. The conclusion follows by sending $t \ra 0$. \pfbox

\iffalse
\begin{deff}
    For any $\eps>0$, consider a partition $\X = \sqcup_{m=1}^M\X_m$ satisfying $\operatorname{diam}(\X_m)\leq \delta_M$ for all $m$  and $\delta_M\ra0$ as $M\ra\infty$
\end{deff}

For any $\eps>0$, we can construct a partition $\X = \sqcup_{m=1}^M\X_m$ with $\operatorname{diam}(\X_m)\leq \delta_M \ra0$ as $M\ra\infty$. For each $m$, pick a node $x^{(m)}\in\X_m$. Let $\psi_m^\eps(x) =\int_{\X_m}\eta_\eps(x-y)dy$. Then we can construct $\K_t^Mh(x)=  \sum_{m=1}^M \psi_m^\eps(x) \K_th(x^{(m)})$, and show that 
\begin{equation}
    \begin{split}
        & |\K_t^Mh(x)-\K_t^\eps h(x)|_\infty \\
        = & \left|\sum_{m=1}^M\left(\int_{\X_i}\eta_\eps(x-y)dy\right)\K_th(x^{(m)}) - \J_\eps\K_th(x) \right|_\infty\\
        = & \left|\int_{\X}\eta_\eps(x-y) \left(\sum_{m=1}^M\K_th(x^{(m)})\mathds{1}_{\X_m}(y)dy\right) - \J_\eps\K_th(x) \right|_\infty\\
        = & \left|\J_\eps\left(\sum_{m=1}^M\K_th(x^{(m)})\mathds{1}_{\X_m}(x)\right) - \J_\eps\K_th(x) \right|_\infty\\
        \leq  & \|\J_\eps\| \cdot \sup_{x\in\X}|\K_th(x^{(m)})\mathds{1}_{\X_m}(x)- \K_th(x)|\\
        \leq & \|\J_\eps\|\|\K_t\|L_h\delta_M  \lesssim \delta_M. 
    \end{split}
\end{equation}
\fi

\ym{

\begin{lem}\label{lem: multi_converge}
 Let $\B:\mathcal{C}(\X)\ra \mathcal{C}(\X)$ be any bounded linear operator. Recall  $\bar{G}$, $\bar{H}_g$ of some $g\in\mathcal{C}(\X)$, and $\bar{H}_{\B}$  from Definition \ref{def: projection}. Recall $\J_\eps$ from Proposition \ref{prop: smoothen}.  Define $Z_N^\eps(x) = \J_\eps  Z_N(x)$, $\Pi^{\eps,N}  g = Z^\eps_N(x)^\trans \bar{G}^{\dagger}\bar{H}_g$, and  $\Pi^{\eps,N}\B Z_N(x)^\trans = Z_N^\eps(x)^\trans \bar{G}^{\dagger}\bar{H}_{\B}$. Let $\B^{\eps,N} = \Pi^{\eps,N}\B$. Then, $\|\B^{\eps,N} h -\J_\eps\B h\|\ra 0$ as $N\ra\infty$ for each $\eps>0$. Now, let $\B^N h = \Pi^N\B h$, then $\|\B^{\eps,N} h -\B^N h\| \ra 0$ as $\eps \ra 0$. 
\end{lem}
\begin{pf}
For each $\eps>0$, let $k_x^\eps(y) = \eta_\eps(x-y)$. By definition,   
\begin{equation}\label{E: smooth_projection}
    \begin{split}
        \B^{\eps,N} h(x) %& = Z_N^\eps(x)^\trans \hat{X}^{\dagger}\int_\X Z_N(y)(\B h)(y)dy\\
        & = (\bar{G}^{\dagger}Z_N^\eps(x))^\trans  \int_\X Z_N(y)(\B h)(y)dy\\
        & = \langle Z_N(\cdot)\bar{G}^{\dagger}Z_N^\eps(x), \B h(\cdot)\rangle\\
        & = \langle \Pi_N k_x^\eps(\cdot), \B h(\cdot)\rangle, 
    \end{split}
\end{equation}
%where $\Pi^N$ is from  Definition \ref{def: projection}. 
   % where  $\Pi^Nk_x^\eps(\cdot) = Z_N(y)^\trans\hat{X}^\dagger Z_N^\eps(x)$ is by definition. 
   Therefore, $ |{\B}^{\eps,N} h(x) -\J_\eps\B h(x)|   =  |\langle (\id-\Pi^N)k_x^\eps(\cdot), \B h(\cdot)\rangle|\\
             \leq  \|(\id-\Pi^N)k_x^\eps(\cdot)\|_{L^2}\|\B\|\|h\|_{L^2}$,
             \iffalse
    \begin{equation}
        \begin{split}
              |\widehat{\B}^N h(x) -\J_\eps\B h(x)|   = &|\langle (\id-\Pi^N)k_x^\eps(\cdot), \B h(\cdot)\rangle|\\
             \leq &\|(\id-\Pi^N)k_x^\eps(\cdot)\|_{L^2}\|\B\|\|h\|_{L^2},
        \end{split}
    \end{equation}\fi
    and $\|{\B}^{\eps,N} h -\J_\eps\B h\|\leq \sup_{x\in\X}\|(\id-\Pi^N)k_x^\eps(\cdot)\|_{L^2}\|\B\|\|h\|_{L^2}$. It is well known that the family  $\mathcal{D}:=\{k_x^\eps(\cdot):x\in\X\}$ is equicontinuous and therefore forms a compact subset of $L^2(\X)$. Additionally, since for each $k_x^\eps\in\mathcal{D}$ we have $\|(\id-\Pi^N)k_x^\eps(\cdot)\|_{L^2}\ra 0$ as $N\ra \infty$ by the construction of $\Pi^N$ \cite{williams2015data}, the \textit{uniform boundedness principle} implies $\sup_{x\in\X}\|(\id-\Pi^N)k_x^\eps(\cdot)\|_{L^2}\ra 0$. 

    For the second part, by the standard mollifier property,   $\sup_{x\in\X}|\langle k^\eps_x(\cdot), \zk_i(\cdot)\rangle - \zk_i(x)| \lesssim \eps\cdot \|\nabla z_i(x)\|_\infty$, where the omitted constant is the mean distance under the mollifier and $\|\nabla z_i(x)\|_\infty:=\sup_{x\in\X}|\nabla z_i(x)|$.  This implies that $\sup_{x\in\X}|Z_N^\eps(x) - Z_N(x)| \lesssim \eps \sqrt{N}\max_{1\leq i\leq N}\|\nabla z_i(x)\|_\infty$. Then, by  direct comparison, we have $\|\B^{\eps, N} h -\B^N h\|\leq \sup_{x\in\X}|Z_N^\eps(x) - Z_N(x)|\|\bar{G}^{\dagger}\||\bar{H}_h|\ra 0$ as $\eps\ra 0$. 
\end{pf}

    %%Note that one may attempt to directly use $\Pi^N \mathcal{B} Z_N = Z_N(x)^\trans \hat{X}^\dagger\hat{Y_\B}$ for learning purposes, since this form represents the analytical limit of the data-driven version as the number of samples tends to infinity, and $\Pi^N\B h$ converges to $\B h$ in the $L^2$ norm. However, such convergence does not guarantee convergence w.r.t. $|\cdot|_\infty$ due to the lack of compactness of $\B(\Zk_N)$. Indeed, for any  $h\in \mathcal{C}(\X)$. Then, $|\Pi^Nh-h|_\infty\leq \inf_{v\in\Zk_N}|\Pi^Nh - v|_\infty + |v-h|_\infty \leq (1+\|\Pi^N\|)|\inf_{v\in\Zk_N}|v-h|_\infty$. Although $\inf_{v\in\Zk_N}|v-h|$ converges as $N\ra\infty$, obtaining convergence of the r.h.s. is not straightforward if $\|\Pi^{N}\|$ is not uniformly bounded.

}

\section{Modified Generator Learning Formula}\label{sec: mod}

\textbf{Proof of Lemma \ref{lem: unique_proj}:} \ym{For simplicity, we denote $R_i(x) = \hrr_{\mu, T}^N\zk_i(x)\in\Zk_N$ as a fixed function for each $i$. 
%Then, by definition, $|E_N^i|_\infty\lesssim \inf_{v\in \K_T(\Zk_N)}|v-\K_T\zk_i|_\infty$. We also define $E_N^\text{max}: = \max_{i}|E_N^i|_\infty$. Furthermore, 
Let $\{\psi_k\}$ denote the orthornormal eigenfunctions of $\K_T^N$, which satisfies $\K_T^N\psi_k = e^{\alpha_k T}\psi_k$. It follows that  $e^{\operatorname{Re}({\alpha_k}) T}\|\psi_k\| \lesssim \|\K_T\psi_k\|+ E_N^\text{max}\leq (1 + E_N^\text{max} e^{-\omega T})e^{\omega T}$,  
\iffalse
\begin{small}
    \begin{equation}
        \begin{split}
            |\K_T\psi_k|_\infty& \lesssim e^{\operatorname{Re}({\alpha_k}) T}|\psi_k|_\infty + E_N^\text{max}\\
            & \leq e^{\omega T}  + E_N^\text{max}\\
            & = (1 + E_N^\text{max} e^{-\omega T})e^{\omega T}
        \end{split}
    \end{equation}
\end{small}\fi
 %On the other hand, $e^{\operatorname{Re}({\alpha_k}) T}|\psi_k|_\infty\leq |\K_T\psi_k|_\infty + r\leq e^{\omega t} + 2r$, 
 which implies that $\operatorname{Re}({\alpha_k})\leq \omega + \frac{1}{t}\log(1+E_N^\text{max}e^{-\omega t})\leq \omega + \frac{e^{-\omega T}E_N^\text{max}}{T}$ for all $k$. Therefore, $\|\K_T^N\|\leq Me^{\tilde{\omega}T}$ where $|\tilde{\omega}-\omega|\lesssim E_N^\text{max}$. Then, for all $\mu\in(\tilde{\omega}, \infty)$, one can prove in the same manner as in Proposition \ref{prop: contraction}  (with a standard proof provided in the supplementary material) that  $\mathcal{T}V_i^N(x) = R_i(x)+ e^{\mu T}\K_T^NV_i^N(x)$ is a contraction mapping, indicating the uniqueness of the equation in the statement. 
Furthermore, 
\begin{equation}
    \begin{split}
        & \|V_i^N-V_i\|\\
        \leq & \|R_i-\rr(\mu)\zk_i\| +  e^{-\mu T}\|\K_TV_i -\K_T^NV_i\| \\
        &+ e^{-\mu T}\|\K_T^NV_i-\K_T^NV_i^N\|\\
        \lesssim & E_{\sss \sum}+  \|E_N^i\|+ e^{(\tilde{\omega}-\mu) T}\|V_i^N-V^N\|, 
    \end{split}
\end{equation}
implying that $\|V_i^N-V^N\|\lesssim E_{\sss \sum}+ E_N^\text{max}$.} \pfbox

\textbf{Proof of Theorem \ref{thm: approx_modify}:} \ym{We give a sketch of the proof, since all error bounds follow from triangle inequalities and the convergence analysis in Section \ref{sec: sub_compact}. For simplicity, we   define  $\A:=(\lambda-\mu)\rr(\mu)+\id$,  $\B:=\lambda\mu\rr(\mu)-\lambda\id$,  $\bar{A} = \bar{G}^{\dagger}\bar{H}_\A$, and $\bar{A}_\delta^+ = (\bar{A}^\trans \bar{A}+\delta \id)^{-1}\bar{A}^\trans$ for any $\delta>0$. %Suppose the exact value of $\rr(\mu)$, we can correspondingly define   $\hat{X}$, $\hat{Y}_{\A}$, and $\hat{Y}_{\B}$   defined as in Definition \ref{def: projection}, as well as $\A_N$ and $\B_N$ similarly. 
Since $\A$ is invertible, for any $V(x) = Z_N(x)^\trans v\in \Zk_N$, we set $g = \A^{-1} V$, $\xi =\hat{A}_\delta^\dagger v$, and $\hat{g}_{\sss N} = Z_N(x)^\trans \xi$. %Then, it can be verified that $\hat{g}_{\sss N}$ is the optimal candidate in $\Zk_N$ such that $\widehat{\A}\hat{g}_N$ provides the best approximation of $h$. 
 Additionally,  we notice that $\A g = \Pi^N\A g = V$. We also denote $g_{\sss N} = \Pi^N g$, and $\bar{g}_{\sss N}= Z_N(x)^\trans\bar{A}_\delta^+\theta$. 
Then $\bar{g}_{\sss N} = \arg\min_{W\in\Zk_N}(\|\Pi^N \A W-V\|_{L^2} + \delta
 \|W\|_{L^2}^2)$, and $\|\Pi^N \A \bar{g}_{\sss N}-V\|_{L^2} + \delta\|\bar{g}_{\sss N}\|_{L^2}^2\leq \|\Pi^N \A g_{\sss N}-V\|_{L^2}+ \delta\|g_{\sss N}\|_{L^2}^2 = \|\Pi^N \A (g_{\sss N}-g)\|_{L^2} + \delta\|g_{\sss N}\|_{L^2}^2\lesssim E_N(g) + \delta\|g_{\sss N}\|^2_{L^2}$, where $E_N(g)$ of the last inequality is due to the strong convergence of $\Pi^N \A$ to $\A$. Therefore, $\|\Pi^N \A \bar{g}_{\sss N}-\Pi^N\A g\|_{L^2}\lesssim E_N(g) + \delta(\|g_{\sss N}-\bar{g}_{\sss N}\|)$. 
 One can then follow the same procedure as in the proof of Lemma \ref{lem: projection} to show that $\|\Pi^N \A \bar{g}_{\sss N}-\Pi^N\A g\|\lesssim E_N(g) + \delta(\|g_{\sss N}-\bar{g}_{\sss N}\|)$\footnote{\ym{One can let $u = \bar{g}_{\sss N}-g$.  Similar to \eqref{E: smooth_projection},  $\|\Pi^{\eps, N}\A u\|\leq \sup_{x\in\X}\|k_x^\eps(\cdot)\|_{L^2}\|\Pi^N u\|_{L^2}\lesssim\|u\|_{L^2}$ and $\|\Pi^{\eps, N}\A u- \Pi^N\A u\|\lesssim \eps$ for any $\eps>0$. The convergence follows by the same argument as in Remark~\ref{rem: construction}, and we omit the details due to repetition.}}. Then,  $\|\bar{g}_{\sss N}-g\|\leq \|\A^{-1}\|\|\A g - \A\bar{g}_{\sss N}\|\leq \|\A^{-1}\|(\|V - \Pi^N\A\bar{g}_{\sss N}\|+ \|\Pi^N\A \bar{g}_{\sss N}-  \A\bar{g}_{\sss N}\|)\lesssim E_N(g) + \delta\|\bar{g}_{\sss N}-g\| + E_N(\A\bar{g}_{\sss N})$, which implies $\|\bar{g}_{\sss N}-g\|\lesssim E_N(g) + \delta  + E_N(\A\bar{g}_{\sss N})$ for sufficiently small $\delta$. One can let $E_N^\Pi:=E_N(g)+E_N(\A\bar{g}_{\sss N})$, whence $\|\bar{g}_{\sss N}-g\| \lesssim E_N^\Pi + \delta$.

 %Then $\bar{g}_{\sss N}, = \arg\min_{W\in\Zk_N}(|\Pi^N \A W-V|_\infty + \delta|W|_\infty^2)$, and $|\Pi^N \A \bar{g}_{\sss N}-V|_\infty + \delta|\bar{g}_{\sss N}|_\infty^2\leq |\Pi^N \A g_{\sss N}-V|_\infty + \delta|g_{\sss N}|_\infty^2 = |\Pi^N \A (g_{\sss N}-g)|_\infty + \delta|g_{\sss N}|_\infty^2\lesssim E_N(g) + \delta|g_{\sss N}|^2$, where $E_N(g)$ of the last inequality is due to the strong convergence of $\Pi^N \A$ to $\A$. Therefore, $|\Pi^N \A \bar{g}_{\sss N}-V|_\infty\lesssim E_N(g) + \delta(|g_{\sss N}|^2-|\bar{g}_{\sss N}|^2)$, and $|\bar{g}_{\sss N}-g|_\infty\leq \|\A^{-1}\||\A g - \A\bar{g}_{\sss N}|_\infty\leq \|\A^{-1}\|(|V - \Pi^N\A\bar{g}_{\sss N}|_\infty+ |\Pi^N\A \bar{g}_{\sss N}-  \A\bar{g}_{\sss N}|_\infty)\lesssim E_N(g) + \delta|\bar{g}_{\sss N}-g|_\infty + E_N(\A\bar{g}_{\sss N})$, which implies $|\bar{g}_{\sss N}-g|_\infty \lesssim E_N(g) + \delta  + E_N(\A\bar{g}_{\sss N})$ for sufficiently small $\delta$. 

On the other hand, since $\widehat{A}$ and $\A$ differ only through the evaluation of $\rr(\mu)$, with error depending continuously on $E_{\sss \sum} + E_N^{\text{max}}$ (as in Lemma \ref{lem: unique_proj}), we obtain $E_A:=\|\bar{A}-\hat{A}\|\lesssim E_{\sss \sum} + E_N^{\text{max}}$ as well. We then have $\|\hat{A}_\delta^+-\bar{A}_\delta^+\|\leq (\frac{1}{\delta}+\frac{(\|\bar{A}\|+\|\hat{A}\|)\max\{\|\bar{A}\|, \|\hat{A}\|\}}{\delta^2})E_A$ \cite{tikhonov1977solutions}. It can be verified that with $\delta = \mathcal{O} (E_A^{1/3})$, we have a H\"{o}lder bound $\|\hat{A}_\delta^+-\bar{A}_\delta^+\|\lesssim E_A^{1/3}$, and thereby $\|\hat{g}_{\sss N}-\bar{g}_{\sss N}\|\lesssim  E_A^{1/3}$  for sufficiently small $E_A$. 
 
Now, let $V = \widehat{B}h_\theta$, $v = \hat{B}\theta$, and $g = \A^{-1}V$. Then, $\ll_{\lambda}^N h_\theta(x) = Z_N(x)^\trans\hat{A}^{\dagger}_\delta \hat{B}\theta   \hat{g}_{\sss N}$,  $\ll_\lambda h_\theta(x) =  \A^{-1}\B h_\theta$, and $\|\ll_\lambda^Nh_\theta-\ll_\lambda h_\theta\|  \leq \|\hat{g}_{\sss N}- g\|+  \|g-\A^{-1}\B h_\theta\| \lesssim \|\hat{g}_{\sss N}- g\| + \|\widehat{\B}h_\theta - \B h_\theta\|
      \lesssim  E_N^\Pi+  (E_{\sss \sum} + E_N^{\text{max}})^{1/3} + (E_{\sss \sum} + E_N^{\text{max}})\lesssim  E_N^\Pi+ (E_{\sss \sum} + E_N^{\text{max}})^{1/3}$,   where the third term in the third inequality again follows from the fact that $\widehat{B}$ and $\B$ differ only through the evaluation of $\rr(\mu)$, with error  continuous to $E_{\sss \sum} + E_N^{\text{max}}$ (as in Lemma \ref{lem: unique_proj}). The conclusion   follows. }\pfbox

\section{Supplementary Technical Proofs}

\subsection{Technical Proofs of Fundamental Properties}
\label{sec: proof_fundamental}

We first present the following facts about $\{\K_t\}$
\begin{prop}\label{prop: fact1}
   %For system \eqref{E: sys}, there exist  constants $\omega\geq 0$ and $\Sigma\geq 1$ such that $\norm{\K_t}\leq Ce^{\omega t}$ for all $t\geq 0$. In addition, 
   For any $\lambda\in\C$, the family $\set{\K_{t,\lambda}}_{t\geq 0}$, where  
   \begin{equation}\label{E: semi_new}
       \K_{t,\lambda}: =e^{\lambda t}\K_t
   \end{equation}
   is a $\mathcal{C}_0$-semigroup with generator $\ll+\lambda \id: \dom(\ll)\ra \mathcal{C}(\X)$. 
\end{prop}

\begin{pf}
    %The first part follows directly by \cite[Theorem 1.2.2]{pazy2012semigroups}. 
    The semigroup property of $\set{\K_{t, \lambda}}_{t\geq 0}$ follows easily by the definition. To verify its generator, we have %%$\frac{\K_{t, \lambda} h - h}{t} =\frac{e^{\lambda t}\K_{t, \lambda}h - e^{\lambda t} h }{t} + \frac{ e^{\lambda t} h -h }{t}.$
    \begin{equation}%\label{E: semigroup_new}
        \begin{split}
            \frac{\K_{t, \lambda} h - h}{t} =\frac{e^{\lambda t}\K_{t}h - e^{\lambda t} h }{t} + \frac{ e^{\lambda t} h -h }{t}.
        \end{split}
    \end{equation}
    For all $h\in\dom(\ll)$, the limit exists as $t\downarrow 0$ for each  r.h.s. term. It then follows that $\lim_{t\downarrow0}\frac{\K_{t, \lambda} h - h}{t} = \ll + \lambda \id$. 
\end{pf}

\textbf{Proof of Theorem \ref{thm: conv}}:
For any $\tomg\geq \omega$, we consider $\set{\K_{t,-\tomg}}_{t\geq 0}$, where $\K_{t,-\tomg}$ is defined as \eqref{E: semi_new}. It is clear that $\|\K_{t,-\tomg}\|\leq \Sigma$ for all $t\geq 0$. 
%$$\|\K_{t,-\tomg}\|\leq M,\;\;t\geq 0.$$
By Proposition \ref{prop: fact1}, $\set{\K_{t,-\tomg}}_{t\geq 0}$ also admits the generator $\ll-\tomg \id$.  Within $C^1(\X)$, for any $\lambda>\tomg$, we have
\begin{small}
        \begin{equation}\label{E: L_split}
        \begin{split}
            \ll  - \lll    %=&\ll  - \ll_{\lambda + \omega} + \ll_{\lambda + \omega} - \lll  \\
             = &\ll  - \tomg \id + \tomg\id - \lll\\
              = &(\ll  - \tomg \id - (\ll  - \tomg \id)_{\lambda- \tomg})\\
            & + ((\ll  - \tomg \id)_{\lambda- \tomg}+\tomg\id- \lll)\\
             =:&  \O_1 + \O_2, 
        \end{split}
    \end{equation}
\end{small}

\noindent where $(\ll  - \tomg \id)_{\lambda- \tomg}$ is the Yosida approximation of $\ll  - \tomg \id$. It suffices to show the bound for each of the  two terms on the r.h.s. of \eqref{E: L_split}. 

To show the bound for $\O_1$, we consider an alternative  norm  for $h \in \mathcal{C}^1(\X)$, defined as $\norm{h}_\infty=\sup_{t\geq 0}\|\K_{t,-\tomg}h\|$. 
It can be verified that $\|h\|\leq \|h\|_\infty\leq \Sigma\|h\|$.  For $\|h\|_\ll$, we  define the alternative norm $\|h\|_{\ll,\text{alt}} =  \norm{h}_\infty + \norm{\ll h}_\infty.$
\iffalse
\begin{equation}
    \|h\|_\ll =  \norm{h}_\infty + \norm{\ll h}_\infty. 
\end{equation}\fi
In addition, 
$\|\K_{t,-\tomg}h\|_\infty= \sup_{s\geq 0}\|\K_{s,-\tomg}\K_{t,-\tomg}h\|\leq \sup_{s\geq 0}\|\K_{s,-\tomg}h\|=\|h\|_\infty,$
which demonstrates the  contraction property  w.r.t. $\norm{\cdot}_\infty$. Then, by the definition of the resolvent operator, we have $((\lambda - \tomg)I - (\ll-\tomg\id))\rr(\lambda-\tomg; \ll  - \tomg \id)  = \id$, which implies
\begin{small}
    \begin{equation}\label{E: O1}
    ((\lambda - \tomg)I - (\ll-\tomg\id))\rr(\lambda-\tomg; \ll  - \tomg \id)(\ll-\tomg\id)  = \ll-\tomg\id.
\end{equation}
\end{small}

\noindent Recall the definition of the Yoshida approximation for $(\ll  - \tomg \id)_{\lambda- \tomg}$. 
For any $h\in \mathcal{C}^2(\X)$, by expanding \eqref{E: O1}, we can obtain 
$\O_1h = -(\ll  - \tomg \id) \rr(\lambda-\tomg; \ll  - \tomg \id) (\ll  - \tomg \id) h$. Consequently,  
\begin{small}
    \begin{equation*}
    \begin{split}
        |\O_1h|_\infty & \leq \|O_1h\|_\infty\leq (\lambda-\tomg)^{-1}\|(\ll  - \tomg \id) h\|_{\ll-\tomg\id, \text{alt}}\\
        & \leq \Sigma(\lambda-\tomg)^{-1}\|(\ll  - \tomg \id) h\|_{\ll-\tomg\id}\\
        & \leq \Sigma(\lambda-\tomg)^{-1}(\|\ll h \|_\ll + \tomg \|h\|_\ll + \tomg\|h\|)\\
        & \leq \Sigma(\lambda-\tomg)^{-1}(\|\ll h \|_\ll + (2\tomg + 1)\|h\|_\ll),
    \end{split}
\end{equation*}
\end{small}

\noindent where the second inequality can be proved in the same way as in \cite[Chap. I, Lemma 3.2]{pazy2012semigroups}, and the last two inequalities are based on the definitions of the norms $\|\cdot\|_{\ll-\tomg\id}$ and $\|\cdot\|_{\ll}$, obtained by expanding them. Since $\mathcal{C}^2(\X)$ is dense in $\mathcal{C}^1(\X)$, and the operator $\O_1$ is uniformly bounded \cite[Theorem 1.3.1]{pazy2012semigroups} and hence continuous on $\mathcal{C}^1(\X)$, 
we have  
\begin{equation}\label{E: conv_omg}
    \ll  - \tomg \id  = \slim_{\lambda\ra\infty}(\ll  - \tomg \id)_{\lambda- \tomg}
\end{equation}
on~$(\mathcal{C}^1(\X), \|\cdot\|_{\ll, \text{alt}})$.  
Due to the norm equivalence between~$\|\cdot\|_{\ll, \text{alt}}$ and $ \|\cdot\|_\ll$, the  convergence \eqref{E: conv_omg} also holds on $(\mathcal{C}^1(\X), \|\cdot\|_\ll)$. 

We now work on the bound for $\O_2$. One can show by a direct calculation that, for any $h\in \mathcal{C}^1(\X)$, 
\begin{small}
    \begin{equation}
    \begin{split}
        \|\O_2h\| & = \|2\tomg h- \tomg( 2\lambda-\tomg)R(\lambda;\ll)h\|\\
        & = \|\tomg(\tomg R(\lambda;\ll) - 2\ll R(\lambda; \ll))h\|\\
        & \leq \Sigma(\lambda-\tomg)^{-1}(\tomg^2\|h\|+2\tomg \|\ll h\|)\\
        & \leq \Sigma_0(\lambda-\tomg)^{-1} \|h\|_\ll,
    \end{split}
\end{equation}
\end{small}

\noindent where $\Sigma_0=\Sigma\cdot\max\{\tomg^2,2\tomg\}$.

The conclusion follows by combining both parts and considering $\lambda\ra\infty$. \pfbox

\begin{rem}
    In the proof of Theorem \ref{thm: conv}, we have implicitly demonstrated the effects of $C$ and $\omega$ in the semigroup estimation $\norm{\K_t}\leq \Sigma e^{\omega t}$. Intuitively, $\Sigma$ represents the uniform scaling of the magnitude of the Koopman operator, while $\omega$  indicates the dominant exponential growth or decay rate of the flow on $\mathcal{C}(\X)$. 

    Specifically, one can shift the original generator $\ll$ to a stable generator $\ll-\tomg\id$, which generates a contraction semigroup w.r.t. a norm equivalent to $\|\cdot\|$, and then shift back to $\ll$ by adding $\tomg\id$. The analysis in the proof determines the error ($\O_1$) between the stable generator $\ll-\tomg\id$ and its Yosida approximation, as well as the error ($\O_2$) between the   approximation $\ll_\lambda$ of the original generator  and the   approximation $(\ll  - \tomg \id)_{\lambda- \tomg}+\tomg\id$ of the ``shift-back'' operator.  
    
    Note that this shifting strategy used in the proof is advantageous for applying the analysis of contraction semigroups, but the shift inevitably increases the cost. Nonetheless, the overall error of using the Yoshida approximation on $\rho(\ll)$ is always dominated by a reciprocal term.\Qed
    
 %%Working within the semigroup topology, one can always convert a $C_0$ semigroup into a contraction semigroup and approximate the generator using this connection, i.e., $(\ll  - \tomg \id)_{\lambda- \tomg}+\tomg\id$. However, this is \textit{not necessary} given that, for arbitrarily large $\lambda$, $\lll$ is arbitrarily close to $(\ll  - \tomg \id)_{\lambda- \tomg}+\tomg\id$ in a strong sense.  The convergence rate remains similar. 
\end{rem}

\begin{prop}\label{prop: pseudo}
    %Recall  $\rr(\lambda)$ in \eqref{E: pseudo}. 
    The  revolvent identity holds for $\{\rr(\lambda)\}$:
    \begin{equation*}\label{E: resolvent_id}
    \rr(\lambda)-\rr(\mu) =  (\mu-\lambda)\rr(\mu)\rr(\lambda), \;\;\lambda,\mu \in (\omega, \infty). 
\end{equation*}
\end{prop}
\begin{pf}
     Note that $\int_0^t(\mu-\lambda)e^{-\lambda s}e^{-\mu(t-s)}\dd s = e^{-\lambda t}- e^{-\mu t}$. 
    Then,  
    \begin{small}
        \begin{equation}
        \begin{split}
            &\rr(\lambda)h(x)-\rr(\mu)h(x)\\ = & \int_0^\infty (e^{-\lambda t}- e^{-\mu t})\K_th(x) \dd t\\
             = & \int_0^\infty \left(\int_0^t(\mu-\lambda)e^{-\lambda s}e^{-\mu(t-s)}ds\right)\K_th(x)  \dd s \dd t\\
             %= & (\mu-\lambda) \int_0^\infty\int_0^t e^{-\lambda s}e^{-\mu(t-s)}\K_th(x) \dd s \dd t\\
             = & (\mu-\lambda) \int_0^\infty\int_s^\infty e^{-\lambda s}e^{-\mu(t-s)}\K_th(x) \dd t \dd s\\
             = & (\mu-\lambda) \int_0^\infty\int_0^\infty e^{-\lambda s}e^{-\mu r}\K_{r+s}h(x) \dd r \dd s\\
             = & (\mu-\lambda) \int_0^\infty e^{-\lambda s} \K_s \int_0^\infty e^{-\mu r}\K_{r}h(x) \dd r \dd s\\
            % = &(\mu-\lambda) \int_0^\infty e^{-\lambda s} \K_s \rr(\mu) h(x)  \dd s\\
             = & (\mu-\lambda) \rr(\lambda) \rr(\mu) h(x),
        \end{split}
    \end{equation}
    \end{small}
    
   \noindent where the interchange of the order of integration is valid by Fubini’s theorem, since $\int_0^\infty\int_0^t \|e^{-\lambda s}e^{-\mu(t-s)}\K_th(x) \|\dd s \dd t \leq C \int_0^\infty\int_0^t e^{-\lambda s}e^{-\mu(t-s)} e^{\omega t} \dd s \dd t <\infty$. 
\end{pf}
%\fi

\vspace{0.5em}
\textbf{Proof of Proposition \ref{prop: R_form}}:
    For all $h\in  \dom(\ll)$ and $x\in \X$, we have
    \begin{small}
    \begin{equation*}
        \begin{split}
            & \rr(\lambda)(\lambda \id-\ll)h(x) \\
            = & \lambda \rr(\lambda)h(x) - \rr(\lambda)\ll h(x)\\
            = & \lambda \int_0^\infty e^{-\lambda t}h(\phi(t,x)) dt - \int_0^\infty e^{-\lambda t}\ll h(\phi(t, x) dt\\
            = & \left.-h(\phi(t,x))e^{-\lambda t}\right|_0^\infty + \int_0^\infty e^{-\lambda t} d(h(\phi(t, x)) \\
            & - \int_0^\infty e^{-\lambda t}\ll h(\phi(t, x) dt\\
            = & h(\phi(0,x)) = h(x),
        \end{split}
    \end{equation*}     
    \end{small}
    
 \noindent    where we have used the fact that the time derivative along the trajectories of \eqref{E: sys} is such that $dh(\phi(t,x))/dt=\ll h(\phi(t,x))$. 

    Recall   \eqref{E: semi_new}.  For all $h\in \mathcal{C}(\X)$, all $x\in\X$, and all $t\geq 0$,
    \begin{small}
    \begin{equation*}
        \begin{split}
            &\rr(\lambda)h(x)\\
            =&\int_0^t \K_{s,-\lambda} h(x) ds + \int_t^\infty \K_{s,-\lambda} h(x) ds\\
            =&\int_0^t \K_{s,-\lambda} h(x) ds + \int_0^\infty \K_{t+s,-\lambda} h(x) ds\\
            %=&\int_0^t \K_{s,-\lambda} h(x) ds + \int_0^\infty e^{-(t+s)\lambda}\K_{s} h(\phi(t,x)) ds\\
            =&\int_0^t e^{-\lambda s}\K_{s} h(x) ds + e^{-\lambda t} \rr(\lambda)h(\phi(t,x)).
        \end{split}
    \end{equation*}        
    \end{small}
However, 
\begin{small}
\begin{equation*}
    \begin{split}
        \K_t\rr(\lambda)h(x) & = \int_0^\infty e^{-\lambda s} \K_s h(\phi(t,x))ds= \rr(\lambda) \K_t h(x). 
    \end{split}
\end{equation*}  
\end{small}
Therefore, 
\begin{small}
\begin{equation*}
\begin{split}
        &\frac{\K_t\rr(\lambda)h(x)-\rr(\lambda)h(x)}{t}\\
        =&\frac{\K_tR(\lambda)h(x)-e^{-\lambda t} \rr(\lambda)h(\phi(t,x))}{t}\\
        &+\frac{e^{-\lambda t} \rr(\lambda)h(\phi(t,x))-\rr(\lambda)h(x)}{t}\\
        =& \frac{\K_t\rr(\lambda)h(x)-e^{-\lambda t} \K_t\rr(\lambda)h(x)}{t} - \frac{1}{t}\int_0^t e^{-\lambda s}\K_{s} h(x) ds.
\end{split}
\end{equation*}
\end{small}

\noindent With $t\downarrow0$ on both sides,
we have $\ll \rr(\lambda)h(x) = \lambda \rr(\lambda)h(x)- h(x)$,
which is equivalent as $(\lambda\id-\ll)\rr(\lambda)h(x) = h(x)$. \pfbox

\textbf{Proof of Proposition \ref{prop: contraction}}:  
    Note that $\rr_\mu h(x) =\rr_{\mu, T}h(x) + \int_T^\infty e^{-\mu t} \K_th(x)dt$. By change of variable, 
    \begin{equation}
    \begin{split}
                \int_T^\infty e^{-\mu t} \K_th(x)dt & = \int_0^\infty e^{-\mu (T+s)}\K_{T+s}h(x)ds\\
                & = e^{-\mu T} \int_0^\infty e^{-\mu s}\K_T\K_sh(x)ds\\
                & = e^{-\mu T} \K_T \int_0^\infty e^{-\mu s}\K_sh(x)ds,
    \end{split}
    \end{equation}
    which proves the first part of the statement. 

    Now let $V, W\in\D(\mathcal{T})$. Then, 
    \begin{equation}
        \begin{split}
            \|\T_h V-\T_h W\| & = \sup_{x\in\Omega}e^{-\mu T}|\K_T V(x)-\K_T W(x)|\\
            & \leq  e^{-\mu T}\|\K_T\|\cdot  \| V(x)- W(x)\|\\
            & \leq e^{(\omega-\mu) T}\|V-W\|. 
        \end{split}
    \end{equation}
    The contraction   property follows since $e^{(\omega-\mu) T}<1$. \pfbox

\subsection{Finite-Rank Projection of Bounded Linear Operators}
We show the following projection result and its error.

\begin{lem}\label{lem: projection}
For any $v(\cdot) = Z_N(\cdot)^\trans\theta\in\Zk_N$ with some column vector $\theta$, we have $\Pi^N v(x) = v(x)$. %Furthermore, let $\mathcal{O}_N(h) = \min_\theta|h-Z_N(x)^\trans\theta|_\infty$    for any  $h\in \mathcal{C}(\X)$. Then, $|\Pi^Nh-h|_\infty\lesssim \mathcal{O}_N(h)$. Particularly, for any bounded linear operator $\B:\mathcal{C}(\X)\ra\mathcal{C}(\X)$,  we have $\sup_{x\in\Omega}|\Pi^N\B Z_N(x) - \B Z_N(x)|\lesssim \max_{h\in \B Z_N}\mathcal{O}_N(h)$
\end{lem}
\begin{pf}
For $v(x) = Z_N(x)^\trans\theta$, $\bar{Y}_v = [\int_\X Z_N(y)Z_N^\trans (y)dy]\theta= \bar{X} \theta$ and  $\Pi^N v(x) =   Z_N(x)^\trans \bar{X}^\dagger \bar{X} \theta$. It is clearly that $\Pi^N$ is bounded. 
    Suppose $\theta \in \operatorname{ran}(\bar{X})$. Then we have $\bar{X}^\dagger \bar{X}\theta = \theta$, which implies $ \Pi^N h(x) = h(x)$. Otherwise, we have $\bar{X}\theta = 0$, which implies that   $Z_N(x)^\trans \theta$ is orthogonal to every $\zk_i$. Consequently,    $Z_N(x)^\trans\theta = 0$ a.e.. Due to the continuity of each $\zk_i$, this further implies that $Z_N(x)^\trans\theta = 0$ pointwise on $\X$.  By the positive semidefinite property of $\bar{X}$, we then have   $0 = \Pi^Nh(x) = h(x) = 0$. 
\iffalse
    For any $h\in \mathcal{C}(\X)$, let $\theta^\star$ be the minimizer of $\mathcal{O}_N(h)$, and set $v(x) = Z_N(x)^\trans\theta^\star$. Then \begin{small}
        $|S_N h -h|_\infty \leq |S_N(h-v)|_\infty + |v-h|_\infty\lesssim \mathcal{O}_N(h)$.
    \end{small} The last assertion follows by enumerating the error bound for each element in the vector of functions.\fi
\end{pf}

\begin{lem}\label{lem: error_sup}
   For any  $h\in \mathcal{C}(\X)$, recall that $E_N(h) =\inf_{v\in\Zk_N}\|v-h\|$. Then,  $\|\Pi^Nh-h\|\leq   (1+\|\Pi^N\|)E_N(h)$.
\end{lem}
\begin{pf}
    $\|\Pi^Nh-h\|\leq \inf_{v\in\Zk_N}\|\Pi^Nh - v\| + \|v-h\|\leq (1+\|\Pi^N\|)|\inf_{v\in\Zk_N}\|v-h\|$, where we have used the property that $v=\Pi^Nv$ for any $v\in\Zk_N$ from Lemma \ref{lem: projection}. 
\end{pf}

\subsection{Error Estimation of Gauss–Legendre Quadrature}
We first rephrase the technical result from \cite[Chapter 5]{kahaner1989numerical}  on the mathematical error of Gauss–Legendre integration in the following lemma.
\begin{lem}[Mathematical error]
    Consider Gauss-Legendre method for the integral $\int_0^T g(t)dt$ with $\Gamma$ interpolation points. Then, the error of the Gauss-Legendre integration is bounded by
    \begin{equation}\label{E: gauss-legendre_quadrature_error}
\begin{gathered}
   T^{2\Gamma + 1}\cdot\mathcal{E}(\Gamma) \cdot  \sup_{t\in[0, T]}|g^{(2\Gamma)}(t)|, 
\end{gathered}
\end{equation}
where $\mathcal{E}(\Gamma)=\frac{(\Gamma!)^4}{(2\Gamma + 1) [(2\Gamma)!]^3}$. 
\end{lem}

In our case,  $g(t) = e^{-\mu t}\K_th(x)$ for  fixed $x\in\Omega$. In light of Appendix \ref{sec: facts} and the fact that, for any Lipschitz $f$, one can always find a smooth approximation with arbitrarily small error, whose corresponding flow is smooth and arbitrarily close to the true flow. We will, for simplicity, abuse notation and let 
$f$ denote its smooth approximation throughout this section. This setup allows for easier estimation of the error bound. We   work on monomial dictionary functions to  exemplify the estimation, i.e., we work on $h(x)=x^\alpha$, where $\alpha=(\alpha_1, \cdots, \alpha_d)$ is a multi-index and $x^\alpha=\prod_{i=1}^dx_j^{\alpha_j}$. We say the monomial function $h$ is of order $N$, denoted as $h\in P^N(\Omega)$, if $|\alpha|=\sum_{i=1}^d\alpha_i=N$, and $\alpha_i\in\mathbb{N}$ for all $i\in\set{1, 2, \cdots, d}$.

In the formula \eqref{E: gauss-legendre_quadrature_error}, the error bound can be separated into three   terms, where the first and second terms are straightforward to estimate. Particularly, The error term $\mathcal{E}(\Gamma)$  is sufficiently small for $\Gamma\geq 10$. The following lemma establishes a tight bound for 
$\mathcal{E}(\Gamma)$  with a simplified expression.  

\begin{prop}Denote $\overline{\mathcal{E}}(k)=\left(\frac{1}{8k^2}\right)^k$. 
Then, for each $k\geq 1$, 
    $\mathcal{E}(k)\leq \overline{\mathcal{E}}(k)$. 
\end{prop}
\begin{pf}
By Robbins' bounds \cite{robbins1955remark} on the Stirling's formula for   factorials, $k! =  \sqrt{2\pi k}\left(\frac{k}{e}\right)^ke^{r_k}$,  
where $r_k$ satisfies $\frac{1}{12k+1}< r_k<\frac{1}{12k}$. 
Then, 
\begin{small}
        \begin{equation}\label{E: bound_num}
        (k!)^4=\left[\sqrt{2\pi k}\left(\frac{k}{e}\right)^k e^{\frac{1}{12k}}\right]^4 = (2\pi k)^2\left(\frac{k}{e}\right)^{4k}e^{4r_k}  
    \end{equation}
\end{small}

\noindent and 
\begin{small}
    \begin{equation}\label{E: bound_den}
\begin{split}
         [(2k)!]^3 =   \left[\sqrt{2\pi(2k)}\left(\frac{2k}{e}\right)^{2k}e^{r_{2k}}\right]^3 
        =   (4\pi k)^{3/2}\left(\frac{2k}{e}\right)^{6k}e^{3r_{2k}}. 
\end{split}
\end{equation}
\end{small}

\noindent Combining \eqref{E: bound_num}, \eqref{E: bound_den}, and the fact that $2k+1>2k$, we have
\begin{equation}
\begin{split}
        \mathcal{E}(k)&<\frac{(2\pi k)^2\left(\frac{k}{e}\right)^{4k}e^{4r_k}}{(2k)\cdot (4\pi k)^{3/2}\left(\frac{2k}{e}\right)^{6k}e^{3r_{2k}}}\\
        & = \frac{(2\pi k)^2}{(2k)\cdot (4\pi k)^{3/2}} \cdot \frac{k^{4k}}{(2k)^{6k}}\cdot e^{2k}\cdot e^{4r_k-3r_{2k}}\\
        & < \frac{\pi^{3/2}}{k^{1/2}}\cdot\left(\frac{e^2}{64k^2}\right)^k\cdot e^{\frac{1}{3k}-\frac{3}{24k+1}}\\
        & < \frac{\pi^{3/2}}{k^{1/2}}\cdot\left(\frac{e^2}{64k^2}\right)^k\cdot e^{\frac{15k+1}{72k^2}}. 
\end{split}
\end{equation}
Note that, since $\frac{\pi^{3/2}}{k^{1/2}} e^{\frac{15k+1}{72k^2}}$ monotonically decreases and $\left(\frac{8}{e^2}\right)^k$ monotonically increases as $k$ increases, there exists a unique zero for $\frac{\pi^{3/2}}{k^{1/2}} e^{\frac{15k+1}{72k^2}}=\left(\frac{8}{e^2}\right)^k$, which is strictly less than $10$. Then, for $k\geq 10$, we have $\frac{\pi^{3/2}}{k^{1/2}} e^{\frac{15k+1}{72k^2}}<\left(\frac{8}{e^2}\right)^k$ and thereby $\mathcal{E}(k)<\left(\frac{8}{e^2}\right)^k\left(\frac{e^2}{64k^2}\right)^k=\left(\frac{1}{8k^2}\right)^k$. For $k=1, 2, \cdots, 9$, one can verify by direct calculation that $\mathcal{E}(k)< \left(\frac{1}{8k^2}\right)^k$ also holds. 
\end{pf}

We now derive the derivatives of the integrand.
 \begin{prop}\label{prop: derivative}
For any $k\in\mathbb{N}_0$, we have
\begin{equation}\label{E: formula_integrand}
    \frac{d^k}{dt^k}(e^{-\mu t}\K_th) = e^{-\mu t}\K_t((\ll-\mu\id)^kh).
\end{equation}
 Expanding $(\ll-\mu\id)^k$, it is also equivalent to the following explicit formula
    \begin{equation}\label{E: fomula_explicit}
        \frac{d^k}{dt^k}(e^{-\mu t}\K_th)(x) = e^{-\mu t}\sum_{i=0}^k  \begin{pmatrix}
        k \\ i
    \end{pmatrix}(-\mu)^{k-i}\K_t(\ll^i h)(x).
    \end{equation}
\end{prop}
 \begin{pf}
    By the chain rule, 
    \begin{equation}\label{E: derivative_K}
        \frac{d}{dt}\K_th(x)  = \ll (\K_th(x)) = \K_t(\ll h)(x). 
    \end{equation}
    %where the associativity can be implied by \cite[Proposition 3.5]{meng2024koopman}. 
    Then, 
    \begin{equation}
    \begin{split}
                \frac{d^2}{dt^2}\K_th(x) & = \frac{d}{dt} \K_t(\ll h)(x)\\
                & = (\nabla(\nabla h\cdot f)\cdot f)(\phi(t,x))\\
                & = \ll^2h(\phi(t,x))\\
                & = \K_t(\ll^2 h)(x).
    \end{split}
    \end{equation}
    By induction, the $k$-th derivative is  $\frac{d^k}{dt^k}(\K_th)(x) =  \K_t(\ll^kh)(x)$, where 
    $\ll^k h = \ll(\ll^{k-1}h) = \nabla(\ll^{k-1}h)\cdot f.$

    Now, by the product rule, 
\begin{equation}
\begin{split}
        \frac{d}{dt}(e^{-\mu t}\K_th)(x) & =-\mu e^{-\mu t} \K_t h(x) + e^{-\mu t}\frac{d}{dt}\K_t h(x)\\
        & =-\mu e^{-\mu t} \K_t h(x) + e^{-\mu t}\K_t(\ll h)(x)\\
        & = - e^{-\mu t} \K_t (\ll h-\mu h)(x).
\end{split}
\end{equation}
For the second-order derivative, we have
\begin{equation}
    \begin{split}
        \frac{d^2}{dt^t}(e^{-\mu t}\K_th)(x)  =&\frac{d}{dt}(e^{-\mu t}\K_t(\ll h-\mu h)(x))\\
         = &-\mu e^{-\mu t}\K_t(\ll h -\mu h)(x)\\ &+e^{-\mu t}\K_t(\ll(\ll h-\mu h))(x)\\
         = & e^{-\mu t}\K_t(\ll^2 h-2\mu\ll h + \mu^2h)(x)\\
         = & e^{-\mu t}\K_t(\ll-\mu\id)^2(x)
    \end{split}
\end{equation}

By induction, the $k$-th derivative takes the form:
\begin{equation}
\frac{d^k}{dt^k}(e^{-\mu t}\K_th)(x) = e^{-\mu t}\K_t((\ll-\mu\id)^kh)(x),\;\;x\in\Omega,    
\end{equation}
which completes the proof. 
\end{pf}

\begin{thm}\label{thm_12}
    Let $h\in P^N(\Omega)$ and $k, N\in\mathbb{N}_0$. %Suppose $|\nabla^kh|_\infty\leq M_k^h$ for $k\in\{1, 2, \cdots, n\}$. 
    Let $L_f$ denote the (local) Lipschitz constant of $f$ on $\Omega$. Then, for all $t\geq 0$ and $x\in\Omega$,
    $$   \left| \frac{d^k}{dt^k}(e^{-\mu t}\K_th)(x)\right|   
    \leq    e^{(NL_f-\mu) t}(\mu+NL_f)^k|x|^N. $$
%where $C_N$ is such that $C_N^k\geq L_f^k\cdot \left(\prod_{j=0}^{k-1}(N+j)\right)$ for all $k\geq 1$. 
\end{thm}
\begin{pf}
Given the assumption that $f(\bf 0) = \bf 0$, we have $|f(x)| \leq L_f|x|$ for all $x\in\Omega$. % and $|f(x)|_\infty\leq L_f R=:C_f$. 

Then, $\ll h(x) = \sum_{i=1}^df_i(x)\partial_{x_i} h(x) = \sum_{i=1}^df_i(x)\alpha_ix^{\alpha-e_i}$, where $e_i$ is the $i$-th standard basis vector, and we have used the succinct expression $ x^{\alpha-e_i}:= x^{\alpha_i-1}_i\cdot\prod_{j\neq i}x_j^{\alpha_j}$. Therefore, 
\begin{equation}
    \begin{split}
        |\ll h(x)|  &\leq L_f|x|\cdot\sum_{i=1}^d\alpha_i|x^{\alpha-e_i}|\\ 
        & \leq   L_f|x|\cdot\sum_{i=1}^d\alpha_i \left[|x|^{\alpha_i-1}\cdot\prod_{j\neq i}|x|^{\alpha_j}\right]\\
        & =   L_f|x|\cdot\sum_{i=1}^d\alpha_i  |x|^{N-1}\\
        & = NL_f|x|^N, \;\forall x\in\Omega,
    \end{split}
\end{equation}
where we have assume that $\sum_{i=1}^d \alpha_i= N$. Consequently, due to the smoothness of $h$ and $h(\mathbf{0})=0$, we have $|\ll^2h(x)| \leq |f(x)|\cdot |\nabla \ll h(x)|  \leq L_f|x|\cdot (NL_f)N|x|^{N-1}  = (NL_f)^2|x|^N$ 
By induction, $|\ll^kh(x)|\leq (NL_f)^k|x|^N$, for all $x\in\X$.  
Using the formula \eqref{E: fomula_explicit}, we have 
\begin{equation*}
    \begin{split}
        \left| \frac{d^k}{dt^k}(e^{-\mu t}\K_th)(x)\right| & = \left| e^{-\mu t}\sum_{i=0}^k  \begin{pmatrix}
        k \\ i
    \end{pmatrix}(-\mu)^{k-i}|\K_t(\ll^i h)(x)|\right|\\
    & \leq e^{-\mu t}\sum_{i=0}^k  \begin{pmatrix}
        k \\ i
    \end{pmatrix}\mu^{k-i}|\ll^ih(\phi(t,x))|\\
    & \leq  e^{-\mu t}\sum_{i=0}^k  \begin{pmatrix}
        k \\ i
    \end{pmatrix}\mu^{k-i} (NL_f)^i|\phi(t,x)|^N.
    \end{split}
\end{equation*}
By Gr\"{o}wnwall's inequality, $|\phi(t,x)|\leq e^{L_ft}|x|$.  
Combining the above, we have 
\begin{equation}
\begin{split}
        & \left| \frac{d^k}{dt^k}(e^{-\mu t}\K_th)(x)\right|\\  \leq &e^{-\mu t}\sum_{i=0}^k  \begin{pmatrix}
        k \\ i
    \end{pmatrix}\mu^{k-i} (NL_f)^ie^{NL_ft}|x|^N\\
    \leq  & e^{(NL_f-\mu) t}|x|^N\sum_{i=0}^k  \begin{pmatrix}
        k \\ i
    \end{pmatrix}\mu^{k-i} (NL_f)^i.
\end{split}
\end{equation}
The conclusion follows immediately. 
\end{pf}

Denote $q(\mu, t)=e^{(NL_f-\mu) t}(\mu+NL_f)^k$. Then, one can easily verify that,  for $\mu\geq NL_f$,  $\sup_{t\in[0, T]}q(\mu, t)=(\mu+NL_f)^k$.   Combining the above estimation, the error bound of Gauss-Legendre quadrature for monomial observables at each $x\in\X$ for $\mu\geq NL_f$ is given by $|x|^NT^{2\Gamma + 1} \left(\frac{\mu+NL_f}{8\Gamma^2}\right)^{2\Gamma}$.

\section{Supplementary Numerical Experiments}
\subsection{Test of Parameters of the RTM}%\label{sec: van der pol}
\label{sec:vdp}
We use this example to thoroughly explain how   the parameters $\mu, \lambda$, $\taus$ and $\Gamma$ impact the precision of the RTM.  
 We revisit   the reversed Van der Pol oscillator example $\dot \xb_1 =-\xb_2$, $\dot \xb_2 = \xb_1 - (1 - \xb_1^2)\xb_2$ 
\iffalse\begin{equation*}
\dot \xb_1 =-\xb_2, \quad 
\dot \xb_2 = \xb_1 - (1 - \xb_1^2)\xb_2, 
\end{equation*} \fi
with $x:=[x_1, x_2]$.

As the impact of $\lambda$ and $\taus$ has been shown analytically in Theorem \ref{thm: conv} and \ref{thm: conv_t},  we can set $\lambda$ to be a very large number while fixing $\taus$ to any reasonable value. The main focus of this part of the experiment is to see how tuning $\mu$ can impact performance under the constraint of $\Gamma$.

\begin{figure}[hb!]
    \centering
    \includegraphics[width=0.9\linewidth]{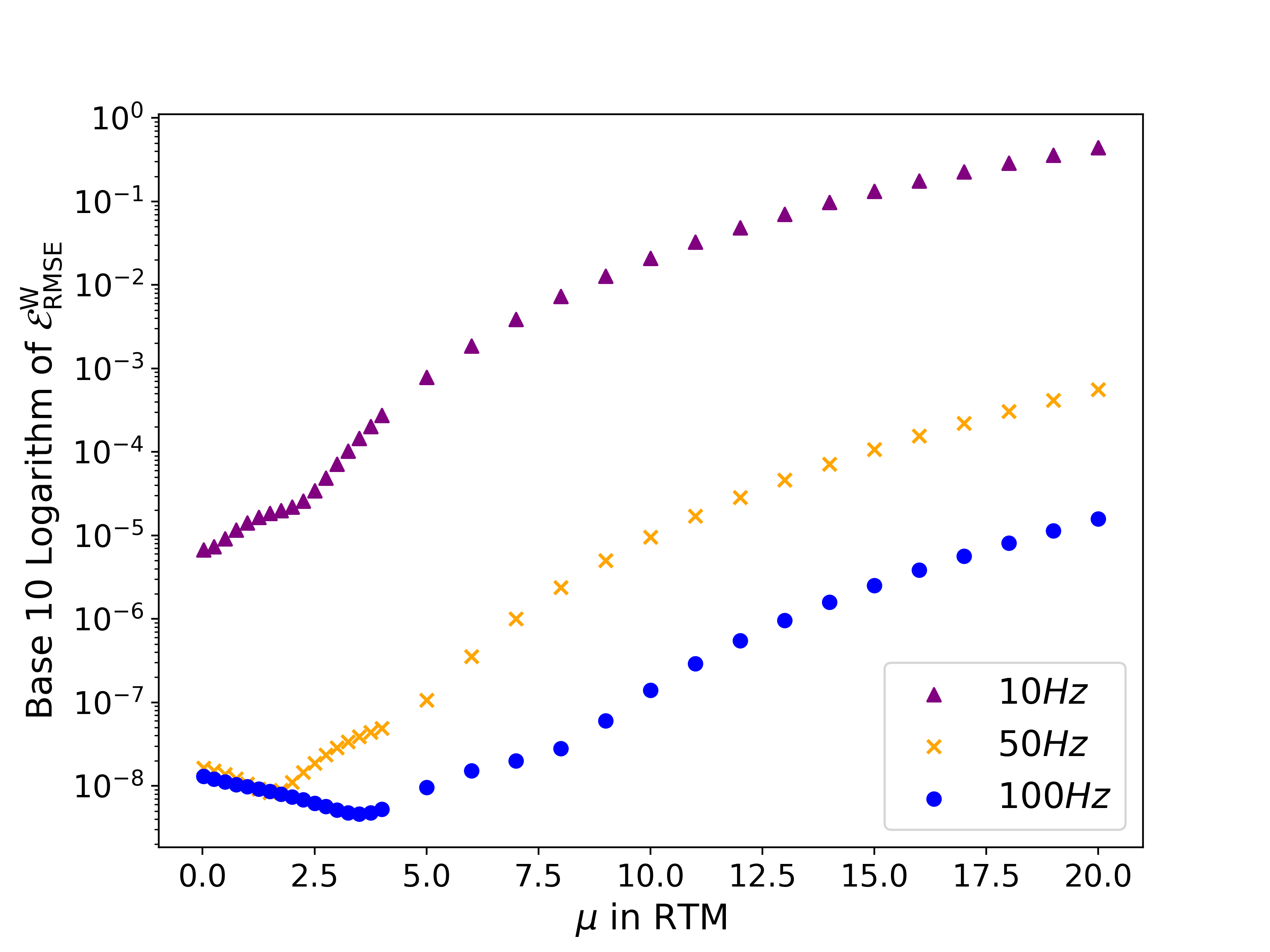}
    \caption{RMSE of weights ($\mathcal{E}_{\operatorname{RMSE}}^{\text{W}}$) using RTM for reversed Van der Pol oscillator with different sampling frequencies.}
    \label{fig:mu}
\end{figure}

In the experiments, we set $\lambda=1e8$ and   $\taus = 5$, with $N = 3\times 3$ monomials.  %$\{\zk_{i,j}\}$ as the dictionary functions where $i\in\set{0, 1, 2}$ and $j\in\set{0, 1, 2}$ \ym{(confirm this part with Yiming)}. 
We first use $\mu$ to evaluate the integral $\I_{\mu, T}^{\text{quad}}$ under the sampling frequency $\gamma = 100$   (equivalently, $\Gamma=500$), and then use  \eqref{E: approx_Lh} to learn the generator. Similarly, we perform the experiments for the cases of $\gamma = 50$ and $10$. In Fig.~\ref{fig:mu}, the relationship between $\taus$, $\mu$, and $\mathcal{E}_{\operatorname{RMSE}}^{\text{W}}$ %($\log_{10}$-scale) 
is illustrated.

Under a sampling rate of $100$, the highest accuracy $(\approx 10^{-8})$ is achieved at $\mu\approx 3.75$. A similar pattern is observed for the
cases of $\gamma = 50$   and $10$. However, the optimal accuracy decreases as the sampling frequency lowers, as expected. For a fixed tuple $(\taus, \lambda,  \Gamma)$, it is clear that the accuracy of the numerical integral $\I_{\mu, T}^{\text{quad}}$ increases as $\mu$ decreases, since the integrand $e^{-\mu s}\K_s \zk_{i} $ becomes less steep. %By observing the form of the inference formula in \eqref{E: mod_L} (or equivalently \eqref{E: mod_eq}), 
The non-monotonic trend of the learning accuracy with respect to changes in $\mu$ can be heuristically attributed to the combined effects of $\lambda$, $\mu$, and the %analytical and numerical 
errors of the  numerical integral.

\begin{rem}
    A similar pattern has also been found in other numerical experiments, but we omit presenting them to avoid repetitive expression of the same idea. The optimal tuning of parameters can be roughly determined by the local Lipschitz constant of the system. %However, since we have no knowledge of how fast the system evolves, estimating the optimal tuning may not be feasible at this point. 
    It is of the authors' interest to pursue further numerical analysis to obtain a reasonable estimate of the tuning under mild assumptions about the system in future work. 
    
    %this type of numerical analysis to obtain a reasonable estimate of the tuning under mild assumptions about the system in future work. 

   % Even with this missing part of the analysis for now, the method can still be applied in real-world scenarios. 
   The method has also been applied in real-world scenarios, such as identifying the car-following model  %Examples can be found in 
   \cite{meng2024koopmanitsc}, particularly in applications such as identifying the car-following model, where the acceleration/deceleration (which directly determines the Lipschitz constant of the model) is known to be bounded. Numerical results provide a range for the choice of $\mu$ that leads to acceptable accuracy under other fixed parameters.\Qed
\end{rem}

\begin{table*}[hbt!p]
\centering
\caption{Comparisons of RMSE of weights and flow over 100 trajectories for the system at the bifurcation point}
\begin{adjustbox}{max width=\textwidth, scale=1}
\begin{tabular}{ccc|c|cccc|cccc}
\hline
\hline
\multirow{2}{*}{System} & \multirow{2}{*}{$M$} & \multirow{2}{*}{$N$} & \multirow{2}{*}{\makecell{$\gamma$}} & \multicolumn{4}{c|}{RMSE of weights ($\mathcal{E}_{\operatorname{RMSE}}^{\text{W}}$) }    & \multicolumn{4}{c}{RMSE of flow ($\mathcal{E}_{\operatorname{RMSE}}^{\text{F}}$) }      
\\ \cline{5-12} & & &  & SINDy & FDM & KLM & RTM & SINDy & FDM & KLM & RTM \\ \hline
\multirow{3}{*}{1D Cubic}     & \multirow{3}{*}{$10$}   & \multirow{3}{*}{5}    &   10 & \textbf{3.98e-5} &  5.94e-2  &  1.61e-3  &  1.05e-4 &  5.18e-6 &  5.11e-3  & 1.23e-4  &  \textbf{3.18e-6}  \\ 
\cline{4-12}  &  & & 50 &  1.09e-7 &  1.36e-2  & 7.37e-5   &  \textbf{7.78e-8} & 1.41e-8 & 1.10e-3    &  5.70e-6 & \textbf{6.89e-9}   \\ 
\cline{4-12}  &  & & 100 & \textbf{7.31e-9} &  6.92e-3  & 1.90e-5   &  6.03e-8 & \textbf{9.43e-10} & 5.53e-4    &  1.52e-6 & 6.40e-9 \\ \hline
\hline
\hline
\end{tabular}
\end{adjustbox}
\label{tab:1D cubic}
\end{table*}

\subsection{System at the Bifurcation Point}
In this example,  
we consider the following scalar system that undergoes a pitchfork bifurcation $\dot{x} = \alpha x - x^{3}$. 
Apparently, the origin is an equilibrium point for all $\alpha$, which is stable for $\alpha < 0$ and unstable for $\alpha > 0$. When $\alpha > 0$, we have two extra stable equilibria branching from the origin i.e., $x_{1,2} = \pm\sqrt{\alpha}$. At the bifurcation point $\alpha = 0$, the origin is a non-hyperbolic equilibrium.  We then consider the system at the bifurcation point $\alpha = 0$:
\begin{equation}\label{E: cubic}
    \dot x = -x^3, \quad x(0) = x_0.
\end{equation}
Then, $x(t) = \frac{x_0}{\sqrt{1+2x_0^2t}}$. For all $h\in \mathcal{C}^1(\Omega)$, $\ll h(x) = -h'(x)\cdot x^3$. Specifically, $\ll h(0) = 0$.

 As pointed out in Remark \ref{rem: inversibility}, $\mathcal{C}^1(\X) \subsetneq \dom(\ll)$ and $\rho(\ll|_{\scriptscriptstyle\mathcal{C}^1(\Omega)})=\emptyset$.  %We provide an alternative way to show this fact.  Suppose $\lambda \in\rho(\ll|_{\scriptscriptstyle\mathcal{C}^1(\Omega)})$, then $\lambda\id-\ll|_{\scriptscriptstyle\mathcal{C}^1(\Omega)}$ is invertible and for any $h\in\mathcal{C}^1(\Omega)$, $\lambda u+x^3u'(x) = h(x)$ has a solution. Solving this, one has $u(x) = \exp\left(\frac{\lambda}{2x^2}\right)\left(\int_\Omega \frac{h(x)}{x^3}\exp\left(\frac{\lambda}{2x^2}\right)dx\right)$. However, it is only valid  on $\Omega\setminus\{0\}$, with a singular value at $0$.  
On the other hand, one can verify that   $\dom(\ll) = \{h\in \mathcal{C}(\Omega): h\in \mathcal{C}^1(\Omega\setminus\{0\}), \ll h\in\mathcal{C}(\Omega), \lim_{x\ra 0}\ll h(x) = 0\}$. Then $\lambda \id - \ll$ is injective on this domain, $\overline{\operatorname{Range}(\lambda\id-\ll)}=\mathcal{C}(\Omega)$ (in fact, $\operatorname{Range}(\lambda\id-\ll)=\mathcal{C}(\Omega)$, which implies surjectivity), and $\ker(\rr(\lambda))=\{0\}$. The closedness of $\ll$ also holds on this domain Therefore, $\rr(\lambda) = (\lambda\id-\ll)^{-1}$, meaning that the pseudo equals the true resolvent operator.

 For the purpose of system identification, we generally do not know and therefore do not require knowledge of the singular point. Consequently, we may not find a function space in which $\ll$ has the desired properties as stated above. Nonetheless,  we always have $\lambda^2\rr(\lambda)-\lambda\id \sra \ll$ on regular function space.

 We then use monomials up to order 4 as the dictionary functions, i.e., $\{1, x, x^2, x^3, x^4\}$, and $M=10$ initial points randomly sampled within $\X=(-1, 1)$. We set  $\tau_s = 1$, $\mu=0.02$ and  $\lambda=1e8$ for RTM. The comparisons with SINDy, WINDy, FDM, and KLM are presented in Table~\ref{tab:1D cubic}. It is evident that RTM outperforms the other two Koopman-based methods and reaches comparable performance as SINDy.

\subsection{Discussion of Other Koopman-Based Methods for Finding the RoA}
Most of existing Koopman-based method for learning Lyapunov functions and the corresponding RoA estimates are based on the spectrum of the Koopman operator, i.e., constructing the Lyapunov functions from the eigenfunctions of the Koopman operator \cite{mauroy2016global, deka2022koopman, zheng2022data}. However, there are two main weaknesses with the method that directly constructs Lyapunov functions from the eigenfunctions. First, the eigenfunctions themselves are only approximations due to the finite-dimensional approximation of the Koopman operators within a limited set of functions. This can lead to errors and no formal guarantees of correctness. Second, even when many valid eigenfunctions exist, one must manually choose which one (or combination) gives the  ``best'' Lyapunov function—usually the one that yields the largest RoA. Since any eigenfunction with a negative real part of its eigenvalue technically works, and weighted sums or products of them can also qualify, finding the optimal candidate becomes cumbersome and uncertain.

In contrast, unlike the Koopman eigenfunction-based approaches, the proposed method avoids guessing or tuning over combinations of the Koopman eigenfunctions for learning the (maximal) Lyapunov functions. It directly solves the Zubov's equation, ensuring the Lyapunov conditions are satisfied by construction. This approach is more general, as it does not rely on specially chosen observable functions or handcrafted dictionaries, which are sometimes needed for eigenfunction methods. As a result, it can handle a broader range of nonlinear systems, including those with unknown dynamics, while still achieving comparable or better performance without heavy customization or manual parameter tuning.

%\section*{References}
\bibliographystyle{ieeetr}        % Include this if you use bibtex 
\bibliography{TAC}

\end{document}